
\ifx\shlhetal\undefinedcontrolsequence\let\shlhetal\relax\fi
\def\fmtname{AmS-TeX}

\def\fmtversion{2.1}
\catcode`\@=11
\ifx\amstexloaded@\relax\catcode`\@=\active
  \endinput\else\let\amstexloaded@\relax\fi
\newlinechar=`\^^J
\def\W@{\immediate\write\sixt@@n}
\def\CR@{\W@{^^J\fmtname - Version \fmtversion^^J%
COPYRIGHT 1985, 1990, 1991 - AMERICAN MATHEMATICAL SOCIETY^^J%
Use of this macro package is not restricted provided^^J%
each use is acknowledged upon publication.^^J}}
\CR@ \everyjob{\CR@}
\message{Loading definitions for}
\message{misc utility macros,}
\toksdef\toks@@=2
\long\def\rightappend@#1\to#2{\toks@{\\{#1}}\toks@@
 =\expandafter{#2}\xdef#2{\the\toks@@\the\toks@}\toks@{}\toks@@{}}
\def\alloclist@{}
\newif\ifalloc@
\def\showallocations{{\def\\{\immediate\write\m@ne}\alloclist@}\alloc@true}
\def\alloc@#1#2#3#4#5{\global\advance\count1#1by\@ne
 \ch@ck#1#4#2\allocationnumber=\count1#1
 \global#3#5=\allocationnumber
 \edef\next@{\string#5=\string#2\the\allocationnumber}%
 \expandafter\rightappend@\next@\to\alloclist@}
\newcount\count@@
\newcount\count@@@
\def\FN@{\futurelet\next}
\def\DN@{\def\next@}
\def\DNii@{\def\nextii@}
\def\RIfM@{\relax\ifmmode}
\def\RIfMIfI@{\relax\ifmmode\ifinner}
\def\setboxz@h{\setbox\z@\hbox}
\def\wdz@{\wd\z@}
\def\boxz@{\box\z@}
\def\setbox@ne{\setbox\@ne}
\def\wd@ne{\wd\@ne}
\def\iterate{\body\expandafter\iterate\else\fi}
\def\err@#1{\errmessage{AmS-TeX error: #1}}
\newhelp\defaulthelp@{Sorry, I already gave what help I could...^^J
Maybe you should try asking a human?^^J
An error might have occurred before I noticed any problems.^^J
``If all else fails, read the instructions.''}
\def\Err@{\errhelp\defaulthelp@\err@}
\def\eat@#1{}
\def\in@#1#2{\def\in@@##1#1##2##3\in@@{\ifx\in@##2\in@false\else\in@true\fi}%
 \in@@#2#1\in@\in@@}
\newif\ifin@
\def\space@.{\futurelet\space@\relax}
\space@. %
\newhelp\athelp@
{Only certain combinations beginning with @ make sense to me.^^J
Perhaps you wanted \string\@\space for a printed @?^^J
I've ignored the character or group after @.}
{\catcode`\~=\active 
 \lccode`\~=`\@ \lowercase{\gdef~{\FN@\at@}}}
\def\at@{\let\next@\at@@
 \ifcat\noexpand\next a\else\ifcat\noexpand\next0\else
 \ifcat\noexpand\next\relax\else
   \let\next\at@@@\fi\fi\fi
 \next@}
\def\at@@#1{\expandafter
 \ifx\csname\space @\string#1\endcsname\relax
  \expandafter\at@@@ \else
  \csname\space @\string#1\expandafter\endcsname\fi}
\def\at@@@#1{\errhelp\athelp@ \err@{\Invalid@@ @}}
\def\atdef@#1{\expandafter\def\csname\space @\string#1\endcsname}
\newhelp\defahelp@{If you typed \string\define\space cs instead of
\string\define\string\cs\space^^J
I've substituted an inaccessible control sequence so that your^^J
definition will be completed without mixing me up too badly.^^J
If you typed \string\define{\string\cs} the inaccessible control sequence^^J
was defined to be \string\cs, and the rest of your^^J
definition appears as input.}
\newhelp\defbhelp@{I've ignored your definition, because it might^^J
conflict with other uses that are important to me.}
\def\define{\FN@\define@}
\def\define@{\ifcat\noexpand\next\relax
 \expandafter\define@@\else\errhelp\defahelp@                               
 \err@{\string\define\space must be followed by a control
 sequence}\expandafter\def\expandafter\nextii@\fi}                          
\def\undefined@@@@@@@@@@{}
\def\preloaded@@@@@@@@@@{}
\def\next@@@@@@@@@@{}
\def\define@@#1{\ifx#1\relax\errhelp\defbhelp@                              
 \err@{\string#1\space is already defined}\DN@{\DNii@}\else
 \expandafter\ifx\csname\expandafter\eat@\string                            
 #1@@@@@@@@@@\endcsname\undefined@@@@@@@@@@\errhelp\defbhelp@
 \err@{\string#1\space can't be defined}\DN@{\DNii@}\else
 \expandafter\ifx\csname\expandafter\eat@\string#1\endcsname\relax          
 \global\let#1\undefined\DN@{\def#1}\else\errhelp\defbhelp@
 \err@{\string#1\space is already defined}\DN@{\DNii@}\fi
 \fi\fi\next@}

\def\predefine#1#2{\let#1#2}
\def\undefine#1{\let#1\undefined}
\message{page layout,}
\newdimen\captionwidth@
\captionwidth@\hsize
\advance\captionwidth@-1.5in
\def\pagewidth#1{\hsize#1\relax
 \captionwidth@\hsize\advance\captionwidth@-1.5in}
\def\pageheight#1{\vsize#1\relax}
\def\hcorrection#1{\advance\hoffset#1\relax}
\def\vcorrection#1{\advance\voffset#1\relax}
\message{accents/punctuation,}

\let\graveaccent\`
\let\acuteaccent\'
\let\tildeaccent\~
\let\hataccent\^
\let\underscore\_
\let\B\=
\let\D\.
\let\ic@\/
\def\/{\unskip\ic@}
\def\textfonti{\the\textfont\@ne}
\def\t#1#2{{\edef\next@{\the\font}\textfonti\accent"7F \next@#1#2}}
\def~{\unskip\nobreak\ \ignorespaces}
\def\.{.\spacefactor\@m}
\atdef@;{\leavevmode\null;}
\atdef@:{\leavevmode\null:}
\atdef@?{\leavevmode\null?}
\edef\@{\string @}
\def\relaxnext@{\let\next\relax}
\atdef@-{\relaxnext@\leavevmode
 \DN@{\ifx\next-\DN@-{\FN@\nextii@}\else
  \DN@{\leavevmode\hbox{-}}\fi\next@}%
 \DNii@{\ifx\next-\DN@-{\leavevmode\hbox{---}}\else
  \DN@{\leavevmode\hbox{--}}\fi\next@}%
 \FN@\next@}
\def\srdr@{\kern.16667em}
\def\drsr@{\kern.02778em}
\def\sldl@{\drsr@}
\def\dlsl@{\srdr@}
\atdef@"{\unskip\relaxnext@
 \DN@{\ifx\next\space@\DN@. {\FN@\nextii@}\else
  \DN@.{\FN@\nextii@}\fi\next@.}%
 \DNii@{\ifx\next`\DN@`{\FN@\nextiii@}\else
  \ifx\next\lq\DN@\lq{\FN@\nextiii@}\else
  \DN@####1{\FN@\nextiv@}\fi\fi\next@}%
 \def\nextiii@{\ifx\next`\DN@`{\sldl@``}\else\ifx\next\lq
  \DN@\lq{\sldl@``}\else\DN@{\dlsl@`}\fi\fi\next@}%
 \def\nextiv@{\ifx\next'\DN@'{\srdr@''}\else
  \ifx\next\rq\DN@\rq{\srdr@''}\else\DN@{\drsr@'}\fi\fi\next@}%
 \FN@\next@}

\def\textfontii{\the\textfont\tw@}
\def\lbrace@{\delimiter"4266308 }
\def\rbrace@{\delimiter"5267309 }
\def\{{\RIfM@\lbrace@\else{\textfontii f}\spacefactor\@m\fi}
\def\}{\RIfM@\rbrace@\else
 \let\@sf\empty\ifhmode\edef\@sf{\spacefactor\the\spacefactor}\fi
 {\textfontii g}\@sf\relax\fi}
\let\lbrace\{
\let\rbrace\}
\def\AmSTeX{{\textfontii A\kern-.1667em%
  \lower.5ex\hbox{M}\kern-.125emS}-\TeX}
\message{line and page breaks,}
\def\vmodeerr@#1{\Err@{\string#1\space not allowed between paragraphs}}
\def\mathmodeerr@#1{\Err@{\string#1\space not allowed in math mode}}
\def\linebreak{\RIfM@\mathmodeerr@\linebreak\else
 \ifhmode\unskip\unkern\break\else\vmodeerr@\linebreak\fi\fi}

\newskip\saveskip@
\def\allowlinebreak{\RIfM@\mathmodeerr@\allowlinebreak\else
 \ifhmode\saveskip@\lastskip\unskip
 \allowbreak\ifdim\saveskip@>\z@\hskip\saveskip@\fi
 \else\vmodeerr@\allowlinebreak\fi\fi}
\def\nolinebreak{\RIfM@\mathmodeerr@\nolinebreak\else
 \ifhmode\saveskip@\lastskip\unskip
 \nobreak\ifdim\saveskip@>\z@\hskip\saveskip@\fi
 \else\vmodeerr@\nolinebreak\fi\fi}
\def\newline{\relaxnext@
 \DN@{\RIfM@\expandafter\mathmodeerr@\expandafter\newline\else
  \ifhmode\ifx\next\par\else
  \expandafter\unskip\expandafter\null\expandafter\hfill\expandafter\break\fi
  \else
  \expandafter\vmodeerr@\expandafter\newline\fi\fi}%
 \FN@\next@}
\def\dmatherr@#1{\Err@{\string#1\space not allowed in display math mode}}
\def\nondmatherr@#1{\Err@{\string#1\space not allowed in non-display math
 mode}}
\def\onlydmatherr@#1{\Err@{\string#1\space allowed only in display math mode}}
\def\nonmatherr@#1{\Err@{\string#1\space allowed only in math mode}}
\def\mathbreak{\RIfMIfI@\break\else
 \dmatherr@\mathbreak\fi\else\nonmatherr@\mathbreak\fi}
\def\nomathbreak{\RIfMIfI@\nobreak\else
 \dmatherr@\nomathbreak\fi\else\nonmatherr@\nomathbreak\fi}
\def\allowmathbreak{\RIfMIfI@\allowbreak\else
 \dmatherr@\allowmathbreak\fi\else\nonmatherr@\allowmathbreak\fi}
\def\pagebreak{\RIfM@
 \ifinner\nondmatherr@\pagebreak\else\postdisplaypenalty-\@M\fi
 \else\ifvmode\removelastskip\break\else\vadjust{\break}\fi\fi}
\def\nopagebreak{\RIfM@
 \ifinner\nondmatherr@\nopagebreak\else\postdisplaypenalty\@M\fi
 \else\ifvmode\nobreak\else\vadjust{\nobreak}\fi\fi}
\def\nonvmodeerr@#1{\Err@{\string#1\space not allowed within a paragraph
 or in math}}
\def\vnonvmode@#1#2{\relaxnext@\DNii@{\ifx\next\par\DN@{#1}\else
 \DN@{#2}\fi\next@}%
 \ifvmode\DN@{#1}\else
 \DN@{\FN@\nextii@}\fi\next@}
\def\newpage{\vnonvmode@{\vfill\break}{\nonvmodeerr@\newpage}}
\def\smallpagebreak{\vnonvmode@\smallbreak{\nonvmodeerr@\smallpagebreak}}
\def\medpagebreak{\vnonvmode@\medbreak{\nonvmodeerr@\medpagebreak}}
\def\bigpagebreak{\vnonvmode@\bigbreak{\nonvmodeerr@\bigpagebreak}}
\def\NoBlackBoxes{\global\overfullrule\z@}
\def\BlackBoxes{\global\overfullrule5\p@}
\def\Invalid@#1{\def#1{\Err@{\Invalid@@\string#1}}}
\def\Invalid@@{Invalid use of }
\message{figures,}
\Invalid@\caption
\Invalid@\captionwidth
\newdimen\smallcaptionwidth@
\def\topspace{\mid@false\ins@}
\def\midspace{\mid@true\ins@}
\newif\ifmid@
\def\captionfont@{}
\def\ins@#1{\relaxnext@\allowbreak
 \smallcaptionwidth@\captionwidth@\gdef\thespace@{#1}%
 \DN@{\ifx\next\space@\DN@. {\FN@\nextii@}\else
  \DN@.{\FN@\nextii@}\fi\next@.}%
 \DNii@{\ifx\next\caption\DN@\caption{\FN@\nextiii@}%
  \else\let\next@\nextiv@\fi\next@}%
 \def\nextiv@{\vnonvmode@
  {\ifmid@\expandafter\midinsert\else\expandafter\topinsert\fi
   \vbox to\thespace@{}\endinsert}
  {\ifmid@\nonvmodeerr@\midspace\else\nonvmodeerr@\topspace\fi}}%
 \def\nextiii@{\ifx\next\captionwidth\expandafter\nextv@
  \else\expandafter\nextvi@\fi}%
 \def\nextv@\captionwidth##1##2{\smallcaptionwidth@##1\relax\nextvi@{##2}}%
 \def\nextvi@##1{\def\thecaption@{\captionfont@##1}%
  \DN@{\ifx\next\space@\DN@. {\FN@\nextvii@}\else
   \DN@.{\FN@\nextvii@}\fi\next@.}%
  \FN@\next@}%
 \def\nextvii@{\vnonvmode@
  {\ifmid@\expandafter\midinsert\else
  \expandafter\topinsert\fi\vbox to\thespace@{}\nobreak\smallskip
  \setboxz@h{\noindent\ignorespaces\thecaption@\unskip}%
  \ifdim\wdz@>\smallcaptionwidth@\centerline{\vbox{\hsize\smallcaptionwidth@
   \noindent\ignorespaces\thecaption@\unskip}}%
  \else\centerline{\boxz@}\fi\endinsert}
  {\ifmid@\nonvmodeerr@\midspace
  \else\nonvmodeerr@\topspace\fi}}%
 \FN@\next@}
\message{comments,}
\def\newcodes@{\catcode`\\12\catcode`\{12\catcode`\}12\catcode`\#12%
 \catcode`\%12\relax}
\def\oldcodes@{\catcode`\\0\catcode`\{1\catcode`\}2\catcode`\#6%
 \catcode`\%14\relax}
\def\comment{\newcodes@\endlinechar=10 \comment@}
{\lccode`\0=`\\
\lowercase{\gdef\comment@#1^^J{\comment@@#10endcomment\comment@@@}%
\gdef\comment@@#10endcomment{\FN@\comment@@@}%
\gdef\comment@@@#1\comment@@@{\ifx\next\comment@@@\let\next\comment@
 \else\def\next{\oldcodes@\endlinechar=`\^^M\relax}%
 \fi\next}}}
\def\pr@m@s{\ifx'\next\DN@##1{\prim@s}\else\let\next@\egroup\fi\next@}
\def\prime{{\null\prime@\null}}
\mathchardef\prime@="0230
\let\dsize\displaystyle

\let\ssize\scriptstyle

\message{math spacing,}
\def\,{\RIfM@\mskip\thinmuskip\relax\else\kern.16667em\fi}
\def\!{\RIfM@\mskip-\thinmuskip\relax\else\kern-.16667em\fi}
\let\thinspace\,
\let\negthinspace\!
\def\medspace{\RIfM@\mskip\medmuskip\relax\else\kern.222222em\fi}
\def\negmedspace{\RIfM@\mskip-\medmuskip\relax\else\kern-.222222em\fi}
\def\thickspace{\RIfM@\mskip\thickmuskip\relax\else\kern.27777em\fi}
\let\;\thickspace
\def\negthickspace{\RIfM@\mskip-\thickmuskip\relax\else
 \kern-.27777em\fi}
\atdef@,{\RIfM@\mskip.1\thinmuskip\else\leavevmode\null,\fi}
\atdef@!{\RIfM@\mskip-.1\thinmuskip\else\leavevmode\null!\fi}
\atdef@.{\RIfM@&&\else\leavevmode.\spacefactor3000 \fi}
\def\and{\DOTSB\;\mathchar"3026 \;}

\message{fractions,}
\def\frac#1#2{{#1\over#2}}

\newdimen\ex@
\ex@.2326ex
\Invalid@\thickness
\def\thickfrac{\relaxnext@
 \DN@{\ifx\next\thickness\let\next@\nextii@\else
 \DN@{\nextii@\thickness1}\fi\next@}%
 \DNii@\thickness##1##2##3{{##2\above##1\ex@##3}}%
 \FN@\next@}

\def\thickfracwithdelims#1#2{\relaxnext@\def\ldelim@{#1}\def\rdelim@{#2}%
 \DN@{\ifx\next\thickness\let\next@\nextii@\else
 \DN@{\nextii@\thickness1}\fi\next@}%
 \DNii@\thickness##1##2##3{{##2\abovewithdelims
 \ldelim@\rdelim@##1\ex@##3}}%
 \FN@\next@}

\def\:{\nobreak\hskip.1111em\mathpunct{}\nonscript\mkern-\thinmuskip{:}\hskip
 .3333emplus.0555em\relax}
\def\snug{\unskip\kern-\mathsurround}
\message{smash commands,}
\def\topsmash{\top@true\bot@false\smash@}
\def\botsmash{\top@false\bot@true\smash@}
\newif\iftop@
\newif\ifbot@
\def\smash{\top@true\bot@true\smash@}
\def\smash@{\RIfM@\expandafter\mathpalette\expandafter\mathsm@sh\else
 \expandafter\makesm@sh\fi}
\def\finsm@sh{\iftop@\ht\z@\z@\fi\ifbot@\dp\z@\z@\fi\leavevmode\boxz@}
\message{large operator symbols,}
\def\LimitsOnSums{\global\let\slimits@\displaylimits}
\def\NoLimitsOnSums{\global\let\slimits@\nolimits}
\LimitsOnSums
\mathchardef\coprod@="1360       \def\coprod{\DOTSB\coprod@\slimits@}
\mathchardef\bigvee@="1357       \def\bigvee{\DOTSB\bigvee@\slimits@}
\mathchardef\bigwedge@="1356     \def\bigwedge{\DOTSB\bigwedge@\slimits@}
\mathchardef\biguplus@="1355     \def\biguplus{\DOTSB\biguplus@\slimits@}
\mathchardef\bigcap@="1354       \def\bigcap{\DOTSB\bigcap@\slimits@}
\mathchardef\bigcup@="1353       \def\bigcup{\DOTSB\bigcup@\slimits@}
\mathchardef\prod@="1351         \def\prod{\DOTSB\prod@\slimits@}
\mathchardef\sum@="1350          \def\sum{\DOTSB\sum@\slimits@}
\mathchardef\bigotimes@="134E    \def\bigotimes{\DOTSB\bigotimes@\slimits@}
\mathchardef\bigoplus@="134C     \def\bigoplus{\DOTSB\bigoplus@\slimits@}
\mathchardef\bigodot@="134A      \def\bigodot{\DOTSB\bigodot@\slimits@}
\mathchardef\bigsqcup@="1346     \def\bigsqcup{\DOTSB\bigsqcup@\slimits@}
\message{integrals,}
\def\LimitsOnInts{\global\let\ilimits@\displaylimits}
\def\NoLimitsOnInts{\global\let\ilimits@\nolimits}
\NoLimitsOnInts
\def\int{\DOTSI\intop\ilimits@}
\def\oint{\DOTSI\ointop\ilimits@}
\def\intic@{\mathchoice{\hskip.5em}{\hskip.4em}{\hskip.4em}{\hskip.4em}}
\def\negintic@{\mathchoice
 {\hskip-.5em}{\hskip-.4em}{\hskip-.4em}{\hskip-.4em}}
\def\intkern@{\mathchoice{\!\!\!}{\!\!}{\!\!}{\!\!}}
\def\intdots@{\mathchoice{\plaincdots@}
 {{\cdotp}\mkern1.5mu{\cdotp}\mkern1.5mu{\cdotp}}
 {{\cdotp}\mkern1mu{\cdotp}\mkern1mu{\cdotp}}
 {{\cdotp}\mkern1mu{\cdotp}\mkern1mu{\cdotp}}}
\newcount\intno@
\def\iint{\DOTSI\intno@\tw@\FN@\ints@}
\def\iiint{\DOTSI\intno@\thr@@\FN@\ints@}
\def\iiiint{\DOTSI\intno@4 \FN@\ints@}
\def\idotsint{\DOTSI\intno@\z@\FN@\ints@}
\def\ints@{\findlimits@\ints@@}
\newif\iflimtoken@
\newif\iflimits@
\def\findlimits@{\limtoken@true\ifx\next\limits\limits@true
 \else\ifx\next\nolimits\limits@false\else
 \limtoken@false\ifx\ilimits@\nolimits\limits@false\else
 \ifinner\limits@false\else\limits@true\fi\fi\fi\fi}
\def\multint@{\int\ifnum\intno@=\z@\intdots@                                
 \else\intkern@\fi                                                          
 \ifnum\intno@>\tw@\int\intkern@\fi                                         
 \ifnum\intno@>\thr@@\int\intkern@\fi                                       
 \int}                                                                      
\def\multintlimits@{\intop\ifnum\intno@=\z@\intdots@\else\intkern@\fi
 \ifnum\intno@>\tw@\intop\intkern@\fi
 \ifnum\intno@>\thr@@\intop\intkern@\fi\intop}
\def\ints@@{\iflimtoken@                                                    
 \def\ints@@@{\iflimits@\negintic@\mathop{\intic@\multintlimits@}\limits    
  \else\multint@\nolimits\fi                                                
  \eat@}                                                                    
 \else                                                                      
 \def\ints@@@{\iflimits@\negintic@
  \mathop{\intic@\multintlimits@}\limits\else
  \multint@\nolimits\fi}\fi\ints@@@}
\def\LimitsOnNames{\global\let\nlimits@\displaylimits}
\def\NoLimitsOnNames{\global\let\nlimits@\nolimits@}
\LimitsOnNames
\def\nolimits@{\relaxnext@
 \DN@{\ifx\next\limits\DN@\limits{\nolimits}\else
  \let\next@\nolimits\fi\next@}%
 \FN@\next@}
\message{operator names,}
\def\newmcodes@{\mathcode`\'"27\mathcode`\*"2A\mathcode`\."613A%
 \mathcode`\-"2D\mathcode`\/"2F\mathcode`\:"603A }
\def\operatorname#1{\mathop{\newmcodes@\kern\z@\fam\z@#1}\nolimits@}
\def\operatornamewithlimits#1{\mathop{\newmcodes@\kern\z@\fam\z@#1}\nlimits@}
\def\qopname@#1{\mathop{\fam\z@#1}\nolimits@}
\def\qopnamewl@#1{\mathop{\fam\z@#1}\nlimits@}
\def\arccos{\qopname@{arccos}}
\def\arcsin{\qopname@{arcsin}}
\def\arctan{\qopname@{arctan}}
\def\arg{\qopname@{arg}}
\def\cos{\qopname@{cos}}
\def\cosh{\qopname@{cosh}}
\def\cot{\qopname@{cot}}
\def\coth{\qopname@{coth}}
\def\csc{\qopname@{csc}}
\def\deg{\qopname@{deg}}
\def\det{\qopnamewl@{det}}
\def\dim{\qopname@{dim}}
\def\exp{\qopname@{exp}}
\def\gcd{\qopnamewl@{gcd}}
\def\hom{\qopname@{hom}}
\def\inf{\qopnamewl@{inf}}
\def\injlim{\qopnamewl@{inj\,lim}}
\def\ker{\qopname@{ker}}
\def\lg{\qopname@{lg}}
\def\lim{\qopnamewl@{lim}}
\def\liminf{\qopnamewl@{lim\,inf}}
\def\limsup{\qopnamewl@{lim\,sup}}
\def\ln{\qopname@{ln}}
\def\log{\qopname@{log}}
\def\max{\qopnamewl@{max}}
\def\min{\qopnamewl@{min}}
\def\Pr{\qopnamewl@{Pr}}
\def\projlim{\qopnamewl@{proj\,lim}}
\def\sec{\qopname@{sec}}
\def\sin{\qopname@{sin}}
\def\sinh{\qopname@{sinh}}
\def\sup{\qopnamewl@{sup}}
\def\tan{\qopname@{tan}}
\def\tanh{\qopname@{tanh}}
\def\varinjlim{\mathop{\vtop{\ialign{##\crcr
 \hfil\rm lim\hfil\crcr\noalign{\nointerlineskip}\rightarrowfill\crcr
 \noalign{\nointerlineskip\kern-\ex@}\crcr}}}}
\def\varprojlim{\mathop{\vtop{\ialign{##\crcr
 \hfil\rm lim\hfil\crcr\noalign{\nointerlineskip}\leftarrowfill\crcr
 \noalign{\nointerlineskip\kern-\ex@}\crcr}}}}
\def\varliminf{\mathop{\underline{\vrule height\z@ depth.2exwidth\z@
 \hbox{\rm lim}}}}

\newdimen\buffer@
\buffer@\fontdimen13 \tenex
\newdimen\buffer
\buffer\buffer@

\def\ResetBuffer{\fontdimen13 \tenex\buffer@\global\buffer\buffer@}
\def\shave#1{\mathop{\hbox{$\m@th\fontdimen13 \tenex\z@                     
 \displaystyle{#1}$}}\fontdimen13 \tenex\buffer}

\message{multilevel sub/superscripts,}
\Invalid@\\
\def\Let@{\relax\iffalse{\fi\let\\=\cr\iffalse}\fi}
\Invalid@\vspace
\def\vspace@{\def\vspace##1{\crcr\noalign{\vskip##1\relax}}}
\def\multilimits@{\bgroup\vspace@\Let@
 \baselineskip\fontdimen10 \scriptfont\tw@
 \advance\baselineskip\fontdimen12 \scriptfont\tw@
 \lineskip\thr@@\fontdimen8 \scriptfont\thr@@
 \lineskiplimit\lineskip
 \vbox\bgroup\ialign\bgroup\hfil$\m@th\scriptstyle{##}$\hfil\crcr}
\def\Sb{_\multilimits@}
\def\endSb{\crcr\egroup\egroup\egroup}
\def\Sp{^\multilimits@}

\def\spreadlines#1{\RIfMIfI@\onlydmatherr@\spreadlines\else
 \openup#1\relax\fi\else\onlydmatherr@\spreadlines\fi}
\def\Mathstrut@{\copy\Mathstrutbox@}
\newbox\Mathstrutbox@
\setbox\Mathstrutbox@\null
\setboxz@h{$\m@th($}
\ht\Mathstrutbox@\ht\z@
\dp\Mathstrutbox@\dp\z@
\message{matrices,}
\newdimen\spreadmlines@
\def\spreadmatrixlines#1{\RIfMIfI@
 \onlydmatherr@\spreadmatrixlines\else
 \spreadmlines@#1\relax\fi\else\onlydmatherr@\spreadmatrixlines\fi}
\def\matrix{\null\,\vcenter\bgroup\Let@\vspace@
 \normalbaselines\openup\spreadmlines@\ialign
 \bgroup\hfil$\m@th##$\hfil&&\quad\hfil$\m@th##$\hfil\crcr
 \Mathstrut@\crcr\noalign{\kern-\baselineskip}}
\def\endmatrix{\crcr\Mathstrut@\crcr\noalign{\kern-\baselineskip}\egroup
 \egroup\,}
\def\format{\crcr\egroup\iffalse{\fi\ifnum`}=0 \fi\format@}
\newtoks\hashtoks@
\hashtoks@{#}
\def\format@#1\\{\def\preamble@{#1}%
 \def\l{$\m@th\the\hashtoks@$\hfil}%
 \def\c{\hfil$\m@th\the\hashtoks@$\hfil}%
 \def\r{\hfil$\m@th\the\hashtoks@$}%
 \edef\preamble@@{\preamble@}\ifnum`{=0 \fi\iffalse}\fi
 \ialign\bgroup\span\preamble@@\crcr}
\def\smallmatrix{\null\,\vcenter\bgroup\vspace@\Let@
 \baselineskip9\ex@\lineskip\ex@
 \ialign\bgroup\hfil$\m@th\scriptstyle{##}$\hfil&&\thickspace\hfil
 $\m@th\scriptstyle{##}$\hfil\crcr}
\def\endsmallmatrix{\crcr\egroup\egroup\,}

\newmuskip\dotsspace@
\dotsspace@1.5mu
\def\strip@#1 {#1}
\def\spacehdots#1\for#2{\multispan{#2}\xleaders
 \hbox{$\m@th\mkern\strip@#1 \dotsspace@.\mkern\strip@#1 \dotsspace@$}\hfill}
\def\hdotsfor#1{\spacehdots\@ne\for{#1}}
\def\multispan@#1{\omit\mscount#1\unskip\loop\ifnum\mscount>\@ne\sp@n\repeat}
\def\spaceinnerhdots#1\for#2\after#3{\multispan@{\strip@#2 }#3\xleaders
 \hbox{$\m@th\mkern\strip@#1 \dotsspace@.\mkern\strip@#1 \dotsspace@$}\hfill}
\def\innerhdotsfor#1\after#2{\spaceinnerhdots\@ne\for#1\after{#2}}
\def\cases{\bgroup\spreadmlines@\jot\left\{\,\matrix\format\l&\quad\l\\}
\def\endcases{\endmatrix\right.\egroup}
\message{multiline displays,}
\newif\ifinany@
\newif\ifinalign@
\newif\ifingather@
\def\strut@{\copy\strutbox@}
\newbox\strutbox@
\setbox\strutbox@\hbox{\vrule height8\p@ depth3\p@ width\z@}
\def\topaligned{\null\,\vtop\aligned@}
\def\botaligned{\null\,\vbox\aligned@}
\def\aligned{\null\,\vcenter\aligned@}
\def\aligned@{\bgroup\vspace@\Let@
 \ifinany@\else\openup\jot\fi\ialign
 \bgroup\hfil\strut@$\m@th\displaystyle{##}$&
 $\m@th\displaystyle{{}##}$\hfil\crcr}
\def\endaligned{\crcr\egroup\egroup}

\def\alignedat#1{\null\,\vcenter\bgroup\doat@{#1}\vspace@\Let@
 \ifinany@\else\openup\jot\fi\ialign\bgroup\span\preamble@@\crcr}
\newcount\atcount@
\def\doat@#1{\toks@{\hfil\strut@$\m@th
 \displaystyle{\the\hashtoks@}$&$\m@th\displaystyle
 {{}\the\hashtoks@}$\hfil}
 \atcount@#1\relax\advance\atcount@\m@ne                                    
 \loop\ifnum\atcount@>\z@\toks@=\expandafter{\the\toks@&\hfil$\m@th
 \displaystyle{\the\hashtoks@}$&$\m@th
 \displaystyle{{}\the\hashtoks@}$\hfil}\advance
  \atcount@\m@ne\repeat                                                     
 \xdef\preamble@{\the\toks@}\xdef\preamble@@{\preamble@}}

\def\gathered{\null\,\vcenter\bgroup\vspace@\Let@
 \ifinany@\else\openup\jot\fi\ialign
 \bgroup\hfil\strut@$\m@th\displaystyle{##}$\hfil\crcr}
\def\endgathered{\crcr\egroup\egroup}
\newif\iftagsleft@
\def\TagsOnLeft{\global\tagsleft@true}
\def\TagsOnRight{\global\tagsleft@false}
\TagsOnLeft
\newif\ifmathtags@
\def\TagsAsMath{\global\mathtags@true}
\def\TagsAsText{\global\mathtags@false}
\TagsAsText
\def\tagform@#1{\hbox{\rm(\ignorespaces#1\unskip)}}
\def\thetag{\leavevmode\tagform@}
\def\tag#1$${\iftagsleft@\leqno\else\eqno\fi                                
 \maketag@#1\maketag@                                                       
 $$}                                                                        
\def\maketag@{\FN@\maketag@@}
\def\maketag@@{\ifx\next"\expandafter\maketag@@@\else\expandafter\maketag@@@@
 \fi}
\def\maketag@@@"#1"#2\maketag@{\hbox{\rm#1}}                                
\def\maketag@@@@#1\maketag@{\ifmathtags@\tagform@{$\m@th#1$}\else
 \tagform@{#1}\fi}
\interdisplaylinepenalty\@M
\def\allowdisplaybreaks{\RIfMIfI@
 \onlydmatherr@\allowdisplaybreaks\else
 \interdisplaylinepenalty\z@\fi\else\onlydmatherr@\allowdisplaybreaks\fi}
\Invalid@\allowdisplaybreak
\Invalid@\displaybreak
\Invalid@\intertext
\def\allowdisplaybreak@{\def\allowdisplaybreak{\crcr\noalign{\allowbreak}}}
\def\displaybreak@{\def\displaybreak{\crcr\noalign{\break}}}
\def\intertext@{\def\intertext##1{\crcr\noalign{%
 \penalty\postdisplaypenalty \vskip\belowdisplayskip
 \vbox{\normalbaselines\noindent##1}%
 \penalty\predisplaypenalty \vskip\abovedisplayskip}}}
\newskip\centering@
\centering@\z@ plus\@m\p@
\def\align{\relax\ifingather@\DN@{\csname align (in
  \string\gather)\endcsname}\else
 \ifmmode\ifinner\DN@{\onlydmatherr@\align}\else
  \let\next@\align@\fi
 \else\DN@{\onlydmatherr@\align}\fi\fi\next@}
\newhelp\andhelp@
{An extra & here is so disastrous that you should probably exit^^J
and fix things up.}
\newif\iftag@
\newcount\and@
\def\align@{\inalign@true\inany@true
 \vspace@\allowdisplaybreak@\displaybreak@\intertext@
 \def\tag{\global\tag@true\ifnum\and@=\z@\DN@{&&}\else
          \DN@{&}\fi\next@}%
 \iftagsleft@\DN@{\csname align \endcsname}\else
  \DN@{\csname align \space\endcsname}\fi\next@}
\def\Tag@{\iftag@\else\errhelp\andhelp@\err@{Extra & on this line}\fi}
\newdimen\lwidth@
\newdimen\rwidth@
\newdimen\maxlwidth@
\newdimen\maxrwidth@
\newdimen\totwidth@
\def\measure@#1\endalign{\lwidth@\z@\rwidth@\z@\maxlwidth@\z@\maxrwidth@\z@
 \global\and@\z@                                                            
 \setbox@ne\vbox                                                            
  {\everycr{\noalign{\global\tag@false\global\and@\z@}}\Let@                
  \halign{\setboxz@h{$\m@th\displaystyle{\@lign##}$}
   \global\lwidth@\wdz@                                                     
   \ifdim\lwidth@>\maxlwidth@\global\maxlwidth@\lwidth@\fi                  
   \global\advance\and@\@ne                                                 
   &\setboxz@h{$\m@th\displaystyle{{}\@lign##}$}\global\rwidth@\wdz@        
   \ifdim\rwidth@>\maxrwidth@\global\maxrwidth@\rwidth@\fi                  
   \global\advance\and@\@ne                                                
   &\Tag@
   \eat@{##}\crcr#1\crcr}}
 \totwidth@\maxlwidth@\advance\totwidth@\maxrwidth@}                       
\def\displ@y@{\global\dt@ptrue\openup\jot
 \everycr{\noalign{\global\tag@false\global\and@\z@\ifdt@p\global\dt@pfalse
 \vskip-\lineskiplimit\vskip\normallineskiplimit\else
 \penalty\interdisplaylinepenalty\fi}}}
\def\black@#1{\noalign{\ifdim#1>\displaywidth
 \dimen@\prevdepth\nointerlineskip                                          
 \vskip-\ht\strutbox@\vskip-\dp\strutbox@                                   
 \vbox{\noindent\hbox to#1{\strut@\hfill}}
 \prevdepth\dimen@                                                          
 \fi}}
\expandafter\def\csname align \space\endcsname#1\endalign
 {\measure@#1\endalign\global\and@\z@                                       
 \ifingather@\everycr{\noalign{\global\and@\z@}}\else\displ@y@\fi           
 \Let@\tabskip\centering@                                                   
 \halign to\displaywidth
  {\hfil\strut@\setboxz@h{$\m@th\displaystyle{\@lign##}$}
  \global\lwidth@\wdz@\boxz@\global\advance\and@\@ne                        
  \tabskip\z@skip                                                           
  &\setboxz@h{$\m@th\displaystyle{{}\@lign##}$}
  \global\rwidth@\wdz@\boxz@\hfill\global\advance\and@\@ne                  
  \tabskip\centering@                                                       
  &\setboxz@h{\@lign\strut@\maketag@##\maketag@}
  \dimen@\displaywidth\advance\dimen@-\totwidth@
  \divide\dimen@\tw@\advance\dimen@\maxrwidth@\advance\dimen@-\rwidth@     
  \ifdim\dimen@<\tw@\wdz@\llap{\vtop{\normalbaselines\null\boxz@}}
  \else\llap{\boxz@}\fi                                                    
  \tabskip\z@skip                                                          
  \crcr#1\crcr                                                             
  \black@\totwidth@}}                                                      
\newdimen\lineht@
\expandafter\def\csname align \endcsname#1\endalign{\measure@#1\endalign
 \global\and@\z@
 \ifdim\totwidth@>\displaywidth\let\displaywidth@\totwidth@\else
  \let\displaywidth@\displaywidth\fi                                        
 \ifingather@\everycr{\noalign{\global\and@\z@}}\else\displ@y@\fi
 \Let@\tabskip\centering@\halign to\displaywidth
  {\hfil\strut@\setboxz@h{$\m@th\displaystyle{\@lign##}$}%
  \global\lwidth@\wdz@\global\lineht@\ht\z@                                 
  \boxz@\global\advance\and@\@ne
  \tabskip\z@skip&\setboxz@h{$\m@th\displaystyle{{}\@lign##}$}%
  \global\rwidth@\wdz@\ifdim\ht\z@>\lineht@\global\lineht@\ht\z@\fi         
  \boxz@\hfil\global\advance\and@\@ne
  \tabskip\centering@&\kern-\displaywidth@                                  
  \setboxz@h{\@lign\strut@\maketag@##\maketag@}%
  \dimen@\displaywidth\advance\dimen@-\totwidth@
  \divide\dimen@\tw@\advance\dimen@\maxlwidth@\advance\dimen@-\lwidth@
  \ifdim\dimen@<\tw@\wdz@
   \rlap{\vbox{\normalbaselines\boxz@\vbox to\lineht@{}}}\else
   \rlap{\boxz@}\fi
  \tabskip\displaywidth@\crcr#1\crcr\black@\totwidth@}}
\expandafter\def\csname align (in \string\gather)\endcsname
  #1\endalign{\vcenter{\align@#1\endalign}}
\Invalid@\endalign
\newif\ifxat@
\def\alignat{\RIfMIfI@\DN@{\onlydmatherr@\alignat}\else
 \DN@{\csname alignat \endcsname}\fi\else
 \DN@{\onlydmatherr@\alignat}\fi\next@}
\newif\ifmeasuring@
\newbox\savealignat@
\expandafter\def\csname alignat \endcsname#1#2\endalignat                   
 {\inany@true\xat@false
 \def\tag{\global\tag@true\count@#1\relax\multiply\count@\tw@
  \xdef\tag@{}\loop\ifnum\count@>\and@\xdef\tag@{&\tag@}\advance\count@\m@ne
  \repeat\tag@}%
 \vspace@\allowdisplaybreak@\displaybreak@\intertext@
 \displ@y@\measuring@true                                                   
 \setbox\savealignat@\hbox{$\m@th\displaystyle\Let@
  \attag@{#1}
  \vbox{\halign{\span\preamble@@\crcr#2\crcr}}$}%
 \measuring@false                                                           
 \Let@\attag@{#1}
 \tabskip\centering@\halign to\displaywidth
  {\span\preamble@@\crcr#2\crcr                                             
  \black@{\wd\savealignat@}}}                                               
\Invalid@\endalignat
\def\xalignat{\RIfMIfI@
 \DN@{\onlydmatherr@\xalignat}\else
 \DN@{\csname xalignat \endcsname}\fi\else
 \DN@{\onlydmatherr@\xalignat}\fi\next@}
\expandafter\def\csname xalignat \endcsname#1#2\endxalignat
 {\inany@true\xat@true
 \def\tag{\global\tag@true\def\tag@{}\count@#1\relax\multiply\count@\tw@
  \loop\ifnum\count@>\and@\xdef\tag@{&\tag@}\advance\count@\m@ne\repeat\tag@}%
 \vspace@\allowdisplaybreak@\displaybreak@\intertext@
 \displ@y@\measuring@true\setbox\savealignat@\hbox{$\m@th\displaystyle\Let@
 \attag@{#1}\vbox{\halign{\span\preamble@@\crcr#2\crcr}}$}%
 \measuring@false\Let@
 \attag@{#1}\tabskip\centering@\halign to\displaywidth
 {\span\preamble@@\crcr#2\crcr\black@{\wd\savealignat@}}}
\def\attag@#1{\let\Maketag@\maketag@\let\TAG@\Tag@                          
 \let\Tag@=0\let\maketag@=0
 \ifmeasuring@\def\llap@##1{\setboxz@h{##1}\hbox to\tw@\wdz@{}}%
  \def\rlap@##1{\setboxz@h{##1}\hbox to\tw@\wdz@{}}\else
  \let\llap@\llap\let\rlap@\rlap\fi                                         
 \toks@{\hfil\strut@$\m@th\displaystyle{\@lign\the\hashtoks@}$\tabskip\z@skip
  \global\advance\and@\@ne&$\m@th\displaystyle{{}\@lign\the\hashtoks@}$\hfil
  \ifxat@\tabskip\centering@\fi\global\advance\and@\@ne}
 \iftagsleft@
  \toks@@{\tabskip\centering@&\Tag@\kern-\displaywidth
   \rlap@{\@lign\maketag@\the\hashtoks@\maketag@}%
   \global\advance\and@\@ne\tabskip\displaywidth}\else
  \toks@@{\tabskip\centering@&\Tag@\llap@{\@lign\maketag@
   \the\hashtoks@\maketag@}\global\advance\and@\@ne\tabskip\z@skip}\fi      
 \atcount@#1\relax\advance\atcount@\m@ne
 \loop\ifnum\atcount@>\z@
 \toks@=\expandafter{\the\toks@&\hfil$\m@th\displaystyle{\@lign
  \the\hashtoks@}$\global\advance\and@\@ne
  \tabskip\z@skip&$\m@th\displaystyle{{}\@lign\the\hashtoks@}$\hfil\ifxat@
  \tabskip\centering@\fi\global\advance\and@\@ne}\advance\atcount@\m@ne
 \repeat                                                                    
 \xdef\preamble@{\the\toks@\the\toks@@}
 \xdef\preamble@@{\preamble@}
 \let\maketag@\Maketag@\let\Tag@\TAG@}                                      
\Invalid@\endxalignat
\def\xxalignat{\RIfMIfI@
 \DN@{\onlydmatherr@\xxalignat}\else\DN@{\csname xxalignat
  \endcsname}\fi\else
 \DN@{\onlydmatherr@\xxalignat}\fi\next@}
\expandafter\def\csname xxalignat \endcsname#1#2\endxxalignat{\inany@true
 \vspace@\allowdisplaybreak@\displaybreak@\intertext@
 \displ@y\setbox\savealignat@\hbox{$\m@th\displaystyle\Let@
 \xxattag@{#1}\vbox{\halign{\span\preamble@@\crcr#2\crcr}}$}%
 \Let@\xxattag@{#1}\tabskip\z@skip\halign to\displaywidth
 {\span\preamble@@\crcr#2\crcr\black@{\wd\savealignat@}}}
\def\xxattag@#1{\toks@{\tabskip\z@skip\hfil\strut@
 $\m@th\displaystyle{\the\hashtoks@}$&%
 $\m@th\displaystyle{{}\the\hashtoks@}$\hfil\tabskip\centering@&}%
 \atcount@#1\relax\advance\atcount@\m@ne\loop\ifnum\atcount@>\z@
 \toks@=\expandafter{\the\toks@&\hfil$\m@th\displaystyle{\the\hashtoks@}$%
  \tabskip\z@skip&$\m@th\displaystyle{{}\the\hashtoks@}$\hfil
  \tabskip\centering@}\advance\atcount@\m@ne\repeat
 \xdef\preamble@{\the\toks@\tabskip\z@skip}\xdef\preamble@@{\preamble@}}
\Invalid@\endxxalignat
\newdimen\gwidth@
\newdimen\gmaxwidth@
\def\gmeasure@#1\endgather{\gwidth@\z@\gmaxwidth@\z@\setbox@ne\vbox{\Let@
 \halign{\setboxz@h{$\m@th\displaystyle{##}$}\global\gwidth@\wdz@
 \ifdim\gwidth@>\gmaxwidth@\global\gmaxwidth@\gwidth@\fi
 &\eat@{##}\crcr#1\crcr}}}
\def\gather{\RIfMIfI@\DN@{\onlydmatherr@\gather}\else
 \ingather@true\inany@true\def\tag{&}%
 \vspace@\allowdisplaybreak@\displaybreak@\intertext@
 \displ@y\Let@
 \iftagsleft@\DN@{\csname gather \endcsname}\else
  \DN@{\csname gather \space\endcsname}\fi\fi
 \else\DN@{\onlydmatherr@\gather}\fi\next@}
\expandafter\def\csname gather \space\endcsname#1\endgather
 {\gmeasure@#1\endgather\tabskip\centering@
 \halign to\displaywidth{\hfil\strut@\setboxz@h{$\m@th\displaystyle{##}$}%
 \global\gwidth@\wdz@\boxz@\hfil&
 \setboxz@h{\strut@{\maketag@##\maketag@}}%
 \dimen@\displaywidth\advance\dimen@-\gwidth@
 \ifdim\dimen@>\tw@\wdz@\llap{\boxz@}\else
 \llap{\vtop{\normalbaselines\null\boxz@}}\fi
 \tabskip\z@skip\crcr#1\crcr\black@\gmaxwidth@}}
\newdimen\glineht@
\expandafter\def\csname gather \endcsname#1\endgather{\gmeasure@#1\endgather
 \ifdim\gmaxwidth@>\displaywidth\let\gdisplaywidth@\gmaxwidth@\else
 \let\gdisplaywidth@\displaywidth\fi\tabskip\centering@\halign to\displaywidth
 {\hfil\strut@\setboxz@h{$\m@th\displaystyle{##}$}%
 \global\gwidth@\wdz@\global\glineht@\ht\z@\boxz@\hfil&\kern-\gdisplaywidth@
 \setboxz@h{\strut@{\maketag@##\maketag@}}%
 \dimen@\displaywidth\advance\dimen@-\gwidth@
 \ifdim\dimen@>\tw@\wdz@\rlap{\boxz@}\else
 \rlap{\vbox{\normalbaselines\boxz@\vbox to\glineht@{}}}\fi
 \tabskip\gdisplaywidth@\crcr#1\crcr\black@\gmaxwidth@}}
\newif\ifctagsplit@
\def\CenteredTagsOnSplits{\global\ctagsplit@true}
\def\TopOrBottomTagsOnSplits{\global\ctagsplit@false}
\TopOrBottomTagsOnSplits
\def\split{\relax\ifinany@\let\next@\insplit@\else
 \ifmmode\ifinner\def\next@{\onlydmatherr@\split}\else
 \let\next@\outsplit@\fi\else
 \def\next@{\onlydmatherr@\split}\fi\fi\next@}
\def\insplit@{\global\setbox\z@\vbox\bgroup\vspace@\Let@\ialign\bgroup
 \hfil\strut@$\m@th\displaystyle{##}$&$\m@th\displaystyle{{}##}$\hfill\crcr}
\def\endsplit{\crcr\egroup\egroup\iftagsleft@\expandafter\lendsplit@\else
 \expandafter\rendsplit@\fi}
\def\rendsplit@{\global\setbox9 \vbox
 {\unvcopy\z@\global\setbox8 \lastbox\unskip}
 \setbox@ne\hbox{\unhcopy8 \unskip\global\setbox\tw@\lastbox
 \unskip\global\setbox\thr@@\lastbox}
 \global\setbox7 \hbox{\unhbox\tw@\unskip}
 \ifinalign@\ifctagsplit@                                                   
  \gdef\split@{\hbox to\wd\thr@@{}&
   \vcenter{\vbox{\moveleft\wd\thr@@\boxz@}}}
 \else\gdef\split@{&\vbox{\moveleft\wd\thr@@\box9}\crcr
  \box\thr@@&\box7}\fi                                                      
 \else                                                                      
  \ifctagsplit@\gdef\split@{\vcenter{\boxz@}}\else
  \gdef\split@{\box9\crcr\hbox{\box\thr@@\box7}}\fi
 \fi
 \split@}                                                                   
\def\lendsplit@{\global\setbox9\vtop{\unvcopy\z@}
 \setbox@ne\vbox{\unvcopy\z@\global\setbox8\lastbox}
 \setbox@ne\hbox{\unhcopy8\unskip\setbox\tw@\lastbox
  \unskip\global\setbox\thr@@\lastbox}
 \ifinalign@\ifctagsplit@                                                   
  \gdef\split@{\hbox to\wd\thr@@{}&
  \vcenter{\vbox{\moveleft\wd\thr@@\box9}}}
  \else                                                                     
  \gdef\split@{\hbox to\wd\thr@@{}&\vbox{\moveleft\wd\thr@@\box9}}\fi
 \else
  \ifctagsplit@\gdef\split@{\vcenter{\box9}}\else
  \gdef\split@{\box9}\fi
 \fi\split@}
\def\outsplit@#1$${\align\insplit@#1\endalign$$}
\newdimen\multlinegap@
\multlinegap@1em
\newdimen\multlinetaggap@
\multlinetaggap@1em
\def\MultlineGap#1{\global\multlinegap@#1\relax}
\def\multlinegap#1{\RIfMIfI@\onlydmatherr@\multlinegap\else
 \multlinegap@#1\relax\fi\else\onlydmatherr@\multlinegap\fi}
\def\nomultlinegap{\multlinegap{\z@}}
\def\multline{\RIfMIfI@
 \DN@{\onlydmatherr@\multline}\else
 \DN@{\multline@}\fi\else
 \DN@{\onlydmatherr@\multline}\fi\next@}
\newif\iftagin@
\def\tagin@#1{\tagin@false\in@\tag{#1}\ifin@\tagin@true\fi}
\def\multline@#1$${\inany@true\vspace@\allowdisplaybreak@\displaybreak@
 \tagin@{#1}\iftagsleft@\DN@{\multline@l#1$$}\else
 \DN@{\multline@r#1$$}\fi\next@}
\newdimen\mwidth@
\def\rmmeasure@#1\endmultline{%
 \def\shoveleft##1{##1}\def\shoveright##1{##1}
 \setbox@ne\vbox{\Let@\halign{\setboxz@h
  {$\m@th\@lign\displaystyle{}##$}\global\mwidth@\wdz@
  \crcr#1\crcr}}}
\newdimen\mlineht@
\newif\ifzerocr@
\newif\ifonecr@
\def\lmmeasure@#1\endmultline{\global\zerocr@true\global\onecr@false
 \everycr{\noalign{\ifonecr@\global\onecr@false\fi
  \ifzerocr@\global\zerocr@false\global\onecr@true\fi}}
  \def\shoveleft##1{##1}\def\shoveright##1{##1}%
 \setbox@ne\vbox{\Let@\halign{\setboxz@h
  {$\m@th\@lign\displaystyle{}##$}\ifonecr@\global\mwidth@\wdz@
  \global\mlineht@\ht\z@\fi\crcr#1\crcr}}}
\newbox\mtagbox@
\newdimen\ltwidth@
\newdimen\rtwidth@
\def\multline@l#1$${\iftagin@\DN@{\lmultline@@#1$$}\else
 \DN@{\setbox\mtagbox@\null\ltwidth@\z@\rtwidth@\z@
  \lmultline@@@#1$$}\fi\next@}
\def\lmultline@@#1\endmultline\tag#2$${%
 \setbox\mtagbox@\hbox{\maketag@#2\maketag@}
 \lmmeasure@#1\endmultline\dimen@\mwidth@\advance\dimen@\wd\mtagbox@
 \advance\dimen@\multlinetaggap@                                            
 \ifdim\dimen@>\displaywidth\ltwidth@\z@\else\ltwidth@\wd\mtagbox@\fi       
 \lmultline@@@#1\endmultline$$}
\def\lmultline@@@{\displ@y
 \def\shoveright##1{##1\hfilneg\hskip\multlinegap@}%
 \def\shoveleft##1{\setboxz@h{$\m@th\displaystyle{}##1$}%
  \setbox@ne\hbox{$\m@th\displaystyle##1$}%
  \hfilneg
  \iftagin@
   \ifdim\ltwidth@>\z@\hskip\ltwidth@\hskip\multlinetaggap@\fi
  \else\hskip\multlinegap@\fi\hskip.5\wd@ne\hskip-.5\wdz@##1}
  \halign\bgroup\Let@\hbox to\displaywidth
   {\strut@$\m@th\displaystyle\hfil{}##\hfil$}\crcr
   \hfilneg                                                                 
   \iftagin@                                                                
    \ifdim\ltwidth@>\z@                                                     
     \box\mtagbox@\hskip\multlinetaggap@                                    
    \else
     \rlap{\vbox{\normalbaselines\hbox{\strut@\box\mtagbox@}%
     \vbox to\mlineht@{}}}\fi                                               
   \else\hskip\multlinegap@\fi}                                             
\def\multline@r#1$${\iftagin@\DN@{\rmultline@@#1$$}\else
 \DN@{\setbox\mtagbox@\null\ltwidth@\z@\rtwidth@\z@
  \rmultline@@@#1$$}\fi\next@}
\def\rmultline@@#1\endmultline\tag#2$${\ltwidth@\z@
 \setbox\mtagbox@\hbox{\maketag@#2\maketag@}%
 \rmmeasure@#1\endmultline\dimen@\mwidth@\advance\dimen@\wd\mtagbox@
 \advance\dimen@\multlinetaggap@
 \ifdim\dimen@>\displaywidth\rtwidth@\z@\else\rtwidth@\wd\mtagbox@\fi
 \rmultline@@@#1\endmultline$$}
\def\rmultline@@@{\displ@y
 \def\shoveright##1{##1\hfilneg\iftagin@\ifdim\rtwidth@>\z@
  \hskip\rtwidth@\hskip\multlinetaggap@\fi\else\hskip\multlinegap@\fi}%
 \def\shoveleft##1{\setboxz@h{$\m@th\displaystyle{}##1$}%
  \setbox@ne\hbox{$\m@th\displaystyle##1$}%
  \hfilneg\hskip\multlinegap@\hskip.5\wd@ne\hskip-.5\wdz@##1}%
 \halign\bgroup\Let@\hbox to\displaywidth
  {\strut@$\m@th\displaystyle\hfil{}##\hfil$}\crcr
 \hfilneg\hskip\multlinegap@}
\def\endmultline{\iftagsleft@\expandafter\lendmultline@\else
 \expandafter\rendmultline@\fi}
\def\lendmultline@{\hfilneg\hskip\multlinegap@\crcr\egroup}
\def\rendmultline@{\iftagin@                                                
 \ifdim\rtwidth@>\z@                                                        
  \hskip\multlinetaggap@\box\mtagbox@                                       
 \else\llap{\vtop{\normalbaselines\null\hbox{\strut@\box\mtagbox@}}}\fi     
 \else\hskip\multlinegap@\fi                                                
 \hfilneg\crcr\egroup}
\def\bmod{\mskip-\medmuskip\mkern5mu\mathbin{\fam\z@ mod}\penalty900
 \mkern5mu\mskip-\medmuskip}
\def\pmod#1{\allowbreak\ifinner\mkern8mu\else\mkern18mu\fi
 ({\fam\z@ mod}\,\,#1)}
\def\pod#1{\allowbreak\ifinner\mkern8mu\else\mkern18mu\fi(#1)}
\def\mod#1{\allowbreak\ifinner\mkern12mu\else\mkern18mu\fi{\fam\z@ mod}\,\,#1}
\message{continued fractions,}
\newcount\cfraccount@
\def\cfrac{\bgroup\bgroup\advance\cfraccount@\@ne\strut
 \iffalse{\fi\def\\{\over\displaystyle}\iffalse}\fi}
\def\lcfrac{\bgroup\bgroup\advance\cfraccount@\@ne\strut
 \iffalse{\fi\def\\{\hfill\over\displaystyle}\iffalse}\fi}
\def\rcfrac{\bgroup\bgroup\advance\cfraccount@\@ne\strut\hfill
 \iffalse{\fi\def\\{\over\displaystyle}\iffalse}\fi}
\def\gloop@#1\repeat{\gdef\body{#1}\iterate}
\def\endcfrac{\gloop@\ifnum\cfraccount@>\z@\global\advance\cfraccount@\m@ne
 \egroup\hskip-\nulldelimiterspace\egroup\repeat}
\message{compound symbols,}
\def\binrel@#1{\setboxz@h{\thinmuskip0mu
  \medmuskip\m@ne mu\thickmuskip\@ne mu$#1\m@th$}%
 \setbox@ne\hbox{\thinmuskip0mu\medmuskip\m@ne mu\thickmuskip
  \@ne mu${}#1{}\m@th$}%
 \setbox\tw@\hbox{\hskip\wd@ne\hskip-\wdz@}}
\def\overset#1\to#2{\binrel@{#2}\ifdim\wd\tw@<\z@
 \mathbin{\mathop{\kern\z@#2}\limits^{#1}}\else\ifdim\wd\tw@>\z@
 \mathrel{\mathop{\kern\z@#2}\limits^{#1}}\else
 {\mathop{\kern\z@#2}\limits^{#1}}{}\fi\fi}
\def\underset#1\to#2{\binrel@{#2}\ifdim\wd\tw@<\z@
 \mathbin{\mathop{\kern\z@#2}\limits_{#1}}\else\ifdim\wd\tw@>\z@
 \mathrel{\mathop{\kern\z@#2}\limits_{#1}}\else
 {\mathop{\kern\z@#2}\limits_{#1}}{}\fi\fi}
\def\oversetbrace#1\to#2{\overbrace{#2}^{#1}}
\def\undersetbrace#1\to#2{\underbrace{#2}_{#1}}
\def\sideset#1\and#2\to#3{%
 \setbox@ne\hbox{$\dsize{\vphantom{#3}}#1{#3}\m@th$}%
 \setbox\tw@\hbox{$\dsize{#3}#2\m@th$}%
 \hskip\wd@ne\hskip-\wd\tw@\mathop{\hskip\wd\tw@\hskip-\wd@ne
  {\vphantom{#3}}#1{#3}#2}}
\def\rightarrowfill@#1{\setboxz@h{$#1-\m@th$}\ht\z@\z@
  $#1\m@th\copy\z@\mkern-6mu\cleaders
  \hbox{$#1\mkern-2mu\box\z@\mkern-2mu$}\hfill
  \mkern-6mu\mathord\rightarrow$}
\def\leftarrowfill@#1{\setboxz@h{$#1-\m@th$}\ht\z@\z@
  $#1\m@th\mathord\leftarrow\mkern-6mu\cleaders
  \hbox{$#1\mkern-2mu\copy\z@\mkern-2mu$}\hfill
  \mkern-6mu\box\z@$}
\def\leftrightarrowfill@#1{\setboxz@h{$#1-\m@th$}\ht\z@\z@
  $#1\m@th\mathord\leftarrow\mkern-6mu\cleaders
  \hbox{$#1\mkern-2mu\box\z@\mkern-2mu$}\hfill
  \mkern-6mu\mathord\rightarrow$}
\def\overrightarrow{\mathpalette\overrightarrow@}
\def\overrightarrow@#1#2{\vbox{\ialign{##\crcr\rightarrowfill@#1\crcr
 \noalign{\kern-\ex@\nointerlineskip}$\m@th\hfil#1#2\hfil$\crcr}}}

\def\overleftarrow{\mathpalette\overleftarrow@}
\def\overleftarrow@#1#2{\vbox{\ialign{##\crcr\leftarrowfill@#1\crcr
 \noalign{\kern-\ex@\nointerlineskip}$\m@th\hfil#1#2\hfil$\crcr}}}
\def\overleftrightarrow{\mathpalette\overleftrightarrow@}
\def\overleftrightarrow@#1#2{\vbox{\ialign{##\crcr\leftrightarrowfill@#1\crcr
 \noalign{\kern-\ex@\nointerlineskip}$\m@th\hfil#1#2\hfil$\crcr}}}
\def\underrightarrow{\mathpalette\underrightarrow@}
\def\underrightarrow@#1#2{\vtop{\ialign{##\crcr$\m@th\hfil#1#2\hfil$\crcr
 \noalign{\nointerlineskip}\rightarrowfill@#1\crcr}}}

\def\underleftarrow{\mathpalette\underleftarrow@}
\def\underleftarrow@#1#2{\vtop{\ialign{##\crcr$\m@th\hfil#1#2\hfil$\crcr
 \noalign{\nointerlineskip}\leftarrowfill@#1\crcr}}}
\def\underleftrightarrow{\mathpalette\underleftrightarrow@}
\def\underleftrightarrow@#1#2{\vtop{\ialign{##\crcr$\m@th\hfil#1#2\hfil$\crcr
 \noalign{\nointerlineskip}\leftrightarrowfill@#1\crcr}}}
\message{various kinds of dots,}
\let\DOTSI\relax
\let\DOTSB\relax

\newif\ifmath@
{\uccode`7=`\\ \uccode`8=`m \uccode`9=`a \uccode`0=`t \uccode`!=`h
 \uppercase{\gdef\math@#1#2#3#4#5#6\math@{\global\math@false\ifx 7#1\ifx 8#2%
 \ifx 9#3\ifx 0#4\ifx !#5\xdef\meaning@{#6}\global\math@true\fi\fi\fi\fi\fi}}}
\newif\ifmathch@
{\uccode`7=`c \uccode`8=`h \uccode`9=`\"
 \uppercase{\gdef\mathch@#1#2#3#4#5#6\mathch@{\global\mathch@false
  \ifx 7#1\ifx 8#2\ifx 9#5\global\mathch@true\xdef\meaning@{9#6}\fi\fi\fi}}}
\newcount\classnum@
\def\getmathch@#1.#2\getmathch@{\classnum@#1 \divide\classnum@4096
 \ifcase\number\classnum@\or\or\gdef\thedots@{\dotsb@}\or
 \gdef\thedots@{\dotsb@}\fi}
\newif\ifmathbin@
{\uccode`4=`b \uccode`5=`i \uccode`6=`n
 \uppercase{\gdef\mathbin@#1#2#3{\relaxnext@
  \DNii@##1\mathbin@{\ifx\space@\next\global\mathbin@true\fi}%
 \global\mathbin@false\DN@##1\mathbin@{}%
 \ifx 4#1\ifx 5#2\ifx 6#3\DN@{\FN@\nextii@}\fi\fi\fi\next@}}}
\newif\ifmathrel@
{\uccode`4=`r \uccode`5=`e \uccode`6=`l
 \uppercase{\gdef\mathrel@#1#2#3{\relaxnext@
  \DNii@##1\mathrel@{\ifx\space@\next\global\mathrel@true\fi}%
 \global\mathrel@false\DN@##1\mathrel@{}%
 \ifx 4#1\ifx 5#2\ifx 6#3\DN@{\FN@\nextii@}\fi\fi\fi\next@}}}
\newif\ifmacro@
{\uccode`5=`m \uccode`6=`a \uccode`7=`c
 \uppercase{\gdef\macro@#1#2#3#4\macro@{\global\macro@false
  \ifx 5#1\ifx 6#2\ifx 7#3\global\macro@true
  \xdef\meaning@{\macro@@#4\macro@@}\fi\fi\fi}}}
\def\macro@@#1->#2\macro@@{#2}
\newif\ifDOTS@
\newcount\DOTSCASE@
{\uccode`6=`\\ \uccode`7=`D \uccode`8=`O \uccode`9=`T \uccode`0=`S
 \uppercase{\gdef\DOTS@#1#2#3#4#5{\global\DOTS@false\DN@##1\DOTS@{}%
  \ifx 6#1\ifx 7#2\ifx 8#3\ifx 9#4\ifx 0#5\let\next@\DOTS@@\fi\fi\fi\fi\fi
  \next@}}}
{\uccode`3=`B \uccode`4=`I \uccode`5=`X
 \uppercase{\gdef\DOTS@@#1{\relaxnext@
  \DNii@##1\DOTS@{\ifx\space@\next\global\DOTS@true\fi}%
  \DN@{\FN@\nextii@}%
  \ifx 3#1\global\DOTSCASE@\z@\else
  \ifx 4#1\global\DOTSCASE@\@ne\else
  \ifx 5#1\global\DOTSCASE@\tw@\else\DN@##1\DOTS@{}%
  \fi\fi\fi\next@}}}
\newif\ifnot@
{\uccode`5=`\\ \uccode`6=`n \uccode`7=`o \uccode`8=`t
 \uppercase{\gdef\not@#1#2#3#4{\relaxnext@
  \DNii@##1\not@{\ifx\space@\next\global\not@true\fi}%
 \global\not@false\DN@##1\not@{}%
 \ifx 5#1\ifx 6#2\ifx 7#3\ifx 8#4\DN@{\FN@\nextii@}\fi\fi\fi
 \fi\next@}}}
\newif\ifkeybin@
\def\keybin@{\keybin@true
 \ifx\next+\else\ifx\next=\else\ifx\next<\else\ifx\next>\else\ifx\next-\else
 \ifx\next*\else\ifx\next:\else\keybin@false\fi\fi\fi\fi\fi\fi\fi}
\def\dots{\RIfM@\expandafter\mdots@\else\expandafter\tdots@\fi}
\def\tdots@{\unskip\relaxnext@
 \DN@{$\m@th\mathinner{\ldotp\ldotp\ldotp}\,
   \ifx\next,\,$\else\ifx\next.\,$\else\ifx\next;\,$\else\ifx\next:\,$\else
   \ifx\next?\,$\else\ifx\next!\,$\else$ \fi\fi\fi\fi\fi\fi}%
 \ \FN@\next@}
\def\mdots@{\FN@\mdots@@}
\def\mdots@@{\gdef\thedots@{\dotso@}
 \ifx\next\boldkey\gdef\thedots@\boldkey{\boldkeydots@}\else                
 \ifx\next\boldsymbol\gdef\thedots@\boldsymbol{\boldsymboldots@}\else       
 \ifx,\next\gdef\thedots@{\dotsc}
 \else\ifx\not\next\gdef\thedots@{\dotsb@}
 \else\keybin@
 \ifkeybin@\gdef\thedots@{\dotsb@}
 \else\xdef\meaning@{\meaning\next..........}\xdef\meaning@@{\meaning@}
  \expandafter\math@\meaning@\math@
  \ifmath@
   \expandafter\mathch@\meaning@\mathch@
   \ifmathch@\expandafter\getmathch@\meaning@\getmathch@\fi                 
  \else\expandafter\macro@\meaning@@\macro@                                 
  \ifmacro@                                                                
   \expandafter\not@\meaning@\not@\ifnot@\gdef\thedots@{\dotsb@}
  \else\expandafter\DOTS@\meaning@\DOTS@
  \ifDOTS@
   \ifcase\number\DOTSCASE@\gdef\thedots@{\dotsb@}%
    \or\gdef\thedots@{\dotsi}\else\fi                                      
  \else\expandafter\math@\meaning@\math@                                   
  \ifmath@\expandafter\mathbin@\meaning@\mathbin@
  \ifmathbin@\gdef\thedots@{\dotsb@}
  \else\expandafter\mathrel@\meaning@\mathrel@
  \ifmathrel@\gdef\thedots@{\dotsb@}
  \fi\fi\fi\fi\fi\fi\fi\fi\fi\fi\fi\fi
 \thedots@}
\def\plainldots@{\mathinner{\ldotp\ldotp\ldotp}}
\def\plaincdots@{\mathinner{\cdotp\cdotp\cdotp}}
\def\dotsi{\!\plaincdots@}
\let\dotsb@\plaincdots@
\newif\ifextra@
\newif\ifrightdelim@
\def\rightdelim@{\global\rightdelim@true                                    
 \ifx\next)\else                                                            
 \ifx\next]\else
 \ifx\next\rbrack\else
 \ifx\next\}\else
 \ifx\next\rbrace\else
 \ifx\next\rangle\else
 \ifx\next\rceil\else
 \ifx\next\rfloor\else
 \ifx\next\rgroup\else
 \ifx\next\rmoustache\else
 \ifx\next\right\else
 \ifx\next\bigr\else
 \ifx\next\biggr\else
 \ifx\next\Bigr\else                                                        
 \ifx\next\Biggr\else\global\rightdelim@false
 \fi\fi\fi\fi\fi\fi\fi\fi\fi\fi\fi\fi\fi\fi\fi}
\def\extra@{%
 \global\extra@false\rightdelim@\ifrightdelim@\global\extra@true            
 \else\ifx\next$\global\extra@true                                          
 \else\xdef\meaning@{\meaning\next..........}
 \expandafter\macro@\meaning@\macro@\ifmacro@                               
 \expandafter\DOTS@\meaning@\DOTS@
 \ifDOTS@
 \ifnum\DOTSCASE@=\tw@\global\extra@true                                    
 \fi\fi\fi\fi\fi}
\newif\ifbold@
\def\dotso@{\relaxnext@
 \ifbold@
  \let\next\delayed@
  \DNii@{\extra@\plainldots@\ifextra@\,\fi}%
 \else
  \DNii@{\DN@{\extra@\plainldots@\ifextra@\,\fi}\FN@\next@}%
 \fi
 \nextii@}
\def\extrap@#1{%
 \ifx\next,\DN@{#1\,}\else
 \ifx\next;\DN@{#1\,}\else
 \ifx\next.\DN@{#1\,}\else\extra@
 \ifextra@\DN@{#1\,}\else
 \let\next@#1\fi\fi\fi\fi\next@}
\def\ldots{\DN@{\extrap@\plainldots@}%
 \FN@\next@}
\def\cdots{\DN@{\extrap@\plaincdots@}%
 \FN@\next@}

\def\dotsc{\relaxnext@
 \DN@{\ifx\next;\plainldots@\,\else
  \ifx\next.\plainldots@\,\else\extra@\plainldots@
  \ifextra@\,\fi\fi\fi}%
 \FN@\next@}
\def\cdot{\mathchar"2201 }
\def\longrightarrow{\DOTSB\relbar\joinrel\rightarrow}

\def\mapsto{\DOTSB\mapstochar\rightarrow}

\message{special superscripts,}
\def\dddot#1{{\mathop{#1}\limits^{\vbox to-1.4\ex@{\kern-\tw@\ex@
 \hbox{\rm...}\vss}}}}
\def\ddddot#1{{\mathop{#1}\limits^{\vbox to-1.4\ex@{\kern-\tw@\ex@
 \hbox{\rm....}\vss}}}}
\def\sphat{^{\mathchoice{}{}%
 {\,\,\botsmash{\hbox{\lower4\ex@\hbox{$\m@th\widehat{\null}$}}}}%
 {\,\botsmash{\hbox{\lower3\ex@\hbox{$\m@th\hat{\null}$}}}}}}

\def\spacute{^{\!\botsmash{\hbox{\lower\@ne ex\hbox{\'{}}}}}}
\def\spgrave{^{\mathchoice{}{}{}{\!}%
 \botsmash{\hbox{\lower\@ne ex\hbox{\`{}}}}}}
\def\spdot{^{\hbox{\raise\ex@\hbox{\rm.}}}}
\def\spddot{^{\hbox{\raise\ex@\hbox{\rm..}}}}
\def\spdddot{^{\hbox{\raise\ex@\hbox{\rm...}}}}
\def\spddddot{^{\hbox{\raise\ex@\hbox{\rm....}}}}
\def\spbreve{^{\!\botsmash{\hbox{\lower4\ex@\hbox{\u{}}}}}}

\message{\string\text,}
\def\textonlyfont@#1#2{\def#1{\RIfM@
 \Err@{Use \string#1\space only in text}\else#2\fi}}
\textonlyfont@\rm\tenrm
\textonlyfont@\it\tenit
\textonlyfont@\sl\tensl
\textonlyfont@\bf\tenbf
\def\oldnos#1{\RIfM@{\mathcode`\,="013B \fam\@ne#1}\else
 \leavevmode\hbox{$\m@th\mathcode`\,="013B \fam\@ne#1$}\fi}
\def\text{\RIfM@\expandafter\text@\else\expandafter\text@@\fi}
\def\text@@#1{\leavevmode\hbox{#1}}
\def\mathhexbox@#1#2#3{\text{$\m@th\mathchar"#1#2#3$}}
\def\dag{{\mathhexbox@279}}
\def\ddag{{\mathhexbox@27A}}
\def\S{{\mathhexbox@278}}
\def\P{{\mathhexbox@27B}}
\newif\iffirstchoice@
\firstchoice@true
\def\text@#1{\mathchoice
 {\hbox{\everymath{\displaystyle}\def\textfonti{\the\textfont\@ne}%
  \def\textfontii{\the\textfont\tw@}\textdef@@ T#1}}
 {\hbox{\firstchoice@false
  \everymath{\textstyle}\def\textfonti{\the\textfont\@ne}%
  \def\textfontii{\the\textfont\tw@}\textdef@@ T#1}}
 {\hbox{\firstchoice@false
  \everymath{\scriptstyle}\def\textfonti{\the\scriptfont\@ne}%
  \def\textfontii{\the\scriptfont\tw@}\textdef@@ S\rm#1}}
 {\hbox{\firstchoice@false
  \everymath{\scriptscriptstyle}\def\textfonti
  {\the\scriptscriptfont\@ne}%
  \def\textfontii{\the\scriptscriptfont\tw@}\textdef@@ s\rm#1}}}
\def\textdef@@#1{\textdef@#1\rm\textdef@#1\bf\textdef@#1\sl\textdef@#1\it}
\def\rmfam{0}
\def\textdef@#1#2{%
 \DN@{\csname\expandafter\eat@\string#2fam\endcsname}%
 \if S#1\edef#2{\the\scriptfont\next@\relax}%
 \else\if s#1\edef#2{\the\scriptscriptfont\next@\relax}%
 \else\edef#2{\the\textfont\next@\relax}\fi\fi}
\scriptfont\itfam\tenit \scriptscriptfont\itfam\tenit
\scriptfont\slfam\tensl \scriptscriptfont\slfam\tensl
\newif\iftopfolded@
\newif\ifbotfolded@
\def\topfoldedtext{\topfolded@true\botfolded@false\foldedtext@}
\def\botfoldedtext{\botfolded@true\topfolded@false\foldedtext@}
\def\foldedtext{\topfolded@false\botfolded@false\foldedtext@}
\Invalid@\foldedwidth
\def\foldedtext@{\relaxnext@
 \DN@{\ifx\next\foldedwidth\let\next@\nextii@\else
  \DN@{\nextii@\foldedwidth{.3\hsize}}\fi\next@}%
 \DNii@\foldedwidth##1##2{\setbox\z@\vbox
  {\normalbaselines\hsize##1\relax
  \tolerance1600 \noindent\ignorespaces##2}\ifbotfolded@\boxz@\else
  \iftopfolded@\vtop{\unvbox\z@}\else\vcenter{\boxz@}\fi\fi}%
 \FN@\next@}
\message{math font commands,}
\def\bold{\RIfM@\expandafter\bold@\else
 \expandafter\nonmatherr@\expandafter\bold\fi}
\def\bold@#1{{\bold@@{#1}}}
\def\bold@@#1{\fam\bffam\relax#1}
\def\slanted{\RIfM@\expandafter\slanted@\else
 \expandafter\nonmatherr@\expandafter\slanted\fi}
\def\slanted@#1{{\slanted@@{#1}}}
\def\slanted@@#1{\fam\slfam\relax#1}
\def\roman{\RIfM@\expandafter\roman@\else
 \expandafter\nonmatherr@\expandafter\roman\fi}
\def\roman@#1{{\roman@@{#1}}}
\def\roman@@#1{\fam\rmfam\relax#1}
\def\italic{\RIfM@\expandafter\italic@\else
 \expandafter\nonmatherr@\expandafter\italic\fi}
\def\italic@#1{{\italic@@{#1}}}
\def\italic@@#1{\fam\itfam\relax#1}
\def\Cal{\RIfM@\expandafter\Cal@\else
 \expandafter\nonmatherr@\expandafter\Cal\fi}
\def\Cal@#1{{\Cal@@{#1}}}
\def\Cal@@#1{\noaccents@\fam\tw@#1}
\mathchardef\Gamma="0000
\mathchardef\Delta="0001
\mathchardef\Theta="0002
\mathchardef\Lambda="0003
\mathchardef\Xi="0004
\mathchardef\Pi="0005
\mathchardef\Sigma="0006
\mathchardef\Upsilon="0007
\mathchardef\Phi="0008
\mathchardef\Psi="0009
\mathchardef\Omega="000A
\mathchardef\varGamma="0100
\mathchardef\varDelta="0101
\mathchardef\varTheta="0102
\mathchardef\varLambda="0103
\mathchardef\varXi="0104
\mathchardef\varPi="0105
\mathchardef\varSigma="0106
\mathchardef\varUpsilon="0107
\mathchardef\varPhi="0108
\mathchardef\varPsi="0109
\mathchardef\varOmega="010A
\let\alloc@@\alloc@
\def\hexnumber@#1{\ifcase#1 0\or 1\or 2\or 3\or 4\or 5\or 6\or 7\or 8\or
 9\or A\or B\or C\or D\or E\or F\fi}
\def\loadmsam{%
 \font@\tenmsa=msam10
 \font@\sevenmsa=msam7
 \font@\fivemsa=msam5
 \alloc@@8\fam\chardef\sixt@@n\msafam
 \textfont\msafam=\tenmsa
 \scriptfont\msafam=\sevenmsa
 \scriptscriptfont\msafam=\fivemsa
 \edef\next{\hexnumber@\msafam}%
 \mathchardef\dabar@"0\next39
 \edef\dashrightarrow{\mathrel{\dabar@\dabar@\mathchar"0\next4B}}%
 \edef\dashleftarrow{\mathrel{\mathchar"0\next4C\dabar@\dabar@}}%
 \let\dasharrow\dashrightarrow
 \edef\ulcorner{\delimiter"4\next70\next70 }%
 \edef\urcorner{\delimiter"5\next71\next71 }%
 \edef\llcorner{\delimiter"4\next78\next78 }%
 \edef\lrcorner{\delimiter"5\next79\next79 }%
 \edef\yen{{\noexpand\mathhexbox@\next55}}%
 \edef\checkmark{{\noexpand\mathhexbox@\next58}}%
 \edef\circledR{{\noexpand\mathhexbox@\next72}}%
 \edef\maltese{{\noexpand\mathhexbox@\next7A}}%
 \global\let\loadmsam\empty}%
\def\loadmsbm{%
 \font@\tenmsb=msbm10 \font@\sevenmsb=msbm7 \font@\fivemsb=msbm5
 \alloc@@8\fam\chardef\sixt@@n\msbfam
 \textfont\msbfam=\tenmsb
 \scriptfont\msbfam=\sevenmsb \scriptscriptfont\msbfam=\fivemsb
 \global\let\loadmsbm\empty
 }
\def\widehat#1{\ifx\undefined\msbfam \DN@{362}%
  \else \setboxz@h{$\m@th#1$}%
    \edef\next@{\ifdim\wdz@>\tw@ em%
        \hexnumber@\msbfam 5B%
      \else 362\fi}\fi
  \mathaccent"0\next@{#1}}
\def\widetilde#1{\ifx\undefined\msbfam \DN@{365}%
  \else \setboxz@h{$\m@th#1$}%
    \edef\next@{\ifdim\wdz@>\tw@ em%
        \hexnumber@\msbfam 5D%
      \else 365\fi}\fi
  \mathaccent"0\next@{#1}}
\message{\string\newsymbol,}
\def\newsymbol#1#2#3#4#5{\define#1{}%
  \count@#2\relax \advance\count@\m@ne 
 \ifcase\count@
   \ifx\undefined\msafam\loadmsam\fi \let\next@\msafam
 \or \ifx\undefined\msbfam\loadmsbm\fi \let\next@\msbfam
 \else  \Err@{\Invalid@@\string\newsymbol}\let\next@\tw@\fi
 \mathchardef#1="#3\hexnumber@\next@#4#5\space}

\def\Bbb{\RIfM@\expandafter\Bbb@\else
 \expandafter\nonmatherr@\expandafter\Bbb\fi}
\def\Bbb@#1{{\Bbb@@{#1}}}
\def\Bbb@@#1{\noaccents@\fam\msbfam\relax#1}
\message{bold Greek and bold symbols,}
\def\loadbold{%
 \font@\tencmmib=cmmib10 \font@\sevencmmib=cmmib7 \font@\fivecmmib=cmmib5
 \skewchar\tencmmib'177 \skewchar\sevencmmib'177 \skewchar\fivecmmib'177
 \alloc@@8\fam\chardef\sixt@@n\cmmibfam
 \textfont\cmmibfam\tencmmib
 \scriptfont\cmmibfam\sevencmmib \scriptscriptfont\cmmibfam\fivecmmib
 \font@\tencmbsy=cmbsy10 \font@\sevencmbsy=cmbsy7 \font@\fivecmbsy=cmbsy5
 \skewchar\tencmbsy'60 \skewchar\sevencmbsy'60 \skewchar\fivecmbsy'60
 \alloc@@8\fam\chardef\sixt@@n\cmbsyfam
 \textfont\cmbsyfam\tencmbsy
 \scriptfont\cmbsyfam\sevencmbsy \scriptscriptfont\cmbsyfam\fivecmbsy
 \let\loadbold\empty
}
\def\boldnotloaded#1{\Err@{\ifcase#1\or First\else Second\fi
       bold symbol font not loaded}}
\def\mathchari@#1#2#3{\ifx\undefined\cmmibfam
    \boldnotloaded@\@ne
  \else\mathchar"#1\hexnumber@\cmmibfam#2#3\space \fi}
\def\mathcharii@#1#2#3{\ifx\undefined\cmbsyfam
    \boldnotloaded\tw@
  \else \mathchar"#1\hexnumber@\cmbsyfam#2#3\space\fi}
\edef\bffam@{\hexnumber@\bffam}
\def\boldkey#1{\ifcat\noexpand#1A%
  \ifx\undefined\cmmibfam \boldnotloaded\@ne
  \else {\fam\cmmibfam#1}\fi
 \else
 \ifx#1!\mathchar"5\bffam@21 \else
 \ifx#1(\mathchar"4\bffam@28 \else\ifx#1)\mathchar"5\bffam@29 \else
 \ifx#1+\mathchar"2\bffam@2B \else\ifx#1:\mathchar"3\bffam@3A \else
 \ifx#1;\mathchar"6\bffam@3B \else\ifx#1=\mathchar"3\bffam@3D \else
 \ifx#1?\mathchar"5\bffam@3F \else\ifx#1[\mathchar"4\bffam@5B \else
 \ifx#1]\mathchar"5\bffam@5D \else
 \ifx#1,\mathchari@63B \else
 \ifx#1-\mathcharii@200 \else
 \ifx#1.\mathchari@03A \else
 \ifx#1/\mathchari@03D \else
 \ifx#1<\mathchari@33C \else
 \ifx#1>\mathchari@33E \else
 \ifx#1*\mathcharii@203 \else
 \ifx#1|\mathcharii@06A \else
 \ifx#10\bold0\else\ifx#11\bold1\else\ifx#12\bold2\else\ifx#13\bold3\else
 \ifx#14\bold4\else\ifx#15\bold5\else\ifx#16\bold6\else\ifx#17\bold7\else
 \ifx#18\bold8\else\ifx#19\bold9\else
  \Err@{\string\boldkey\space can't be used with #1}%
 \fi\fi\fi\fi\fi\fi\fi\fi\fi\fi\fi\fi\fi\fi\fi
 \fi\fi\fi\fi\fi\fi\fi\fi\fi\fi\fi\fi\fi\fi}
\def\boldsymbol#1{%
 \DN@{\Err@{You can't use \string\boldsymbol\space with \string#1}#1}%
 \ifcat\noexpand#1A%
   \let\next@\relax
   \ifx\undefined\cmmibfam \boldnotloaded\@ne
   \else {\fam\cmmibfam#1}\fi
 \else
  \xdef\meaning@{\meaning#1.........}%
  \expandafter\math@\meaning@\math@
  \ifmath@
   \expandafter\mathch@\meaning@\mathch@
   \ifmathch@
    \expandafter\boldsymbol@@\meaning@\boldsymbol@@
   \fi
  \else
   \expandafter\macro@\meaning@\macro@
   \expandafter\delim@\meaning@\delim@
   \ifdelim@
    \expandafter\delim@@\meaning@\delim@@
   \else
    \boldsymbol@{#1}%
   \fi
  \fi
 \fi
 \next@}
\def\mathhexboxii@#1#2{\ifx\undefined\cmbsyfam
    \boldnotloaded\tw@
  \else \mathhexbox@{\hexnumber@\cmbsyfam}{#1}{#2}\fi}
\def\boldsymbol@#1{\let\next@\relax\let\next#1%
 \ifx\next\cdot\mathcharii@201 \else
 \ifx\next\prime{{\null\mathcharii@030 \null}}\else
 \ifx\next\lbrack\mathchar"4\bffam@5B \else
 \ifx\next\rbrack\mathchar"5\bffam@5D \else
 \ifx\next\{\mathcharii@466 \else
 \ifx\next\lbrace\mathcharii@466 \else
 \ifx\next\}\mathcharii@567 \else
 \ifx\next\rbrace\mathcharii@567 \else
 \ifx\next\surd{{\mathcharii@170}}\else
 \ifx\next\S{{\mathhexboxii@78}}\else
 \ifx\next\P{{\mathhexboxii@7B}}\else
 \ifx\next\dag{{\mathhexboxii@79}}\else
 \ifx\next\ddag{{\mathhexboxii@7A}}\else
 \DN@{\Err@{You can't use \string\boldsymbol\space with \string#1}#1}%
 \fi\fi\fi\fi\fi\fi\fi\fi\fi\fi\fi\fi\fi}
\def\boldsymbol@@#1.#2\boldsymbol@@{\classnum@#1 \count@@@\classnum@        
 \divide\classnum@4096 \count@\classnum@                                    
 \multiply\count@4096 \advance\count@@@-\count@ \count@@\count@@@           
 \divide\count@@@\@cclvi \count@\count@@                                    
 \multiply\count@@@\@cclvi \advance\count@@-\count@@@                       
 \divide\count@@@\@cclvi                                                    
 \multiply\classnum@4096 \advance\classnum@\count@@                         
 \ifnum\count@@@=\z@                                                        
  \count@"\bffam@ \multiply\count@\@cclvi
  \advance\classnum@\count@
  \DN@{\mathchar\number\classnum@}%
 \else
  \ifnum\count@@@=\@ne                                                      
   \ifx\undefined\cmmibfam \DN@{\boldnotloaded\@ne}%
   \else \count@\cmmibfam \multiply\count@\@cclvi
     \advance\classnum@\count@
     \DN@{\mathchar\number\classnum@}\fi
  \else
   \ifnum\count@@@=\tw@                                                    
     \ifx\undefined\cmbsyfam
       \DN@{\boldnotloaded\tw@}%
     \else
       \count@\cmbsyfam \multiply\count@\@cclvi
       \advance\classnum@\count@
       \DN@{\mathchar\number\classnum@}%
     \fi
  \fi
 \fi
\fi}
\newif\ifdelim@
\newcount\delimcount@
{\uccode`6=`\\ \uccode`7=`d \uccode`8=`e \uccode`9=`l
 \uppercase{\gdef\delim@#1#2#3#4#5\delim@
  {\delim@false\ifx 6#1\ifx 7#2\ifx 8#3\ifx 9#4\delim@true
   \xdef\meaning@{#5}\fi\fi\fi\fi}}}
\def\delim@@#1"#2#3#4#5#6\delim@@{\if#32%
\let\next@\relax
 \ifx\undefined\cmbsyfam \boldnotloaded\@ne
 \else \mathcharii@#2#4#5\space \fi\fi}
\def\vert{\delimiter"026A30C }
\def\Vert{\delimiter"026B30D }
\let\|\Vert
\def\backslash{\delimiter"026E30F }
\def\boldkeydots@#1{\bold@true\let\next=#1\let\delayed@=#1\mdots@@
 \boldkey#1\bold@false}  
\def\boldsymboldots@#1{\bold@true\let\next#1\let\delayed@#1\mdots@@
 \boldsymbol#1\bold@false}
\message{Euler fonts,}

\def\frak{\mathfont@\frak}

\def\loadmathfont#1{%
   \expandafter\font@\csname ten#1\endcsname=#110
   \expandafter\font@\csname seven#1\endcsname=#17
   \expandafter\font@\csname five#1\endcsname=#15
   \edef\next{\noexpand\alloc@@8\fam\chardef\sixt@@n
     \expandafter\noexpand\csname#1fam\endcsname}%
   \next
   \textfont\csname#1fam\endcsname \csname ten#1\endcsname
   \scriptfont\csname#1fam\endcsname \csname seven#1\endcsname
   \scriptscriptfont\csname#1fam\endcsname \csname five#1\endcsname
   \expandafter\def\csname #1\expandafter\endcsname\expandafter{%
      \expandafter\mathfont@\csname#1\endcsname}%
 \expandafter\gdef\csname load#1\endcsname{}%
}
\def\mathfont@#1{\RIfM@\expandafter\mathfont@@\expandafter#1\else
  \expandafter\nonmatherr@\expandafter#1\fi}
\def\mathfont@@#1#2{{\mathfont@@@#1{#2}}}
\def\mathfont@@@#1#2{\noaccents@
   \fam\csname\expandafter\eat@\string#1fam\endcsname
   \relax#2}
\message{math accents,}
\def\accentclass@{7}
\def\noaccents@{\def\accentclass@{0}}
\def\makeacc@#1#2{\def#1{\mathaccent"\accentclass@#2 }}
\makeacc@\hat{05E}
\makeacc@\check{014}
\makeacc@\tilde{07E}
\makeacc@\acute{013}
\makeacc@\grave{012}
\makeacc@\dot{05F}
\makeacc@\ddot{07F}
\makeacc@\breve{015}
\makeacc@\bar{016}

\newcount\skewcharcount@
\newcount\familycount@
\def\theskewchar@{\familycount@\@ne
 \global\skewcharcount@\the\skewchar\textfont\@ne                           
 \ifnum\fam>\m@ne\ifnum\fam<16
  \global\familycount@\the\fam\relax
  \global\skewcharcount@\the\skewchar\textfont\the\fam\relax\fi\fi          
 \ifnum\skewcharcount@>\m@ne
  \ifnum\skewcharcount@<128
  \multiply\familycount@256
  \global\advance\skewcharcount@\familycount@
  \global\advance\skewcharcount@28672
  \mathchar\skewcharcount@\else
  \global\skewcharcount@\m@ne\fi\else
 \global\skewcharcount@\m@ne\fi}                                            
\newcount\pointcount@
\def\getpoints@#1.#2\getpoints@{\pointcount@#1 }
\newdimen\accentdimen@
\newcount\accentmu@
\def\dimentomu@{\multiply\accentdimen@ 100
 \expandafter\getpoints@\the\accentdimen@\getpoints@
 \multiply\pointcount@18
 \divide\pointcount@\@m
 \global\accentmu@\pointcount@}
\def\Makeacc@#1#2{\def#1{\RIfM@\DN@{\mathaccent@
 {"\accentclass@#2 }}\else\DN@{\nonmatherr@{#1}}\fi\next@}}
\def\unbracefonts@{\let\Cal@\Cal@@\let\roman@\roman@@\let\bold@\bold@@
 \let\slanted@\slanted@@}
\def\mathaccent@#1#2{\ifnum\fam=\m@ne\xdef\thefam@{1}\else
 \xdef\thefam@{\the\fam}\fi                                                 
 \accentdimen@\z@                                                           
 \setboxz@h{\unbracefonts@$\m@th\fam\thefam@\relax#2$}
 \ifdim\accentdimen@=\z@\DN@{\mathaccent#1{#2}}
  \setbox@ne\hbox{\unbracefonts@$\m@th\fam\thefam@\relax#2\theskewchar@$}
  \setbox\tw@\hbox{$\m@th\ifnum\skewcharcount@=\m@ne\else
   \mathchar\skewcharcount@\fi$}
  \global\accentdimen@\wd@ne\global\advance\accentdimen@-\wdz@
  \global\advance\accentdimen@-\wd\tw@                                     
  \global\multiply\accentdimen@\tw@
  \dimentomu@\global\advance\accentmu@\@ne                                 
 \else\DN@{{\mathaccent#1{#2\mkern\accentmu@ mu}%
    \mkern-\accentmu@ mu}{}}\fi                                             
 \next@}\Makeacc@\Hat{05E}
\Makeacc@\Check{014}
\Makeacc@\Tilde{07E}
\Makeacc@\Acute{013}
\Makeacc@\Grave{012}
\Makeacc@\Dot{05F}
\Makeacc@\Ddot{07F}
\Makeacc@\Breve{015}
\Makeacc@\Bar{016}
\def\Vec{\RIfM@\DN@{\mathaccent@{"017E }}\else
 \DN@{\nonmatherr@\Vec}\fi\next@}
\def\accentedsymbol#1#2{\csname newbox\expandafter\endcsname
  \csname\expandafter\eat@\string#1@box\endcsname
 \expandafter\setbox\csname\expandafter\eat@
  \string#1@box\endcsname\hbox{$\m@th#2$}\define
  #1{\copy\csname\expandafter\eat@\string#1@box\endcsname{}}}
\message{roots,}
\def\sqrt#1{\radical"270370 {#1}}
\let\underline@\underline
\let\overline@\overline
\def\underline#1{\underline@{#1}}
\def\overline#1{\overline@{#1}}
\Invalid@\leftroot
\Invalid@\uproot
\newcount\uproot@
\newcount\leftroot@
\def\root{\relaxnext@
  \DN@{\ifx\next\uproot\let\next@\nextii@\else
   \ifx\next\leftroot\let\next@\nextiii@\else
   \let\next@\plainroot@\fi\fi\next@}%
  \DNii@\uproot##1{\uproot@##1\relax\FN@\nextiv@}%
  \def\nextiv@{\ifx\next\space@\DN@. {\FN@\nextv@}\else
   \DN@.{\FN@\nextv@}\fi\next@.}%
  \def\nextv@{\ifx\next\leftroot\let\next@\nextvi@\else
   \let\next@\plainroot@\fi\next@}%
  \def\nextvi@\leftroot##1{\leftroot@##1\relax\plainroot@}%
   \def\nextiii@\leftroot##1{\leftroot@##1\relax\FN@\nextvii@}%
  \def\nextvii@{\ifx\next\space@
   \DN@. {\FN@\nextviii@}\else
   \DN@.{\FN@\nextviii@}\fi\next@.}%
  \def\nextviii@{\ifx\next\uproot\let\next@\nextix@\else
   \let\next@\plainroot@\fi\next@}%
  \def\nextix@\uproot##1{\uproot@##1\relax\plainroot@}%
  \bgroup\uproot@\z@\leftroot@\z@\FN@\next@}
\def\plainroot@#1\of#2{\setbox\rootbox\hbox{$\m@th\scriptscriptstyle{#1}$}%
 \mathchoice{\r@@t\displaystyle{#2}}{\r@@t\textstyle{#2}}
 {\r@@t\scriptstyle{#2}}{\r@@t\scriptscriptstyle{#2}}\egroup}
\def\r@@t#1#2{\setboxz@h{$\m@th#1\sqrt{#2}$}%
 \dimen@\ht\z@\advance\dimen@-\dp\z@
 \setbox@ne\hbox{$\m@th#1\mskip\uproot@ mu$}\advance\dimen@ 1.667\wd@ne
 \mkern-\leftroot@ mu\mkern5mu\raise.6\dimen@\copy\rootbox
 \mkern-10mu\mkern\leftroot@ mu\boxz@}
\def\boxed#1{\setboxz@h{$\m@th\displaystyle{#1}$}\dimen@.4\ex@
 \advance\dimen@3\ex@\advance\dimen@\dp\z@
 \hbox{\lower\dimen@\hbox{%
 \vbox{\hrule height.4\ex@
 \hbox{\vrule width.4\ex@\hskip3\ex@\vbox{\vskip3\ex@\boxz@\vskip3\ex@}%
 \hskip3\ex@\vrule width.4\ex@}\hrule height.4\ex@}%
 }}}
\message{commutative diagrams,}
\let\ampersand@\relax
\newdimen\minaw@
\minaw@11.11128\ex@
\newdimen\minCDaw@
\minCDaw@2.5pc
\def\minCDarrowwidth#1{\RIfMIfI@\onlydmatherr@\minCDarrowwidth
 \else\minCDaw@#1\relax\fi\else\onlydmatherr@\minCDarrowwidth\fi}
\newif\ifCD@
\def\CD{\bgroup\vspace@\relax\let\ampersand@&\iffalse}\fi
 \CD@true\vcenter\bgroup\Let@\tabskip\z@skip\baselineskip20\ex@
 \lineskip3\ex@\lineskiplimit3\ex@\halign\bgroup
 &\hfill$\m@th##$\hfill\crcr}
\def\endCD{\crcr\egroup\egroup\egroup}
\newdimen\bigaw@
\atdef@>#1>#2>{\ampersand@                                                  
 \setboxz@h{$\m@th\ssize\;{#1}\;\;$}
 \setbox@ne\hbox{$\m@th\ssize\;{#2}\;\;$}
 \setbox\tw@\hbox{$\m@th#2$}
 \ifCD@\global\bigaw@\minCDaw@\else\global\bigaw@\minaw@\fi                 
 \ifdim\wdz@>\bigaw@\global\bigaw@\wdz@\fi
 \ifdim\wd@ne>\bigaw@\global\bigaw@\wd@ne\fi                                
 \ifCD@\enskip\fi                                                           
 \ifdim\wd\tw@>\z@
  \mathrel{\mathop{\hbox to\bigaw@{\rightarrowfill@\displaystyle}}%
    \limits^{#1}_{#2}}
 \else\mathrel{\mathop{\hbox to\bigaw@{\rightarrowfill@\displaystyle}}%
    \limits^{#1}}\fi                                                        
 \ifCD@\enskip\fi                                                          
 \ampersand@}                                                              
\atdef@<#1<#2<{\ampersand@\setboxz@h{$\m@th\ssize\;\;{#1}\;$}%
 \setbox@ne\hbox{$\m@th\ssize\;\;{#2}\;$}\setbox\tw@\hbox{$\m@th#2$}%
 \ifCD@\global\bigaw@\minCDaw@\else\global\bigaw@\minaw@\fi
 \ifdim\wdz@>\bigaw@\global\bigaw@\wdz@\fi
 \ifdim\wd@ne>\bigaw@\global\bigaw@\wd@ne\fi
 \ifCD@\enskip\fi
 \ifdim\wd\tw@>\z@
  \mathrel{\mathop{\hbox to\bigaw@{\leftarrowfill@\displaystyle}}%
       \limits^{#1}_{#2}}\else
  \mathrel{\mathop{\hbox to\bigaw@{\leftarrowfill@\displaystyle}}%
       \limits^{#1}}\fi
 \ifCD@\enskip\fi\ampersand@}
\begingroup
 \catcode`\~=\active \lccode`\~=`\@
 \lowercase{%
  \global\atdef@)#1)#2){~>#1>#2>}
  \global\atdef@(#1(#2({~<#1<#2<}}
\endgroup
\atdef@ A#1A#2A{\llap{$\m@th\vcenter{\hbox
 {$\ssize#1$}}$}\Big\uparrow\rlap{$\m@th\vcenter{\hbox{$\ssize#2$}}$}&&}
\atdef@ V#1V#2V{\llap{$\m@th\vcenter{\hbox
 {$\ssize#1$}}$}\Big\downarrow\rlap{$\m@th\vcenter{\hbox{$\ssize#2$}}$}&&}
\atdef@={&\enskip\mathrel
 {\vbox{\hrule width\minCDaw@\vskip3\ex@\hrule width
 \minCDaw@}}\enskip&}
\atdef@|{\Big\Vert&&}
\atdef@\vert{\Big\Vert&&}
\def\pretend#1\haswidth#2{\setboxz@h{$\m@th\scriptstyle{#2}$}\hbox
 to\wdz@{\hfill$\m@th\scriptstyle{#1}$\hfill}}
\message{poor man's bold,}
\def\pmb{\RIfM@\expandafter\mathpalette\expandafter\pmb@\else
 \expandafter\pmb@@\fi}
\def\pmb@@#1{\leavevmode\setboxz@h{#1}%
   \dimen@-\wdz@
   \kern-.5\ex@\copy\z@
   \kern\dimen@\kern.25\ex@\raise.4\ex@\copy\z@
   \kern\dimen@\kern.25\ex@\box\z@
}
\def\binrel@@#1{\ifdim\wd2<\z@\mathbin{#1}\else\ifdim\wd\tw@>\z@
 \mathrel{#1}\else{#1}\fi\fi}
\newdimen\pmbraise@
\def\pmb@#1#2{\setbox\thr@@\hbox{$\m@th#1{#2}$}%
 \setbox4\hbox{$\m@th#1\mkern.5mu$}\pmbraise@\wd4\relax
 \binrel@{#2}%
 \dimen@-\wd\thr@@
   \binrel@@{%
   \mkern-.8mu\copy\thr@@
   \kern\dimen@\mkern.4mu\raise\pmbraise@\copy\thr@@
   \kern\dimen@\mkern.4mu\box\thr@@
}}
\def\documentstyle#1{\W@{}\input #1.sty\relax}
\message{syntax check,}
\font\dummyft@=dummy
\fontdimen1 \dummyft@=\z@
\fontdimen2 \dummyft@=\z@
\fontdimen3 \dummyft@=\z@
\fontdimen4 \dummyft@=\z@
\fontdimen5 \dummyft@=\z@
\fontdimen6 \dummyft@=\z@
\fontdimen7 \dummyft@=\z@
\fontdimen8 \dummyft@=\z@
\fontdimen9 \dummyft@=\z@
\fontdimen10 \dummyft@=\z@
\fontdimen11 \dummyft@=\z@
\fontdimen12 \dummyft@=\z@
\fontdimen13 \dummyft@=\z@
\fontdimen14 \dummyft@=\z@
\fontdimen15 \dummyft@=\z@
\fontdimen16 \dummyft@=\z@
\fontdimen17 \dummyft@=\z@
\fontdimen18 \dummyft@=\z@
\fontdimen19 \dummyft@=\z@
\fontdimen20 \dummyft@=\z@
\fontdimen21 \dummyft@=\z@
\fontdimen22 \dummyft@=\z@
\def\fontlist@{\\{\tenrm}\\{\sevenrm}\\{\fiverm}\\{\teni}\\{\seveni}%
 \\{\fivei}\\{\tensy}\\{\sevensy}\\{\fivesy}\\{\tenex}\\{\tenbf}\\{\sevenbf}%
 \\{\fivebf}\\{\tensl}\\{\tenit}}
\def\font@#1=#2 {\rightappend@#1\to\fontlist@\font#1=#2 }
\def\dodummy@{{\def\\##1{\global\let##1\dummyft@}\fontlist@}}
\def\nopages@{\output{\setbox\z@\box\@cclv \deadcycles\z@}%
 \alloc@5\toks\toksdef\@cclvi\output}
\let\galleys\nopages@
\newif\ifsyntax@
\newcount\countxviii@
\def\syntax{\syntax@true\dodummy@\countxviii@\count18
 \loop\ifnum\countxviii@>\m@ne\textfont\countxviii@=\dummyft@
 \scriptfont\countxviii@=\dummyft@\scriptscriptfont\countxviii@=\dummyft@
 \advance\countxviii@\m@ne\repeat                                           
 \dummyft@\tracinglostchars\z@\nopages@\frenchspacing\hbadness\@M}
\def\first@#1#2\end{#1}
\def\printoptions{\W@{Do you want S(yntax check),
  G(alleys) or P(ages)?}%
 \message{Type S, G or P, followed by <return>: }%
 \begingroup 
 \endlinechar\m@ne 
 \read\m@ne to\ans@
 \edef\ans@{\uppercase{\def\noexpand\ans@{%
   \expandafter\first@\ans@ P\end}}}%
 \expandafter\endgroup\ans@
 \if\ans@ P
 \else \if\ans@ S\syntax
 \else \if\ans@ G\galleys
 \else\message{? Unknown option: \ans@; using the `pages' option.}%
 \fi\fi\fi}
\def\alloc@#1#2#3#4#5{\global\advance\count1#1by\@ne
 \ch@ck#1#4#2\allocationnumber=\count1#1
 \global#3#5=\allocationnumber
 \ifalloc@\wlog{\string#5=\string#2\the\allocationnumber}\fi}
\def\document{\def\alloclist@{}\def\fontlist@{}}
\let\enddocument\bye

\let\proclaim\undefined
\let\footnote\undefined
\let\=\undefined
\let\>\undefined

\catcode`\@=\active
\message{... finished}

\expandafter\ifx\csname mathdefs.tex\endcsname\relax
  \expandafter\gdef\csname mathdefs.tex\endcsname{}
\else \message{Hey!  Apparently you were trying to
  \string\input{mathdefs.tex} twice.   This does not make sense.} 
\errmessage{Please edit your file (probably \jobname.tex) and remove
any duplicate ``\string\input'' lines}\endinput\fi




\catcode`\X=12\catcode`\@=11

\def\n@wcount{\alloc@0\count\countdef\insc@unt}
\def\n@wwrite{\alloc@7\write\chardef\sixt@@n}
\def\n@wread{\alloc@6\read\chardef\sixt@@n}
\def\r@s@t{\relax}\def\v@idline{\par}\def\@mputate#1/{#1}
\def\l@c@l#1X{\firstpart.#1}\def\gl@b@l#1X{#1}\def\t@d@l#1X{{}}

\def\crossrefs#1{\ifx\all#1\let\tr@ce=\all\else\def\tr@ce{#1,}\fi
   \n@wwrite\cit@tionsout\openout\cit@tionsout=\jobname.cit 
   \write\cit@tionsout{\tr@ce}\expandafter\setfl@gs\tr@ce,}
\def\setfl@gs#1,{\def\@{#1}\ifx\@\empty\let\next=\relax
   \else\let\next=\setfl@gs\expandafter\xdef
   \csname#1tr@cetrue\endcsname{}\fi\next}
\def\m@ketag#1#2{\expandafter\n@wcount\csname#2tagno\endcsname
     \csname#2tagno\endcsname=0\let\tail=\all\xdef\all{\tail#2,}
   \ifx#1\l@c@l\let\tail=\r@s@t\xdef\r@s@t{\csname#2tagno\endcsname=0\tail}\fi
   \expandafter\gdef\csname#2cite\endcsname##1{\expandafter
     \ifx\csname#2tag##1\endcsname\relax?\else\csname#2tag##1\endcsname\fi
     \expandafter\ifx\csname#2tr@cetrue\endcsname\relax\else
     \write\cit@tionsout{#2tag ##1 cited on page \folio.}\fi}
   \expandafter\gdef\csname#2page\endcsname##1{\expandafter
     \ifx\csname#2page##1\endcsname\relax?\else\csname#2page##1\endcsname\fi
     \expandafter\ifx\csname#2tr@cetrue\endcsname\relax\else
     \write\cit@tionsout{#2tag ##1 cited on page \folio.}\fi}
   \expandafter\gdef\csname#2tag\endcsname##1{\expandafter
      \ifx\csname#2check##1\endcsname\relax
      \expandafter\xdef\csname#2check##1\endcsname{}%
      \else\immediate\write16{Warning: #2tag ##1 used more than once.}\fi
      \multit@g{#1}{#2}##1/X%
      \write\t@gsout{#2tag ##1 assigned number \csname#2tag##1\endcsname\space
      on page \number\count0.}%
   \csname#2tag##1\endcsname}}

\def\multit@g#1#2#3/#4X{\def\t@mp{#4}\ifx\t@mp\empty%
      \global\advance\csname#2tagno\endcsname by 1 
      \expandafter\xdef\csname#2tag#3\endcsname
      {#1\number\csname#2tagno\endcsnameX}%
   \else\expandafter\ifx\csname#2last#3\endcsname\relax
      \expandafter\n@wcount\csname#2last#3\endcsname
      \global\advance\csname#2tagno\endcsname by 1 
      \expandafter\xdef\csname#2tag#3\endcsname
      {#1\number\csname#2tagno\endcsnameX}
      \write\t@gsout{#2tag #3 assigned number \csname#2tag#3\endcsname\space
      on page \number\count0.}\fi
   \global\advance\csname#2last#3\endcsname by 1
   \def\t@mp{\expandafter\xdef\csname#2tag#3/}%
   \expandafter\t@mp\@mputate#4\endcsname
   {\csname#2tag#3\endcsname\lastpart{\csname#2last#3\endcsname}}\fi}
\def\t@gs#1{\def\all{}\m@ketag#1e\m@ketag#1s\m@ketag\t@d@l p
\let\realscite\scite
\let\realstag\stag
   \m@ketag\gl@b@l r \n@wread\t@gsin
   \openin\t@gsin=\jobname.tgs \re@der \closein\t@gsin
   \n@wwrite\t@gsout\openout\t@gsout=\jobname.tgs }
\outer\def\localtags{\t@gs\l@c@l}
\outer\def\globaltags{\t@gs\gl@b@l}
\outer\def\newlocaltag#1{\m@ketag\l@c@l{#1}}
\outer\def\newglobaltag#1{\m@ketag\gl@b@l{#1}}

\newif\ifpr@ 
\def\m@kecs #1tag #2 assigned number #3 on page #4.%
   {\expandafter\gdef\csname#1tag#2\endcsname{#3}
   \expandafter\gdef\csname#1page#2\endcsname{#4}
   \ifpr@\expandafter\xdef\csname#1check#2\endcsname{}\fi}
\def\re@der{\ifeof\t@gsin\let\next=\relax\else
   \read\t@gsin to\t@gline\ifx\t@gline\v@idline\else
   \expandafter\m@kecs \t@gline\fi\let \next=\re@der\fi\next}
\def\pretags#1{\pr@true\pret@gs#1,,}
\def\pret@gs#1,{\def\@{#1}\ifx\@\empty\let\n@xtfile=\relax
   \else\let\n@xtfile=\pret@gs \openin\t@gsin=#1.tgs \message{#1} \re@der 
   \closein\t@gsin\fi \n@xtfile}

\newcount\sectno\sectno=0\newcount\subsectno\subsectno=0
\newif\ifultr@local \def\ultralocal{\ultr@localtrue}
\def\firstpart{\number\sectno}
\def\lastpart#1{\ifcase#1 \or a\or b\or c\or d\or e\or f\or g\or h\or 
   i\or k\or l\or m\or n\or o\or p\or q\or r\or s\or t\or u\or v\or w\or 
   x\or y\or z \fi}

\def\resetall{\global\advance\sectno by 1\subsectno=0
   \gdef\firstpart{\number\sectno}\r@s@t}
\def\resetsub{\global\advance\subsectno by 1
   \gdef\firstpart{\number\sectno.\number\subsectno}\r@s@t}
\def\newsection#1\par{\resetall\vskip0pt plus.3\vsize\penalty-250
   \vskip0pt plus-.3\vsize\bigskip\bigskip
   \message{#1}\leftline{\bf#1}\nobreak\bigskip}
\def\subsection#1\par{\ifultr@local\resetsub\fi
   \vskip0pt plus.2\vsize\penalty-250\vskip0pt plus-.2\vsize
   \bigskip\smallskip\message{#1}\leftline{\bf#1}\nobreak\medskip}


\newdimen\marginshift

\newdimen\margindelta
\newdimen\marginmax
\newdimen\marginmin

\def\margininit{       
\marginmax=3 true cm                  
				      
\margindelta=0.1 true cm              
\marginmin=0.1true cm                 
\marginshift=\marginmin
}    

\def\t@gsjj#1,{\def\@{#1}\ifx\@\empty\let\next=\relax\else\let\next=\t@gsjj
   \def\@@{p}\ifx\@\@@\else
   \expandafter\gdef\csname#1cite\endcsname##1{\citejj{##1}}
   \expandafter\gdef\csname#1page\endcsname##1{?}
   \expandafter\gdef\csname#1tag\endcsname##1{\tagjj{##1}}\fi\fi\next}
\newif\ifshowstuffinmargin
\showstuffinmarginfalse
\def\jjtags{\ifx\shlhetal\relax 
  \else
\ifx\shlhetal\undefinedcontrolseq
\else
\showstuffinmargintrue
\ifx\all\relax\else\expandafter\t@gsjj\all,\fi\fi \fi
}

\def\tagjj#1{\realstag{#1}\mginpar{\zeigen{#1}}}
\def\citejj#1{\rechnen{#1}\mginpar{\zeigen{#1}}}     

\def\rechnen#1{\expandafter\ifx\csname stag#1\endcsname\relax ??\else
                           \csname stag#1\endcsname\fi}

\newdimen\theight

\def\marginfont{\sevenrm}

\def\trymarginbox#1{\setbox0=\hbox{\marginfont\hskip\marginshift #1}%
		\global\marginshift\wd0 
		\global\advance\marginshift\margindelta}

\def \mginpar#1{%
\ifvmode\setbox0\hbox to \hsize{\hfill\rlap{\marginfont\quad#1}}%
\ht0 0cm
\dp0 0cm
\box0\vskip-\baselineskip
\else 
             \vadjust{\trymarginbox{#1}%
		\ifdim\marginshift>\marginmax \global\marginshift\marginmin
			\trymarginbox{#1}%
                \fi
             \theight=\ht0
             \advance\theight by \dp0    \advance\theight by \lineskip
             \kern -\theight \vbox to \theight{\rightline{\rlap{\box0}}%
\vss}}\fi}


\def\t@gsoff#1,{\def\@{#1}\ifx\@\empty\let\next=\relax\else\let\next=\t@gsoff
   \def\@@{p}\ifx\@\@@\else
   \expandafter\gdef\csname#1cite\endcsname##1{\zeigen{##1}}
   \expandafter\gdef\csname#1page\endcsname##1{?}
   \expandafter\gdef\csname#1tag\endcsname##1{\zeigen{##1}}\fi\fi\next}
\def\verbatimtags{\showstuffinmarginfalse
\ifx\all\relax\else\expandafter\t@gsoff\all,\fi}
\def\zeigen#1{\hbox{$\langle$}#1\hbox{$\rangle$}}

\def\margincite#1{\ifshowstuffinmargin\mginpar{\zeigen{#1}}\fi}

\def\margintag#1{\ifshowstuffinmargin\mginpar{\zeigen{#1}}\fi}

\def\(#1){\edef\dot@g{\ifmmode\ifinner(\hbox{\noexpand\etag{#1}})
   \else\noexpand\eqno(\hbox{\noexpand\etag{#1}})\fi
   \else(\noexpand\ecite{#1})\fi}\dot@g}

\newif\ifbr@ck
\def\eat#1{}
\def\[#1]{\br@cktrue[\br@cket#1'X]}
\def\br@cket#1'#2X{\def\temp{#2}\ifx\temp\empty\let\next\eat
   \else\let\next\br@cket\fi
   \ifbr@ck\br@ckfalse\br@ck@t#1,X\else\br@cktrue#1\fi\next#2X}
\def\br@ck@t#1,#2X{\def\temp{#2}\ifx\temp\empty\let\neext\eat
   \else\let\neext\br@ck@t\def\temp{,}\fi
   \def\teemp{#1}\ifx\teemp\empty\else\rcite{#1}\fi\temp\neext#2X}
\def\resetbr@cket{\gdef\[##1]{[\rtag{##1}]}}
\def\references{\resetbr@cket\newsection References\par}

\newtoks\symb@ls\newtoks\s@mb@ls\newtoks\p@gelist\n@wcount\ftn@mber
    \ftn@mber=1\newif\ifftn@mbers\ftn@mbersfalse\newif\ifbyp@ge\byp@gefalse
\def\defm@rk{\ifftn@mbers\n@mberm@rk\else\symb@lm@rk\fi}
\def\n@mberm@rk{\xdef\m@rk{{\the\ftn@mber}}%
    \global\advance\ftn@mber by 1 }
\def\rot@te#1{\let\temp=#1\global#1=\expandafter\r@t@te\the\temp,X}
\def\r@t@te#1,#2X{{#2#1}\xdef\m@rk{{#1}}}
\def\b@@st#1{{$^{#1}$}}\def\str@p#1{#1}
\def\symb@lm@rk{\ifbyp@ge\rot@te\p@gelist\ifnum\expandafter\str@p\m@rk=1 
    \s@mb@ls=\symb@ls\fi\write\f@nsout{\number\count0}\fi \rot@te\s@mb@ls}
\def\byp@ge{\byp@getrue\n@wwrite\f@nsin\openin\f@nsin=\jobname.fns 
    \n@wcount\currentp@ge\currentp@ge=0\p@gelist={0}
    \re@dfns\closein\f@nsin\rot@te\p@gelist
    \n@wread\f@nsout\openout\f@nsout=\jobname.fns }
\def\m@kelist#1X#2{{#1,#2}}
\def\re@dfns{\ifeof\f@nsin\let\next=\relax\else\read\f@nsin to \f@nline
    \ifx\f@nline\v@idline\else\let\t@mplist=\p@gelist
    \ifnum\currentp@ge=\f@nline
    \global\p@gelist=\expandafter\m@kelist\the\t@mplistX0
    \else\currentp@ge=\f@nline
    \global\p@gelist=\expandafter\m@kelist\the\t@mplistX1\fi\fi
    \let\next=\re@dfns\fi\next}
\def\symbols#1{\symb@ls={#1}\s@mb@ls=\symb@ls} 
\def\bigsymbol{\textstyle}
\symbols{\bigsymbol\ast,\dagger,\ddagger,\sharp,\flat,\natural,\star}
\def\ftnumbers{\ftn@mberstrue} \def\ftsymbols{\ftn@mbersfalse}
\def\paginal{\byp@ge} \def\resetftnumbers{\ftn@mber=1}
\def\ftnote#1{\defm@rk\expandafter\expandafter\expandafter\footnote
    \expandafter\b@@st\m@rk{#1}}

\long\def\jump#1\endjump{}
\def\ssum{\mathop{\lower .1em\hbox{$\textstyle\Sigma$}}\nolimits}

\def\qed{\nobreak\kern 1em \vrule height .5em width .5em depth 0em}
\def\newneq{\hbox{\rlap{\hbox to 1\wd9{\hss$=$\hss}}\raise .1em 
   \hbox to 1\wd9{\hss$\scriptscriptstyle/$\hss}}}
\def\subsetne{\setbox9 = \hbox{$\subset$}\mathrel{\hbox{\rlap
   {\lower .4em \newneq}\raise .13em \hbox{$\subset$}}}}
\def\supsetne{\setbox9 = \hbox{$\subset$}\mathrel{\hbox{\rlap
   {\lower .4em \newneq}\raise .13em \hbox{$\supset$}}}}

\def\vbar{\mathchoice{\vrule height6.3ptdepth-.5ptwidth.8pt\kern-.8pt}
   {\vrule height6.3ptdepth-.5ptwidth.8pt\kern-.8pt}
   {\vrule height4.1ptdepth-.35ptwidth.6pt\kern-.6pt}
   {\vrule height3.1ptdepth-.25ptwidth.5pt\kern-.5pt}}
\def\f@dge{\mathchoice{}{}{\mkern.5mu}{\mkern.8mu}}
\def\b@c#1#2{{\rm \mkern#2mu\vbar\mkern-#2mu#1}}
\def\b@b#1{{\rm I\mkern-3.5mu #1}}
\def\b@a#1#2{{\rm #1\mkern-#2mu\f@dge #1}}
\def\bb#1{{\count4=`#1 \advance\count4by-64 \ifcase\count4\or\b@a A{11.5}\or
   \b@b B\or\b@c C{5}\or\b@b D\or\b@b E\or\b@b F \or\b@c G{5}\or\b@b H\or
   \b@b I\or\b@c J{3}\or\b@b K\or\b@b L \or\b@b M\or\b@b N\or\b@c O{5} \or
   \b@b P\or\b@c Q{5}\or\b@b R\or\b@a S{8}\or\b@a T{10.5}\or\b@c U{5}\or
   \b@a V{12}\or\b@a W{16.5}\or\b@a X{11}\or\b@a Y{11.7}\or\b@a Z{7.5}\fi}}

\catcode`\X=11 \catcode`\@=12




\let\thischap\jobname

\def\partof#1{\csname returnthe#1part\endcsname}
\def\chapof#1{\csname returnthe#1chap\endcsname}

\def\setchapter#1,#2,#3.{%
  \expandafter\def\csname returnthe#1part\endcsname{#2}%
  \expandafter\def\csname returnthe#1chap\endcsname{#3}%
}

\setchapter 300a,A,I.
\setchapter 300b,A,II.
\setchapter 300c,A,III.
\setchapter 300d,A,IV.
\setchapter 300e,A,V.
\setchapter 300f,A,VI.
\setchapter 300g,A,VII.
\setchapter   88,B,I.
\setchapter  600,B,II.
\setchapter  705,B,III.

\def\cprefix#1{
\edef\theotherpart{\partof{#1}}\edef\theotherchap{\chapof{#1}}%
\ifx\theotherpart\thispart
   \ifx\theotherchap\thischap 
    \else 
     \theotherchap%
    \fi
   \else 
     \theotherpart.\theotherchap\fi}

\def\sectioncite[#1]#2{%
     \cprefix{#2}#1}

\edef\thispart{\partof{\thischap}}
\edef\thischap{\chapof{\thischap}}



\expandafter\ifx\csname citeadd.tex\endcsname\relax
\expandafter\gdef\csname citeadd.tex\endcsname{}
\else \message{Hey!  Apparently you were trying to
\string\input{citeadd.tex} twice.   This does not make sense.} 
\errmessage{Please edit your file (probably \jobname.tex) and remove
any duplicate ``\string\input'' lines}\endinput\fi

\sectno=-1   
\localtags
\jjtags
\NoBlackBoxes
\define\mr{\medskip\roster}
\define\sn{\smallskip\noindent}
\define\mn{\medskip\noindent}
\define\bn{\bigskip\noindent}
\define\ub{\underbar}
\define\wilog{\text{without loss of generality}}
\define\ermn{\endroster\medskip\noindent}
\define\dbca{\dsize\bigcap}
\define\dbcu{\dsize\bigcup}
\define \nl{\newline}
\magnification=\magstep 1
\documentstyle{amsppt}

{    
\catcode`@11

\ifx\alicetwothousandloaded@\relax
  \endinput\else\global\let\alicetwothousandloaded@\relax\fi

\gdef\subjclass{\let\savedef@\subjclass
 \def\subjclass##1\endsubjclass{\let\subjclass\savedef@
   \toks@{\def\usualspace{{\rm\enspace}}\eightpoint}%
   \toks@@{##1\unskip.}%
   \edef\thesubjclass@{\the\toks@
     \frills@{{\noexpand\rm2000 {\noexpand\it Mathematics Subject
       Classification}.\noexpand\enspace}}%
     \the\toks@@}}%
  \nofrillscheck\subjclass}
} 


\expandafter\ifx\csname alice2jlem.tex\endcsname\relax
  \expandafter\xdef\csname alice2jlem.tex\endcsname{\the\catcode`@}
\else \message{Hey!  Apparently you were trying to
\string\input{alice2jlem.tex}  twice.   This does not make sense.}
\errmessage{Please edit your file (probably \jobname.tex) and remove
any duplicate ``\string\input'' lines}\endinput\fi

\expandafter\ifx\csname bib4plain.tex\endcsname\relax
  \expandafter\gdef\csname bib4plain.tex\endcsname{}
\else \message{Hey!  Apparently you were trying to \string\input
  bib4plain.tex twice.   This does not make sense.}
\errmessage{Please edit your file (probably \jobname.tex) and remove
any duplicate ``\string\input'' lines}\endinput\fi

\def\renewcommand{\newcommand}	       
\edef\cite{\the\catcode`@}%
\catcode`@ = 11
\let\@oldatcatcode = \cite
\chardef\@letter = 11
\chardef\@other = 12
%
%
%
%
\def\@innerdef#1#2{\edef#1{\expandafter\noexpand\csname #2\endcsname}}%
%
%
\@innerdef\@innernewcount{newcount}%
\@innerdef\@innernewdimen{newdimen}%
\@innerdef\@innernewif{newif}%
\@innerdef\@innernewwrite{newwrite}%
%
%
%
\def\@gobble#1{}%
%
%
%
\ifx\inputlineno\@undefined
   \let\@linenumber = \empty 
\else
   \def\@linenumber{\the\inputlineno:\space}%
\fi
%
%
%
\def\@futurenonspacelet#1{\def\cs{#1}%
   \afterassignment\@stepone\let\@nexttoken=
}%
\begingroup 
\def\\{\global\let\@stoken= }%
\\ 
\endgroup
\def\@stepone{\expandafter\futurelet\cs\@steptwo}%
\def\@steptwo{\expandafter\ifx\cs\@stoken\let\@@next=\@stepthree
   \else\let\@@next=\@nexttoken\fi \@@next}%
\def\@stepthree{\afterassignment\@stepone\let\@@next= }%
%
%
%
\def\@getoptionalarg#1{%
   \let\@optionaltemp = #1%
   \let\@optionalnext = \relax
   \@futurenonspacelet\@optionalnext\@bracketcheck
}%
%
%
\def\@bracketcheck{%
   \ifx [\@optionalnext
      \expandafter\@@getoptionalarg
   \else
      \let\@optionalarg = \empty
      \expandafter\@optionaltemp
   \fi
}%
\def\@@getoptionalarg[#1]{%
   \def\@optionalarg{#1}%
   \@optionaltemp
}%
%
%
%
\def\@nnil{\@nil}%
\def\@fornoop#1\@@#2#3{}%
\def\@for#1:=#2\do#3{%
   \edef\@fortmp{#2}%
   \ifx\@fortmp\empty \else
      \expandafter\@forloop#2,\@nil,\@nil\@@#1{#3}%
   \fi
}%
\def\@forloop#1,#2,#3\@@#4#5{\def#4{#1}\ifx #4\@nnil \else
       #5\def#4{#2}\ifx #4\@nnil \else#5\@iforloop #3\@@#4{#5}\fi\fi
}%
\def\@iforloop#1,#2\@@#3#4{\def#3{#1}\ifx #3\@nnil
       \let\@nextwhile=\@fornoop \else
      #4\relax\let\@nextwhile=\@iforloop\fi\@nextwhile#2\@@#3{#4}%
}%
%
%
%
\@innernewif\if@fileexists
\def\@testfileexistence{\@getoptionalarg\@finishtestfileexistence}%
\def\@finishtestfileexistence#1{%
   \begingroup
      \def\extension{#1}%
      \immediate\openin0 =
         \ifx\@optionalarg\empty\jobname\else\@optionalarg\fi
         \ifx\extension\empty \else .#1\fi
         \space
      \ifeof 0
         \global\@fileexistsfalse
      \else
         \global\@fileexiststrue
      \fi
      \immediate\closein0
   \endgroup
}%
%
%
%
%
\def\bibliographystyle#1{%
   \@readauxfile
   \@writeaux{\string\bibstyle{#1}}%
}%
\let\bibstyle = \@gobble
%
%
\let\bblfilebasename = \jobname
\def\bibliography#1{%
   \@readauxfile
   \@writeaux{\string\bibdata{#1}}%
   \@testfileexistence[\bblfilebasename]{bbl}%
   \if@fileexists
      \nobreak
      \@readbblfile
   \fi
}%
\let\bibdata = \@gobble
%
%
\def\nocite#1{%
   \@readauxfile
   \@writeaux{\string\citation{#1}}%
}%
\@innernewif\if@notfirstcitation
%
%
\def\cite{\@getoptionalarg\@cite}%
%
%
\def\@cite#1{%
   \let\@citenotetext = \@optionalarg
   \printcitestart
   \nocite{#1}%
   \@notfirstcitationfalse
   \@for \@citation :=#1\do
   {%
      \expandafter\@onecitation\@citation\@@
   }%
   \ifx\empty\@citenotetext\else
      \printcitenote{\@citenotetext}%
   \fi
   \printcitefinish
}%
\newif\ifweareinprivate
\weareinprivatetrue
\ifx\shlhetal\undefinedcontrolseq\weareinprivatefalse\fi
\ifx\shlhetal\relax\weareinprivatefalse\fi
\def\@onecitation#1\@@{%
   \if@notfirstcitation
      \printbetweencitations
   \fi
   \expandafter \ifx \csname\@citelabel{#1}\endcsname \relax
      \if@citewarning
         \message{\@linenumber Undefined citation `#1'.}%
      \fi
     \ifweareinprivate
      \expandafter\gdef\csname\@citelabel{#1}\endcsname{%
\strut 
\vadjust{\vskip-\dp\strutbox
\vbox to 0pt{\vss\parindent0cm \leftskip=\hsize 
\advance\leftskip3mm
\advance\hsize 4cm\strut\openup-4pt 
\rightskip 0cm plus 1cm minus 0.5cm ?  #1 ?\strut}}
         {\tt
            \escapechar = -1
            \nobreak\hskip0pt\pfeilsw
            \expandafter\string\csname#1\endcsname
             \pfeilso
            \nobreak\hskip0pt
         }%
      }%
     \else  
      \expandafter\gdef\csname\@citelabel{#1}\endcsname{%
            {\tt\expandafter\string\csname#1\endcsname}
      }%
     \fi  
   \fi
   \csname\@citelabel{#1}\endcsname
   \@notfirstcitationtrue
}%
%
%
\def\@citelabel#1{b@#1}%
%
%
\def\@citedef#1#2{\expandafter\gdef\csname\@citelabel{#1}\endcsname{#2}}%
%
%
%
\def\@readbblfile{%
   \ifx\@itemnum\@undefined
      \@innernewcount\@itemnum
   \fi
   \begingroup
      \def\begin##1##2{%
         \setbox0 = \hbox{\biblabelcontents{##2}}%
         \biblabelwidth = \wd0
      }%
      \def\end##1{}
      %
      %
      \@itemnum = 0
      \def\bibitem{\@getoptionalarg\@bibitem}%
      \def\@bibitem{%
         \ifx\@optionalarg\empty
            \expandafter\@numberedbibitem
         \else
            \expandafter\@alphabibitem
         \fi
      }%
      \def\@alphabibitem##1{%
         \expandafter \xdef\csname\@citelabel{##1}\endcsname {\@optionalarg}%
         \ifx\biblabelprecontents\@undefined
            \let\biblabelprecontents = \relax
         \fi
         \ifx\biblabelpostcontents\@undefined
            \let\biblabelpostcontents = \hss
         \fi
         \@finishbibitem{##1}%
      }%
      \def\@numberedbibitem##1{%
         \advance\@itemnum by 1
         \expandafter \xdef\csname\@citelabel{##1}\endcsname{\number\@itemnum}%
         \ifx\biblabelprecontents\@undefined
            \let\biblabelprecontents = \hss
         \fi
         \ifx\biblabelpostcontents\@undefined
            \let\biblabelpostcontents = \relax
         \fi
         \@finishbibitem{##1}%
      }%
      \def\@finishbibitem##1{%
         \biblabelprint{\csname\@citelabel{##1}\endcsname}%
         \@writeaux{\string\@citedef{##1}{\csname\@citelabel{##1}\endcsname}}%
         \ignorespaces
      }%
      %
      %
      \let\em = \bblem
      \let\newblock = \bblnewblock
      \let\sc = \bblsc
      \frenchspacing
      \clubpenalty = 4000 \widowpenalty = 4000
      \tolerance = 10000 \hfuzz = .5pt
      \everypar = {\hangindent = \biblabelwidth
                      \advance\hangindent by \biblabelextraspace}%
      \bblrm
      \parskip = 1.5ex plus .5ex minus .5ex
      \biblabelextraspace = .5em
      \bblhook
      \input \bblfilebasename.bbl
   \endgroup
}%
%
%
\@innernewdimen\biblabelwidth
\@innernewdimen\biblabelextraspace
%
%
%
\def\biblabelprint#1{%
   \noindent
   \hbox to \biblabelwidth{%
      \biblabelprecontents
      \biblabelcontents{#1}%
      \biblabelpostcontents
   }%
   \kern\biblabelextraspace
}%
%
%
%
\def\biblabelcontents#1{{\bblrm [#1]}}%
%
%
\def\bblrm{\rm}%
%
%
\def\bblem{\it}%
%
%
\def\bblsc{\ifx\@scfont\@undefined
              \font\@scfont = cmcsc10
           \fi
           \@scfont
}%
%
%
\def\bblnewblock{\hskip .11em plus .33em minus .07em }%
%
%
\let\bblhook = \empty
%
%
%
\def\printcitestart{[}
\def\printcitefinish{]}
\def\printbetweencitations{, }
\def\printcitenote#1{, #1}
%
%
%
\let\citation = \@gobble
%
%
%
\@innernewcount\@numparams
%
%
\def\newcommand#1{%
   \def\@commandname{#1}%
   \@getoptionalarg\@continuenewcommand
}%
%
%
\def\@continuenewcommand{%
   \@numparams = \ifx\@optionalarg\empty 0\else\@optionalarg \fi \relax
   \@newcommand
}%
%
%
\def\@newcommand#1{%
   \def\@startdef{\expandafter\edef\@commandname}%
   \ifnum\@numparams=0
      \let\@paramdef = \empty
   \else
      \ifnum\@numparams>9
         \errmessage{\the\@numparams\space is too many parameters}%
      \else
         \ifnum\@numparams<0
            \errmessage{\the\@numparams\space is too few parameters}%
         \else
            \edef\@paramdef{%
               \ifcase\@numparams
                  \empty  No arguments.
               \or ####1%
               \or ####1####2%
               \or ####1####2####3%
               \or ####1####2####3####4%
               \or ####1####2####3####4####5%
               \or ####1####2####3####4####5####6%
               \or ####1####2####3####4####5####6####7%
               \or ####1####2####3####4####5####6####7####8%
               \or ####1####2####3####4####5####6####7####8####9%
               \fi
            }%
         \fi
      \fi
   \fi
   \expandafter\@startdef\@paramdef{#1}%
}%
%
%
%
%
\def\@readauxfile{%
   \if@auxfiledone \else 
      \global\@auxfiledonetrue
      \@testfileexistence{aux}%
      \if@fileexists
         \begingroup
            \endlinechar = -1
            \catcode`@ = 11
            \input \jobname.aux
         \endgroup
      \else
         \message{\@undefinedmessage}%
         \global\@citewarningfalse
      \fi
      \immediate\openout\@auxfile = \jobname.aux
   \fi
}%
%
%
\newif\if@auxfiledone
\ifx\noauxfile\@undefined \else \@auxfiledonetrue\fi
%
%
%
%
\@innernewwrite\@auxfile
\def\@writeaux#1{\ifx\noauxfile\@undefined \write\@auxfile{#1}\fi}%
%
%
%
\ifx\@undefinedmessage\@undefined
   \def\@undefinedmessage{No .aux file; I won't give you warnings about
                          undefined citations.}%
\fi
%
%
\@innernewif\if@citewarning
\ifx\noauxfile\@undefined \@citewarningtrue\fi
%
%
%
\catcode`@ = \@oldatcatcode

\def\pfeilso{\leavevmode
            \vrule width 1pt height9pt depth 0pt\relax
           \vrule width 1pt height8.7pt depth 0pt\relax
           \vrule width 1pt height8.3pt depth 0pt\relax
           \vrule width 1pt height8.0pt depth 0pt\relax
           \vrule width 1pt height7.7pt depth 0pt\relax
            \vrule width 1pt height7.3pt depth 0pt\relax
            \vrule width 1pt height7.0pt depth 0pt\relax
            \vrule width 1pt height6.7pt depth 0pt\relax
            \vrule width 1pt height6.3pt depth 0pt\relax
            \vrule width 1pt height6.0pt depth 0pt\relax
            \vrule width 1pt height5.7pt depth 0pt\relax
            \vrule width 1pt height5.3pt depth 0pt\relax
            \vrule width 1pt height5.0pt depth 0pt\relax
            \vrule width 1pt height4.7pt depth 0pt\relax
            \vrule width 1pt height4.3pt depth 0pt\relax
            \vrule width 1pt height4.0pt depth 0pt\relax
            \vrule width 1pt height3.7pt depth 0pt\relax
            \vrule width 1pt height3.3pt depth 0pt\relax
            \vrule width 1pt height3.0pt depth 0pt\relax
            \vrule width 1pt height2.7pt depth 0pt\relax
            \vrule width 1pt height2.3pt depth 0pt\relax
            \vrule width 1pt height2.0pt depth 0pt\relax
            \vrule width 1pt height1.7pt depth 0pt\relax
            \vrule width 1pt height1.3pt depth 0pt\relax
            \vrule width 1pt height1.0pt depth 0pt\relax
            \vrule width 1pt height0.7pt depth 0pt\relax
            \vrule width 1pt height0.3pt depth 0pt\relax}

\def\pfeilsw{ \leavevmode 
            \vrule width 1pt height0.3pt depth 0pt\relax
            \vrule width 1pt height0.7pt depth 0pt\relax
            \vrule width 1pt height1.0pt depth 0pt\relax
            \vrule width 1pt height1.3pt depth 0pt\relax
            \vrule width 1pt height1.7pt depth 0pt\relax
            \vrule width 1pt height2.0pt depth 0pt\relax
            \vrule width 1pt height2.3pt depth 0pt\relax
            \vrule width 1pt height2.7pt depth 0pt\relax
            \vrule width 1pt height3.0pt depth 0pt\relax
            \vrule width 1pt height3.3pt depth 0pt\relax
            \vrule width 1pt height3.7pt depth 0pt\relax
            \vrule width 1pt height4.0pt depth 0pt\relax
            \vrule width 1pt height4.3pt depth 0pt\relax
            \vrule width 1pt height4.7pt depth 0pt\relax
            \vrule width 1pt height5.0pt depth 0pt\relax
            \vrule width 1pt height5.3pt depth 0pt\relax
            \vrule width 1pt height5.7pt depth 0pt\relax
            \vrule width 1pt height6.0pt depth 0pt\relax
            \vrule width 1pt height6.3pt depth 0pt\relax
            \vrule width 1pt height6.7pt depth 0pt\relax
            \vrule width 1pt height7.0pt depth 0pt\relax
            \vrule width 1pt height7.3pt depth 0pt\relax
            \vrule width 1pt height7.7pt depth 0pt\relax
            \vrule width 1pt height8.0pt depth 0pt\relax
            \vrule width 1pt height8.3pt depth 0pt\relax
            \vrule width 1pt height8.7pt depth 0pt\relax
            \vrule width 1pt height9pt depth 0pt\relax
      }


\def\widestnumber#1#2{}

\def\citewarning#1{\ifx\shlhetal\relax 
    \else
    \par{#1}\par
    \fi
}

\def\rm{\fam0 \tenrm}

\def\fakesubhead#1\endsubhead{\bigskip\noindent{\bf#1}\par}



%
%
%

%

\font\textrsfs=rsfs10
\font\scriptrsfs=rsfs7
\font\scriptscriptrsfs=rsfs5

\newfam\rsfsfam
\textfont\rsfsfam=\textrsfs
\scriptfont\rsfsfam=\scriptrsfs
\scriptscriptfont\rsfsfam=\scriptscriptrsfs

\edef\oldcatcodeofat{\the\catcode`\@}
\catcode`\@11

\def\Cal@@#1{\noaccents@ \fam \rsfsfam #1}

\catcode`\@\oldcatcodeofat


\expandafter\ifx \csname margininit\endcsname \relax\else\margininit\fi

\long\def\red#1\endred{}
\long\def\green#1\endgreen{}
\long\def\blue#1\endblue{}

\def\endred{ \unmatched endred! }
\def\endgreen{ \unmatched endgreen! }
\def\endblue{ \unmatched endblue! }

\ifx\latexcolors\undefinedcs\def\latexcolors{}\fi

\def\emptycs{}
\def\evaluatelatexcolors{%
        \ifx\latexcolors\emptycs\else
        \expandafter\xxevaluate\latexcolors\xxfertig\evaluatelatexcolors\fi}
\def\xxevaluate#1,#2\xxfertig{\setupthiscolor{#1}%
        \def\latexcolors{#2}}

\font\smallfont=cmsl7
\def\rutgerscolor{\ifmmode\else\endgraf\fi\smallfont
\advance\leftskip0.5cm\relax}
\def\setupthiscolor#1{\edef\tmptmpcs{\noexpand\bgroup\noexpand\rutgerscolor
\noexpand\def\noexpand\currentcolor{#1}%
\noexpand}%
\expandafter\let\csname#1\endcsname\tmptmpcs
\def\tmptmpcs{\checkColorUnmatched{#1}\popthecolor}
\expandafter\let\csname end#1\endcsname\tmptmpcs}

\def\checkColorUnmatched#1{\def\expectcolor{#1}%
    \ifx\expectcolor\currentcolor   
    \else \edef\failhere{\noexpand\tryingToClose '\currentcolor' with end\expectcolor}\failhere\fi}

\def\currentcolor{???}

\def\popthecolor{\ifmmode\else\endgraf\fi\egroup}

\expandafter\def\csname#1\endcsname{}

\evaluatelatexcolors

 \let\outerhead\head
 \def\head{\innerhead}
 \let\innerhead\outerhead

 \let\outersubhead\subhead
 \def\subhead{\innersubhead}
 \let\innersubhead\outersubhead

 \let\outersubsubhead\subsubhead
 \def\subsubhead{\innersubsubhead}
 \let\innersubsubhead\outersubsubhead

 \def\proclaim{\innerproclaim}
 \let\innerproclaim\outerproclaim

 %
 %
 %
 %

\def\demo#1{\medskip\noindent{\it #1.\/}}
\def\enddemo{\smallskip}

\def\remark#1{\medskip\noindent{\it #1.\/}}
\def\endremark{\smallskip}

\pageheight{8.5truein}
\topmatter
\title\nofrills{THE PAIR $(\aleph_{\text{\rm n}},\aleph_0)$ MAY FAIL 
$\aleph_0$-COMPACTNESS} \endtitle
\author {Saharon Shelah \thanks {\null\newline I would like to thank 
Alice Leonhardt for the beautiful typing. \null\newline
Publication 604} \endthanks} \endauthor 

\affil{Institute of Mathematics\\
 The Hebrew University\\
 Jerusalem, Israel
 \medskip
 Rutgers University\\
 Mathematics Department\\
 New Brunswick, NJ  USA} \endaffil


\abstract  Let $P$ be a distinguished unary predicate and $\bold K =
\{M:M$ a model of cardinality $\aleph_n$ with $P^M$ of cardinality
$\aleph_0\}$.  We prove that consistently for $n=4$, for some
countable first order theory $T$ we have: $T$ has no model in $\bold
K$ whereas every finite subset of $T$ has a model in $\bold K$.
We then show how we prove it also for $n=2$, too.   \endabstract
\endtopmatter
\document  
 
\newpage

\head {Annotated Content} \endhead  \resetall 
\bn
\S0 $\quad$ Introduction
\mn
\S1 $\quad$ Relevant identities
\mr
\item "{${{}}$}"  [We deal with the 2-identities we shall use.]
\ermn
\S2 $\quad$ Definition of the forcing
\mr
\item "{${{}}$}"  [We define (historically) our forcing notion, which
depends on $\Gamma$, a set of 2-identities and on a model $M^*$ with
universe $\lambda$ and $\aleph_0$ functions. \nl
The program is to force with (the finite support product)
$\dsize \prod_n \Bbb P_{\Gamma_n}$ where the 
forcing $\Bbb P_{\Gamma_n}$ adds a colouring (= a function)
${\underset\tilde {}\to c_n}:[\lambda]^2 \rightarrow
\aleph_0$ satisfying $ID_2({\underset\tilde {}\to c_n}) \cap ID^* =
\Gamma_n$, but no $\underset\tilde {}\to c:[\lambda]^2 \rightarrow
\aleph_0$ has $ID_2(\underset\tilde {}\to c)$ too small.]
\ermn
\S3 $\quad$ Why does the forcing work
\mr
\item "{${{}}$}"  [We state the partition result in the original
universe which we shall use (in \scite{3.1}).  Then we prove that if
$\Gamma$ contains only identities which restricted to $\le m(*)$ elements are
trivial, then this holds for the colouring in any $p \in \Bbb
P_\Gamma$ (see \scite{3.1a}). \nl
We prove that $\Bbb P_\Gamma$ preserves identities from
$ID_2(\lambda,\mu)$ which are in $\Gamma$ (because we allow in the
definition of the forcing appropriate amalgamations (see
\scite{3.2}(1)).  We have weaker results for $\dsize \prod_n \Bbb
P_{\Gamma_n}$, (see \scite{3.2}(2)). \nl
On the other hand, forcing with $\Bbb P_\Gamma$ gives a colouring
showing relevant 2-identities are not in $ID_2(\lambda,\mu)$.  Lastly,
we derive the main theorem; e.g. incompactness
for $(\aleph_4,\aleph_0)$, (see (\scite{3.7}).]
\ermn
\S4 $\quad$ Improvements and Additions
\mr
\item "{${{}}$}"  [We show that we can deal with the pair
$(\aleph_2,\aleph_0)$ (see \scite{4.2} - \scite{4.6}).]  
\ermn
\S5 $\quad$ Open problems and concluding remarks
\mr
\item "{${{}}$}"  [We list some open problems, and note a property of
$ID(\aleph_n,\aleph_0)$ under the assumption MA $+ 2^{\aleph_0} >
\aleph_n$.  We note on when $k$-simple identities suffice and an
alternative proof of $(\aleph_\omega,\aleph_2) \rightarrow 
(2^{\aleph_0},\aleph_0)$.]
\endroster
\newpage

\head {\S0 Introduction} \endhead  \resetall \sectno=0
\bigskip

Interest in two cardinal models comes from the early days of 
model theory, as generalizations of the Lowenheim-Skolem theorem. 
Already Mostowski ~\cite{Mo57} considered a related
problem concerning generalized quantifiers. Let us introduce the problem.
Throughout the paper $\lambda,\mu$ and $\kappa$ stand for infinite
cardinals and $n,k$  for natural numbers.

We consider a countable language vocabulary $\tau$ with a distinguished unary 
relation symbol $P$ and models $M$ for $\tau$; i.e., $\tau$-models. \nl
\margintag{0,1}\ub{\stag{0,1} Notation}:  We let

$$
K_{(\lambda,\mu)} =: \{M:\,||M||=\lambda \and |P^M| =\mu\}.
$$
\bigskip

\definition{\stag{0.2} Definition}   1)  We say that
$K_{(\lambda,\mu)}$ is $(<\kappa)$-compact when every first 
order theory $T$ in the vocabulary $\tau$ (i.e., in the first order
logic $\Bbb L(\tau)$) with $|T| < \kappa$, satisfies:

\ub{if} every finite $t\subseteq T$ has a model in $K_{(\lambda,\mu)}$,
\ub{then} $T$ has a model in $K_{(\lambda,\mu)}$.  
\smallskip

We similarly give the meaning to $(\le\kappa)$-compactness.
We say that $(\lambda,\mu)$ is $(<\kappa)$-compact if
$K_{(\lambda,\mu)}$ is. \nl
2)  We say that

$$
(\lambda,\mu) \rightarrow'_\kappa (\lambda',\mu')
$$
\mn
when for every first order theory $T$ in $\Bbb L(\tau)$ with $|T| < \kappa$,
\ub{if} every finite $t\subseteq T$ has a model in $K_{(\lambda,\mu)}$,
\ub{then} $T$ has a model in $K_{(\lambda',\mu')}$.  Instead
``$\kappa^+$" we may write ``$\le \kappa$".
\nl
3)  We say that

$$
(\lambda,\mu) \rightarrow_\kappa (\lambda',\mu')
$$
\mn
when for every first order theory $T$ of $L$ with $|T| < \kappa$,
\ub{if} $ T$ has a model in $K_{(\lambda,\mu)}$,
\ub{then} $T$ has a model in $K_{(\lambda',\mu')}$. \nl
4) In both $\rightarrow'_\kappa$ and $\rightarrow_\kappa$ we omit 
$\kappa$ if $\kappa=\aleph_0$.
\enddefinition
\bn
\ub{Note}: Note that $\rightarrow_\kappa$ is transitive and 
$\rightarrow'_\kappa$ is as well. 
Also note that
$\rightarrow_{\aleph_0}$ and $\rightarrow'_{\aleph_0}$ are equivalent.

We consider the problem of $K_{(\lambda,\mu)}$ being compact. 
Before we start, we review the history of the problem. Note that a related problem is the 
one of completeness, i.e. if

$$
\{\psi:\,\psi \text{ has a model in } K_{(\lambda,\mu)}\}
$$
\mn
is recursively enumerable and other related problems, see in the end. 
We do not concentrate on those problems here.

We review some of the history of the problem, in an order which
is not necessarily chronological.

Some early results on the compactness are due to Furkhen
\cite{Fu65}. He showed that
\mr
\item "{$(A)$}"   if $\mu^\kappa=\mu$, \ub{then} 
$K_{(\lambda,\mu)}$ is $(\le\kappa)$-compact.
\ermn
The proof is by using ultraproducts over regular ultrafilters on
$\kappa$, generalizing the well known proof of compactness by
ultrapowers.  A related result is
\mr
\item "{$(B)$}"   (\cite{Mo68})  If 
$\mu^{\aleph_0}\le\mu'\le\lambda'\le\lambda$, \ub{then}
$(\lambda,\mu) \rightarrow_{\le \lambda}(\lambda',\mu')$.
\ermn
Next result we mention is one of Silver \cite{Si} concerning Kurepa trees, 
\mr
\item "{$(C)$}"   (Silver \cite{Si71}) 
\relax From the existence of a strongly inaccessible cardinal,
it follows that the following is consistent  with $ZFC$:
$$
GCH + (\aleph_3,\aleph_1) \nrightarrow_{\aleph_0}(\aleph_2,\aleph_0).
$$
\ermn
Using special Aronszajn trees Mitchell showed 
\mr
\item "{$(D)$}"   (Mitchell \cite{Mi}) From the 
existence of a Mahlo cardinal, it follows
that it is consistent  with $ZFC$ to have
$$
(\aleph_1,\aleph_0) \nrightarrow_{\aleph_2}(\aleph_2,\aleph_1).
$$
\ermn
A later negative consistency result is the one of Schmerl in
\cite{Sc1}
\mr
\item "{$(E)$}"   (Schmerl \cite{Sc74}) Con(if $n<m$ 
then $(\aleph_n,\aleph_{n+1}) \nrightarrow (\aleph_m,\aleph_{m+1}))$.
\ermn
Earlier, Vaught proved two positive results 
\mr
\item "{$(F)$}"   (Vaught \cite{MV62}) 
$(\lambda^+,\lambda) \rightarrow'_{\aleph_1} (\aleph_1,\aleph_0)$.
\ermn
Keisler \cite{Ke66} and \cite{Ke66a} has obtained more results 
in this direction.
\mr
\item "{$(G)$}"   (Vaught \cite{Va65}) If $\lambda \ge \beth_\omega(\mu)$
and $\lambda'>\mu'$, \ub{then} $(\lambda,\mu) \rightarrow'_{\le \mu'}(\lambda',\mu')$.
\ermn
In \cite{Mo68} Morley gives another proof of this result, using
Erd\"os-Rado Theorem and indiscernibles. \nl
Another early positive result is the one of Chang:
\mr
\item "{$(H)$}"   (Chang \cite{Ch65})  If $\mu=\mu^{<\mu}$ {\ub{then}
$(\lambda^+,\lambda) \rightarrow'_{\le \mu} (\mu^+,\mu)$.
\ermn
Jensen in \cite{Jn} uses $\square_\mu$ to show
\mr
\item "{$(I)$}"   (Jensen \cite{Jn})  If $\bold V = \bold L$, \ub{then}
$(\lambda^+,\lambda) \rightarrow'_{\le \mu} (\mu^+,\mu).$ (The fact that
$0^\#$ does not exist suffices.)
\ermn
Hence, Jensen's result deals with the case of
$\mu$ is singular, which was left open after the result of Chang. 
For other early consistency results
concerning gap-1 two cardinal theorems, including consistency,
see \cite{Sh:269}, Cummings, Foreman and Magidor \cite{CFM0x}].

In \cite{Jn} there is actually a simplified proof
of (I) due to Silver.  A further result of Jensen, using morasses, is:
\mr
\item "{$(J)$}"   (Jensen, see \cite{De73}) (for $n=2$)) If $\bold V =
\bold L$, \ub{then} 
$(\lambda^{+n},\lambda) \rightarrow'_{\le \mu} (\mu^{+n},\mu)$ for all $n
< \omega$.
\ermn
Note that by Vaught's result \cite{MV62} stated in (F) we have:
the statement in (I), in
the result of Chang etc., $(\lambda^+,\lambda)$ can be
without loss of generality replaced by $(\aleph_1,\aleph_0)$.

Finally, there are many more related results, for example the ones
concerning Chang's conjecture. A survey article on the topic
was written by Schmerl in \cite{Sc74}. \nl
Note that many of the positive rsults above (F)-(J), their proof also
gives compactness of the pair, e.g., $(\aleph_0,\aleph_1)$ by \cite{MV62}.

We now mention
some results of the author which will have a bearing to the present
paper.
\mr
\item "{$(\alpha)$}"  (Shelah \cite{Sh:8} and the abstract
\cite{Sh:E17}).   If $K_{(\lambda,\mu)}$ is
$(\le\aleph_0)$-compact, \ub{then} $K_{(\lambda,\mu)}$ is $(\le\mu)$-compact and
$(\lambda,\mu) \rightarrow_{\le \mu'} (\lambda',\mu')$ when $\lambda \le \lambda'
\le \mu' \le \mu$.
\ermn
More than $(\le\mu)$-compactness cannot hold for trivial reasons.
In the same work we have the analogous result on $\rightarrow'$ and:
\mr
\item "{$(\beta)$}"   (Shelah \cite{Sh:8} and the abstract 
\cite{Sh:E17}) $(\lambda,\mu)
\rightarrow'_{\aleph_1} (\lambda',\mu')$ is actually a problem on
partition relations, (see below), also it implies $(\lambda,\mu)
\rightarrow'_{\mu'} (\lambda',\mu')$ see \scite{0.4a}(1) below.
\ermn
We state a definition from \cite{Sh:8} that will be used here too.
We do not consider the full generality of \cite{Sh:8}, there problems
like considering $K$ with several $\lambda_\ell$-like 
$(P^2_\ell,<_\ell)$ and $|P^1_\ell| = \mu_\ell$ were addressed. \nl
(We can use below only ordered $a$ and increase $h$, it does not
matter much.)
\bigskip

\definition{\stag{0.4} Definition}  1) An \ub{identity}
\footnote{identification in the terminology of \cite{Sh:8}} is a pair $(a,e)$ where
$a$ is a finite set and $e$ is an equivalence relation on the finite subsets
of $a$, having the property

$$
b\, e\, c \Rightarrow |b| = |c|.
$$
\mn
The equivalence class of $b$ with respect to $e$ will be denoted
$b/e$. \nl
2)  We say that $\lambda \rightarrow (a,e)_\mu$, if for every
$f:\,[\lambda]^{<\aleph_0}\to\mu$, there is $h:\,a \overset {1-1}
\to\longrightarrow \lambda$ such that

$$
b\, e\, c \Rightarrow f(h"(b)) = f(h"(c)).
$$
\mn
3)  We define

$$
ID(\lambda,\mu) =: \{(n,e):\,n < \omega \and (n,e)
\text{ is an identity and } \lambda \rightarrow (n,e)_\mu\}
$$
\mn
and for $f:\,[\lambda]^{<\aleph_0} \rightarrow X$ we let

$$
\align
ID(f) =: \{(n,e):&(n,e) \text{ is an identity such that for some
one-to-one function} \\
  &h \text{ from } n = \{0,\dotsc,n-1\} \text{ to } \lambda \text{ we
have} \\
  &(\forall b,c \subseteq n)(b\, e\, c \Rightarrow f(h''(b)) = f(h''(c)))\}
\endalign
$$
\enddefinition
\bigskip

\proclaim{\stag{0.4a} Claim}  (Shelah \cite{Sh:8} 
and the abstract \cite{Sh:E17})
$(\lambda,\mu) \rightarrow'_{\aleph_1}(\lambda',\mu')$
is equivalent to the existence of a function
$f:\,[\lambda']^{<\aleph_0} \rightarrow \mu'$ such that 

$$
ID(f) \subseteq ID(\lambda,\mu)
$$
\mn
(more on this see \cite[Th.3]{Sh:74} statement there on
$\rightarrow'_{\aleph_1}$, see details in \cite{Sh:E28}). 
\endproclaim
\bigskip

\remark{\stag{0.4A} Remark}  The identities of 
$(\beth_\omega,\aleph_0)$ are clearly characterized by 
Morley's proof of Vaught's theorem (see \cite{Mo68}). 
The identities of $(\aleph_\omega,\aleph_0)$ are stated explicitly in 
\cite{Sh:37} and
\cite{Sh:49}, when $\aleph_\omega \le 2^{\aleph_0}$ where it 
is also shown that $(\aleph_\omega,\aleph_0) \rightarrow'
(2^{\aleph_0},\aleph_0)$. For $(\aleph_1,\aleph_0)$, the identities
are characterized in \cite{Sh:74} (for some details see \cite{Sh:E28}).
The identities for $\lambda$-like models, $\lambda$ $\omega$-strongly
$n$-Mahlo are clear, see Schmerl and Shelah \cite{ScSh:20} (for
strongly $n$-Mahlo this gives positive results, subsequently sharpened
(replacing $n+2$ by $n$) and the negative results proved by Schmerl,
see \cite{Sch85}).
\endremark
\bigskip

\demo{Proof}  By the referee request we indicate the proof.
\enddemo
\bigskip

\remark{Remark}  We generally neglect here three cardinal theorems and
$\lambda$-like model (and combinations, see \cite{Sh:8}, \cite{Sh:18},
the positive results (like \scite{0.4a} are similar.  See recently
\cite{ShVa:790}.
\endremark
\bn
In Gilschrist, Shelah \cite{GcSh:491} and
\cite{GcSh:583}, we dealt with 2-identities.
\definition{\stag{0.4b} Definition}  1) A two-identity or $2$-identity
\footnote{it is not an identity as $e$ is an equivalence relation on
too small set} is a pair $(a,e)$ where $a$ is a finite set and $e$ 
is an equivalence relation on $[a]^2$. 
Let $\lambda \rightarrow (a,e)_\mu$ mean $\lambda \rightarrow
(a,e^+)_\mu$ where $be^+ c \leftrightarrow (bec) \vee (b = c \subseteq
a)$ for any $b,c \subseteq a$.
\nl
2) We defined

$$
ID_2(\lambda,\mu) =: \{(n,e):\,(n,e) \text{ is a } 2\text{-identity and }
\lambda \rightarrow(n,e)_\mu\}
$$
\mn
we define $ID_2(f)$ when $f:[\lambda]^2 \rightarrow X$ as

$$
\align
\biggl\{(n,e):&(n,e) \text{ is a two-identity such that for some } h,
\\
  & \text{ a one-to-one function from } 
\{0,\dotsc,n-1\} \text{ into } \lambda \\
  &\text{ we have } \{\ell_1,\ell_2\} e \{k_1,k_2\} \text{ implies that }
\ell_1 \ne \ell_2 \in \{0,\dotsc,n-1\}, \\
  &k_1 \ne k_2 \in \{0,\dotsc,n-1\} \text{ and }
f(\{h(\ell_1),h(\ell_2)\}) = f(\{h(k_1),h(k_2)\}) \biggr\}.
\endalign
$$
\mn
3) Let us define

$$
\align
ID^\circledast_2 =: 
\{&({}^n 2,e):\,({}^n 2,e) \text{ is a two-identity and if} \\
  &\{\eta_1,\eta_2\} \ne \{\nu_1,\nu_2\} \text{ are } \subseteq {}^n
2, \text{ then} \\
  &\{\eta_1,\eta_2\}e\{\nu_1,\nu_2\} \Rightarrow \eta_1 \cap \eta_2
= \nu_1 \cap \nu_2\}.
\endalign
$$
\mn
By \cite{Sh:49}, under the assumption $\aleph_\omega < 2^{\aleph_0}$,
the families $ID_2(\aleph_\omega,\aleph_0)$ and $ID^\circledast_2$ 
coincide (up
to an isomorphism of identities). 
In Gilchrist and Shelah \cite{GcSh:491} and \cite{GcSh:583}
we considered the question of the equality between these
$ID_2(2^{\aleph_0},\aleph_0)$ and $ID^\circledast_2$
under the assumption $2^{\aleph_0} = \aleph_2$. We showed that
consistently the answer may be ``yes" and may be ``no".

Note that $(\aleph_n,\aleph_0) \nrightarrow (\aleph_\omega,\aleph_0)$ 
so $ID(\aleph_2,\aleph_0) \ne ID(\aleph_\omega,\aleph_0)$, but for
identities for pairs (i.e. $ID_2$) the question is meaningful.
\sn
The history of the problem suggested to me that there should be a model where
$K_{(\lambda,\mu)}$ is not $\aleph_0$-compact for some
$\lambda,\mu$; I do not know about the opinion of others and it was
not easy for me as I thought a priori.
As mathematicians do not feel that
a strong expectation makes a proof, I was quite happy to be able to
prove the existence of such a model.
This was part of my lectures in a 1995 seminar in Jerusalem and notes
of the lecture were taken by Mirna Dzamonja and I thank her for this, but
because the proof was not complete, it was delayed.
\enddefinition
\bn
The following is the main result of this paper (proved in \scite{3.7}): 
\proclaim{\stag{0.5} Main Theorem}  {\rm Con}(the pair 
$(\aleph_n,\aleph_0)$ is not $\aleph_0$-compact 
$+ 2^{\aleph_0} \ge \aleph_n$) for $n \ge 4$.
\endproclaim
\bigskip

Later in the paper we deal with the case $n = 2$ which is somewhat more
involved.  This is the simplest case by a reasonable measure: if you
do not like to use large cardinals then assuming that there is no
inaccessible in $\bold L$, all pairs $(\mu^+,\mu)$ are
known to be $\aleph_0$-compact and if $\bold V = \bold L$
also all logic $L(\exists^{\ge
\lambda}),\lambda > \aleph_0$ are (by putting together already known
results;  $\bold V = \bold L$ is used just to imply that there is no
limit, uncountable not strong limit cardinal).

How much this consistency result will mean to a model theorist,
let us not elaborate, but instead say an anecdote about Jensen. 
He is reputed to have said: `` When I started working on the two-cardinal
problem, I was told it was the heart of model theory. Once I succeeded
to prove something, they told me what I did was pure set theory, and
were not very interested"; also, mathematics is not immune to fashion changes.

My feeling is that 
there are probably more positive theorems in this 
subject waiting to be discovered. Anyway, let us state the following
\nl
\ub{Thesis}  Independence results help us 
clear away the waste, so the possible treasures can stand out.

Of course, I have to admit that, having spent 
quite some time on the independence results, I sometimes look for 
the negative of the picture given by this thesis.

The strategy of our proof is as follows.
It seems natural to consider the simplest case, i.e., that of 
two-place functions, and try to get
the incompactness by constructing a sequence 
$\langle f_k:k < \omega\rangle$ of functions from $[\aleph_n]^2$ into $\aleph_0$ such that
for all $n$ we have $ID_2(f_k) \supseteq ID_2(f_{k+1})$, yet for no
$f:\,[\aleph_n]^2 \rightarrow \aleph_0$ do we have
$ID_2(f) \subseteq \dbcu_{k<\omega} ID_2(f_k)$.
This suffices.
Related proofs to our main results were \cite{Sh:522}. \nl
Note that another interpretation of \scite{0.5} is that if we add to
first ordre logic the cardinality quantifiers $(\exists^{\ge \lambda}
x)$ for $\lambda = \aleph_1,\aleph_2,\aleph_3,\aleph_4$ we 
get a noncompact logic.

We thank the referee for many helpful comments and the reader should
thank him also for urging the inclusion of several proofs.

This work is continued in \cite{ShVa:790} and \cite{Sh:F556}.
\newpage

\head {\S1 Relevant Identities} \endhead  \resetall \sectno=1
\bigskip

We commence by several definitions.
For simplicity, for us all identities, colorings etc. will be 2-place.
\definition{\stag{1.0} Definition}  1) For $m,\ell < \omega$ let

$$
\text{dom}_{\ell,m} = \{\eta \in {}^{\ell +1} \omega:\eta \restriction
\ell \in {}^\ell 2 \text{ and } \eta(\ell) < m\}
$$

$$
\align
ID^1_{\ell,m} = \{(\text{dom}_{\ell,m},e):\,&e \text{ is an
equivalence relation on [dom}_{\ell,n}]^2 \\
  &\,\text{such that } \{\eta_1,\eta_2\} e \{\nu_1,\nu_2\} \and
\{\eta_1,\eta_2\} \ne \{\nu_1,\nu_2\} \\
  &\, \Rightarrow \eta_1 \cap \eta_2 = \nu_1 \cap \nu_2\}.
\endalign
$$
\mn
2) Let

$$
ID^1_\ell = \cup\{\text{ID}^1_{\ell,m}:m < \omega\}
$$

$$
ID^1 = \cup\{ID^1_\ell:\ell < \omega\}.
$$
\mn
3) For $\bold s = (\text{dom}_{\ell,m},e) \in ID^1_{\ell,m}$ 
and $\nu \in {}^{\ell \ge} 2$ let

$$
\text{dom}^{[\nu]}_{\ell,m} = \{\rho \in \text{ dom}_{\ell,m}:\nu
\trianglelefteq \rho\}
$$
\mn
and if $\nu \in {}^{\ell >}2$ we let

$$
e_{<\nu>}(\bold s) = e \restriction \{\{\eta_0,\eta_1\}:\nu \char 94 <i>
\triangleleft \eta_i \text{ for } i=0,1\}.
$$
\mn
We use $\bold s$ to denote identities so $\bold s = (\text{dom}_{\bold
s},e(\bold s))$; and if $\bold s \in ID^1$ then let 
$\bold s = (\text{dom}_{\ell(\bold s),m(\bold s)},e(\bold s))$.
\nl
4)  An equivalence class is nontrivial \ub{if} it is not a singleton.
\enddefinition
\bn
Note that it follows that every $e$-equivalence class is an $e_{<\nu>}$
-equivalence class for some $\nu$.
We restrict ourselves to 
\definition{\stag{1.2} Definition}  1) Let $ID^2_{\ell,m}$ be
the set of $\bold s \in ID^1_{\ell,m}$ 
such that for every $\nu \in {}^{\ell \ge} 2$ the
equivalence relation $e_{<\nu>}(\bold s)$ has at most
one non-singleton equivalence class, which we call $e_{[\nu]} =
e_{[\nu]}(\bold s)$. 

So we 
also allow $e_{<\nu>} =$ equality, in which case we choose a representative
equivalence class $e_{[\nu]}$ as the first one under, say,
lexicographical ordering. \nl
2) $ID^2_\ell = \cup\{ID^2_{\ell,m}:m < \omega\}$.
\enddefinition
\bigskip

\definition{\stag{1.3} Definition}  1) We define for $k < \omega$ 
when $\bold s = (\text{dom}_{\ell,m},e)$ is $k$-nice: the demands are
\mr
\item "{$(a)$}"  $\bold s \in ID^1_{\ell,m}$
\sn
\item "{$(b)$}"   if $\nu \in {}^\ell 2$ and $(\nu \restriction i)
\char 94 \langle 1 - \nu(i) \rangle \triangleleft \rho_i \in
\text{ dom}_{\ell,m}$ for each $i < \ell$ then 
$\{\eta:\nu \triangleleft \eta \in \text{
dom}_{\ell,m}$ and for each $i < \ell$ the set 
$\{\rho_i,\eta\}/e$ is not a singleton$\}$ has
at least two members
\sn
\item "{$(c)$}"  the graph $H[e]$, see below, has no cycle $\le k$
(for $k \le 2$ this holds trivially)
\sn
\item "{$(d)$}"  the graph $H[e]$ has a cycle.
\ermn
2) We can interpret $\bold s = (\text{dom}_{\ell,m},e)$ as the graph
$H[\bold s]$ with set of nodes dom$_{\ell,m}$ and set of edges
$\{\{\eta,\nu\}:\{\eta,\nu\}/e$ not a singleton (and of course $\eta
\ne \nu$ are from dom$_{\ell,m})\}$. \nl
3) We may write $e(\bold s)$ instead of $\bold s$ if 
dom$_{\ell,m}$ can be reconstructed from $e$ (e.g. if the graph
has no isolated point (e.g. if it is $0$-nice, see clause (b) of part
(1)).  
Saying nice we mean $[\text{log}_2(m)]$-nice. 
\enddefinition
\bigskip

\proclaim{\stag{1.4} Claim}  1) If $(\lambda,\mu)$ is
$\aleph_0$-compact and $c_n:[\lambda]^{< \aleph_0} \rightarrow \mu$
and $\Gamma_n = \text{ ID}(c_n)$ for $n < \omega$, 
\ub{then} for some $c:[\lambda]^{<
\aleph_0} \rightarrow \mu$ we have $ID(c) \subseteq \dbca_{n < \omega}
\Gamma_n$ (in fact equality holds). \nl
2) Similarly using $ID_2$.
\endproclaim
\bigskip

\remark{Remark}  By the same proof, if we just assume 
$(\lambda_1,\mu_1) \rightarrow'_{\aleph_1}
(\lambda_2,\mu_2)$ and $c_n:[\lambda_1]^{< \aleph_0} \rightarrow
\mu_1$, \ub{then} we can deduce that there is 
$c:[\lambda_2]^{< \aleph_0} \rightarrow \mu_2$ satisfying $ID(c)
\subseteq \dbca_{n < \omega} ID(c_n)$.
\endremark
\bigskip

\demo{Proof}  Straightforward. \nl
1) In details, let $F_m$ be an $m$-place
function symbol and $P$ the distinguished unary predicate and let $T =
\{\psi_n:n < \omega\} \cup \{\neg \psi_{\bold s}:c$ is an identity of
the form $(n,e)$ not from $\dbca_{n < \omega} ID(c_n)\}$ where
\mr
\item "{$(a)$}"  $\psi_n = (\forall x_0)(\forall x_1)\ldots(\forall
x_{n-1})(P(F(x_0,\dotsc,x_{n-1})) \and \wedge \{\forall
x_0)\ldots(\forall x_n)F(x_0,\dotsc,x_{n-1}) =
F(x_{\pi(0)},\dotsc,x_{\pi(n-1)}):\pi$ a permutation of
$\{0,\dotsc,n-1\}\}$
\sn
\item "{$(b)$}"  if $\bold s = (n,e)$ is an identity 
then $\psi_{\bold s} = (\exists
x_0)\ldots(\exists x_{n-1})[\dsize \bigwedge_{\ell < m < n} x_\ell \ne
x_m \and \dsize \bigwedge_{b,c \subseteq n,b e c}
F_{|b|}(\dotsc,x_\ell,\ldots)_{\ell \in b} = F_{|b|}(\dotsc,x_\ell,
\ldots,)_{\ell \in c}]$.
\ermn
Clearly $T$ is a (first order) countable theory so it suffices to prove the
following two statements $\boxtimes_1,\boxtimes_2$.
\mr
\item "{$\boxtimes_1$}"  if $M \in K_{(\lambda,\mu)}$ is a model of
$T$, \ub{then} there is $c:[\lambda]^{< \aleph_0} \rightarrow \mu$ such
that $ID(c) \subseteq \dbca_{n < \omega} \Gamma_n$. \nl
[Why does $\boxtimes_1$ hold?  There is $N \cong M$ such that $N$ has
universe $|N| = \lambda$ and $P^N = \mu$.  Now we define $c$: if $u
\in [\lambda]^{< \aleph_0}$, let $\{\alpha^u_\ell:\ell < |u|\}$
enumerate $u$ in increasing order and let $c(u) =
F^N_{|u|}(\alpha^u_0,\alpha^u_1,\dotsc,\alpha^u_{|u|-1})$.  Note that because
$N \models \psi_n$ for $n < \omega$ clearly $c$ is a function from
$[\lambda]^{< \aleph_0}$ \ub{into} $\mu$.  Also because $N \models
\psi_n$,  if $n < \omega$ and
$\alpha_0,\dotsc,\alpha_{n-1} < \lambda$ are with no repetitions then
$F^N_n(\alpha_0,\dotsc,\alpha_{n-1}) =
c\{\alpha_0,\dotsc,\alpha_{n-1}\}$.  Now if $\bold s \in ID(c)$ let
$\bold s = (n,e)$ and let $u = \{\alpha_0,\dotsc,\alpha_{n-1}\} \in
[\lambda]^n \subseteq [\lambda]^{< \aleph_0}$ exemplify that $\bold s
\in ID(c)$, hence easily $N \models \psi_{\bold s}$ so necessarily $\neg
\psi_{\bold s} \notin T$ hence $\bold s \in \dbca_{n < \omega}
\Gamma_n$.  This implies that $ID(c) \subseteq \dbca_{n < \omega}
\Gamma_n$ is as required.]
\sn
\item "{$\boxtimes_2$}"  if $T' \subseteq T$ is finite then $T'$ has a
model in $K_{(\lambda,\mu)}$.
\ermn
[Why?  So $T'$ is included in $\{\psi_m:m < m^*\} \cup \{\neg \psi_{\bold
s_k}:k < k^*\}$ for some $m^* < \omega,\bold s_k = (n_k,e_k)$ an
identity not from $\dbca_{\ell < \omega} ID(c_\ell)$, so we can find
$\ell(k) < \omega$ such that $\bold s_k \notin ID(c_{\ell(k)})$.  Let
$H$ be a one-to-one function from ${}^{k^*}\mu$ into $\mu$.  We define
a model $M$:
\mr
\item "{$(a)$}"  its universe $|M|$ is $\lambda$
\sn
\item "{$(b)$}"  $P^M = \mu$
\sn
\item "{$(c)$}"  if $n < \omega,\{\alpha_0,\dotsc,\alpha_{n-1}\} \in
[\lambda]^n$ then \nl
$F^M_n(\alpha_0,\dotsc,\alpha_{n-1}) =
H(c_{\ell(0)}\{\alpha_0,\dotsc,\alpha_{n-1}\}$, \nl
$c_{\ell(1)}\{\alpha_0,\dotsc,\alpha_{n-1}\},\dotsc,
c_{\ell(k^*-1)}\{\alpha_0,\dotsc,\alpha_{n-1}\})$.
\ermn
If $n < \omega$ and $\alpha_0,\dotsc,\alpha_{n-1} < \lambda$ is with
repetitions we let $F^M_n(\alpha_0,\dotsc,\alpha_{n-1}) = 0$. Clearly
$M$ is a model from $K_{(\lambda,\mu)}$ of the vocabulary of $T$.
Also $M$ satisfies each sentence $\psi_m$ by the way we have defined
$F^M_m$.  Lastly, for $k<k^*,M \models \psi_{\bold s_k}$ because
$(n_k,e_k) \notin ID(c_{\ell(k)})$ by the choice of the $F_n$'s as $H$
is a one-to-one function.]  \hfill$\square_{\scite{1.4}}$\margincite{1.4}
\enddemo
\bn
Of course
\demo{\stag{1.5} Observation}  1) For every $\ell < \omega,k < \omega$
for some $m$ there is a $k$-nice $\bold s =
(\text{dom}_{\ell,m},e)$. \nl
2) If $\bold s$ is $k$-nie and $m \le k$, \ub{then} $\bold s$ is $m$-nice.
\enddemo
\newpage

\head {\S2 Definition of the Forcing} \endhead  \resetall \sectno=2
\bigskip

We have outlined the intended end
of the proof at the end of the introductory section.
It is to construct a sequence of functions $\langle f_n:\,n < \omega
\rangle$ with certain properties. 
As we have adopted the decision of dealing 
only with 2-identities from $ID_\ell$, 
all our functions will be colorings of pairs, and we shall
generally use the letter $c$ for them.

Our present theorem \scite{0.5} deals with $\aleph_4$, but we
may as well be talking about some $\aleph_{n(\ast)}$ for a
fixed natural number $n(\ast)\ge 2$.
Of course, the set of identities will depend on $n(\ast)$.
We shall henceforth work with $n(\ast)$, keeping in mind that
the relevant case for Theorem \scite{0.5} is $n(\ast)=4$.
Also we fix $\ell(*) = n(*)+1$ on which the identities depend (but
vary $m$).
Another observation about the proof is that we can replace
$\aleph_0$ with an uncountable cardinal $\kappa$ such that
$\kappa=\kappa^{<\kappa}$
replacing $\aleph_n$ by $\kappa^{+n}$.  Of course, the pair
$(\kappa^{+n},\kappa)$ is compact because $[\kappa = \kappa^{\aleph_0} <
\lambda \Rightarrow (\kappa,\lambda)$ is $\le \kappa$-compact],
however, much of the analysis holds.

We may replace $(\aleph_n,\aleph_0)$ by $(\kappa^{+n(*)},\kappa)$ if
$\kappa^{+n(*)} \le 2^{\aleph_0}$; we hope to return to this elsewhere. \nl
To consider $(\kappa^+,\kappa)$ we need large cardinals; even more so
for considering $(\mu^+,\mu),\mu$ strong limit singular of
cofinality $\aleph_0$, and even $(\kappa^{+n},\kappa),\mu \le \kappa <
\kappa^{+n} \le \mu^{\aleph_0}$.

We now describe the idea behind the definition
of the forcing notion we shall be concerned with.  Each ``component"
of the forcing notion is supposed to add a coloring

$$
c:\,[\lambda]^2 \rightarrow \mu
$$
\mn
preserving some of the possible 2-identities, while ``killing" all those
which were not preserved, in other words it is concerned with adding
$f_n$;  specifically we concentrate on the case
$\lambda = \aleph_{n(*)},\mu = \aleph_0$.  
Hence, at first glance the forcing will be defined
so that to preserve an identity we have to work hard proving some kind
of amalgamation for the forcing notion, while killing
an identity is a consequence of adding a colouring exemplifying
it.
By preserving a set $\Gamma$ of identities, we mean that
$\Gamma \subseteq ID(c)$, and more seriously
$\Gamma \subseteq ID_2(\lambda,\mu)$; we restrict
ourselves to some $ID^*$, an infinite set of 2-identity. \nl
We shall choose $ID^* \subseteq ID^\circledast_2$ below small enough such we
can handle the identities in it.

We define the forcing by putting  in its definition, for each
identity that we want to preserve, a clause specifically assuring this.
Naturally this implies that not only the desired identities are preserved,
but also some others so making an identity be not in
$ID(\lambda,\mu)$ becomes now the hard part. 
So, we lower our sights and simply hope that, if
$\Gamma\subseteq ID^\ast$ is the set of identities that we want to
preserve, than no identity $(a,e)\in ID^\ast\setminus\Gamma$ is
preserved; this may depend on $\Gamma$. 

How does this control over the set of identities help to obtain the
non-compactness?  We shall choose sets $\Gamma_n\subseteq ID^\ast$ of possible
identities for $n<\omega$. The forcing we referred to above, let us
call it $\Bbb P^{\Gamma_n}$, add a colouring ${\underset\tilde {}\to
c_n}:[\lambda]^2 \rightarrow \omega$ such that $ID_2({\underset\tilde {}\to
c_n})$ includes $\Gamma_n$ and is disjoint to $ID^* \backslash
\Gamma_n$; also it 
will turn out to have a strong form of the $ccc$. We shall force
with ${\Bbb P} =: \dsize \prod_{n\in\omega} \Bbb P^{\Gamma_n}$, where the product
is taken with finite support.  Because of the
strong version of $ccc$ possessed by each $\Bbb P^{\Gamma_n}$, 
also ${\Bbb P}$ will have $ccc$.
Now, in $\bold V^{{\Bbb P}}$ we have for every $n$ a 
colouring $c_n:\,[\lambda]^2 \rightarrow \omega$ which preserves the
identities in $\Gamma_n$, moreover $\bold V^{\Bbb P} \models \Gamma_n \subseteq
ID(c_n) \cap ID^*$.

We shall in fact obtain that

$$
ID^* = \Gamma_0 \supseteq \Gamma_1 \and \Gamma_1 \supseteq \Gamma_2
\and \ldots \and \dbca_{n < \omega}
\Gamma_n = \emptyset \and ID(c_n) \cap \Gamma_0=\Gamma_n.
$$
\mn
If we have $\aleph_0$-compactness for $(\lambda,\aleph_0)$, then by
\scite{1.4}(2)  there must be a colouring 
$c:\,[\lambda]^2 \rightarrow \omega$ in $\bold V^{\Bbb P}$ such that

$$
ID_2(c) \cap \Gamma_0 \subseteq \dbca_{n < \omega} \Gamma_n =
\emptyset.
$$
\mn
We can find a name $\underset\tilde {}\to c$ in $\bold V$ 
for such $c$, so by $ccc$, for every
$\{\alpha,\beta\} \in [\lambda]^2$, 
the name $\underset\tilde {}\to c
(\{\alpha,\beta\})$ depends only on $\aleph_0$ ``coordinates".
At this point a first approximation to what we do is to apply a
relative of Erd\"os-Rado theorem to prove that there are an $n$, a
large enough $w \subseteq \lambda$ and for every $\{\alpha,\beta\} \in [w]^2$ a condition
$p_{\{\alpha,\beta\}}\in \dsize \prod_{\ell < n} \Bbb P^{\Gamma_l}$,
such that $p_{\{\alpha,\beta\}}$ forces a value to 
$\underset\tilde {}\to c(\{\alpha,\beta\})$ in a ``uniform" enough way.
We shall be able to extend enough of the conditions $p_{\{\alpha,\beta\}}$
by a single condition $p^\ast$ in $\dsize \prod_{\ell < n} \Bbb P^{\Gamma_l}$, which gives
an identity in $ID_2(\underset\tilde {}\to c)$ which belongs to
$\dbca_{\ell < n} \Gamma_\ell \setminus \Gamma_n$, contradiction.

Before we give the definition of the forcing, we need to introduce a 
notion of closure. 
The properties of the closure operation are the ones possible
to obtain for $(\lambda,\aleph_0)$, 
but not for $(\aleph_\omega,\aleph_0)$. We of
course need to use somewhere such a property, as we know in $ZFC$ that
$(\aleph_\omega,\aleph_0)$ has all those identities, i.e. $ID^\circledast_2 =
ID_2(\lambda,\aleph_0)$.
\nl
On a similar proof see \cite{Sh:424} (for $\omega$-place functions) and 
also (2-place functions), \cite{Sh:522}. The definition of the closure 
in \cite{GcSh:491} is close to ours,
but note that the hard clause from \cite{GcSh:491} is not needed here. 
\bigskip

\definition{\stag{2.0} Definition}  Let $ID^*_{\ell(*)} =: 
\{\bold s \in ID^2_{\ell(*)}:\bold s$ is $0$-nice$\}$.
\enddefinition
\bigskip

\remark{Remark}  We can consider $\{\bold s_n:n < \omega\}$, which
hopefully will be independent, i.e. for every $X \subseteq \omega$ for
some c.c.c. forcing notion $\Bbb P$, in $\bold V^{\Bbb P}$ we have
$\lambda \rightarrow (\bold s_n)_\mu$ iff $n \in X$.  It is natural to
try $\{\bold s_n:n < \omega\}$ where 
$\bold s_n = (\text{dom}_{\ell(*),m_n},e_n)$ where $m_n = n$ (or
$2^{2^n}$ may be more convenient) and $e_n$ is [log log$(n)$]-nice.
\endremark
\bigskip

\definition{\stag{2.1} Definition}  [$\lambda$ is our fixed cardinal.]
\nl
1) Let $M^*$ (or $M^*_\lambda$) 
be a model with universe $\lambda$, countable vocabulary,
and its relations and functions are exactly those defined in $({\Cal
H}(\chi),\in,<^*_\chi)$ for $\chi = \lambda^+$ (and some choice of
$<^*_\chi$, a well ordering of ${\Cal H}(\chi)$). \nl
2) For $\bar \alpha \in {}^{\omega >}(M^*_\lambda)$ let 
$c \ell_\ell(\bar \alpha) = \{\beta < \lambda:$ for some first order
$\varphi(y,\bar x)$ we have $M^*_\lambda \models \varphi[\beta,\bar
\alpha] \and (\exists {}^{\le \aleph_\ell} x)\varphi(x,\bar
\alpha)\}$ and $c \ell(\bar \alpha) = \{\beta < \lambda:$ for
some first order $\varphi(y,\bar x)$ we have $M^*_\lambda \models
\varphi[\beta,\bar \alpha] \and (\exists^{\le \aleph_0}
x)\varphi(x,\bar \alpha)\}$. \nl
3) For a model $M$ and $A \subseteq M$ let $c \ell_M(A)$ is the
smallest set of elements of $M$ including $A$ and closed under the
functions of $M$ (so including the individual constants).
\enddefinition
\bn
Note that \nl
\margintag{2.1a}\ub{\stag{2.1a} Fact}:  If $\beta_0,\beta_1 \in 
c \ell_{\ell +1}(\bar \alpha)$ \ub{then}
for some $i \in \{0,1\}$ we have $\beta_i \in c \ell_\ell(\bar \alpha
\char 94 \langle \beta_{1-i} \rangle)$.
\bigskip

\demo{Proof}  Easy.

The idea of our forcing notion is to do historical forcing 
(see \cite{RoSh:733}
for more on historical forcing and its history). That is, we put in only those
conditions which we have to put in order to meet our demands, so every
condition in the forcing has a definite rule of creation. 
In particular, (see below), in the definition of our partial colourings,
we avoid giving the same color to any pairs for which we can afford
this, if the rule of creation is to be respected. We note that the
situation here is not as involved as the one of \cite{RoSh:733}, and we
do not in fact need the actual history of every condition. 
\sn
We proceed to the formal definition of our forcing.

Clearly case 0 for $k \ge 0$ is not necessary from a historical point
of view but it simplifies our treatment later; also case 1 is used
in clause $(\eta)$ of case 3. \nl
Note that in case 2 below we do not require that the conditions are
isomorphic over their common part (which is natural for
historic forcing) as the present choice 
simplifies clause $(\zeta)(iv)$ in case 3.
\enddemo
\bigskip

\definition{\stag{2.2} Main Definition}   Let $n(\ast) \ge 2,n(*) \le
\ell(*) < \omega,\lambda = \aleph_{n(*)},\mu = \aleph_0$ be fixed. 
All closure operation we shall use are understood to refer to 
$M^*_{\aleph_{n(*)}}$ from \scite{2.1}(2).  
Let $\Gamma \subseteq ID^*_{\ell(*)}$ be given.
For two sets $u$ and $v$ of ordinals with $|u| =|v|$, we
let $OP_{v,u}$ stand for the unique order preserving 1-1 
function from $u$ to $v$.

We shall define $\Bbb P =: \Bbb P_\Gamma = 
\Bbb P^\lambda_\Gamma$, it is $\subseteq \Bbb P^*_\lambda$.

Members of $\Bbb P^*_\lambda$ are the pairs of 
the form $p=(u,c) =: (u^p,c^p)$ with

$$
u \in [\lambda]^{<\aleph_0} \text{ and } c:[u]^2 \rightarrow \omega.
$$
\mn
The order in $\Bbb P^*_\lambda$ is defined by

$$
(u_1,c_1) \le (u_2,c_2) \Leftrightarrow (u_1 \subseteq u_2
\and c_1=c_2 \restriction [u_1]^2) 
$$
\mn
For $p \in \Bbb P^*_\lambda$ let $n(p) = \sup(\text{Rang}(c^p)) +
1$; this is $< \omega$. \nl
We now say which pairs $(u,c)$ of the above form (i.e. $(u,c) \in 
\Bbb P^*_\lambda$)
will enter $\Bbb P$. We shall have $\Bbb P= \dbcu_{k < \omega} 
\Bbb P_k$ where $\Bbb P_k =: \Bbb P_k^{\lambda,\Gamma}$ are defined 
by induction on $k < \omega$, as follows.
\enddefinition
\bn
\ub{Case $0$}:   \ub{$k = 4 \ell$}.  If $k = 0$ let 
$\Bbb P_0 =: \{(\emptyset,\emptyset)\}$.

If \ub{$k = 4l > 0$}, a pair $(u,c) \in \Bbb P_k$ \ub{iff} for some 
$(u',c') \in \dbcu_{m<k} \Bbb P_m$ we have $u \subseteq u'$ 
and $c = c'\restriction [u]^2$;  we write $(u,c) = (u',c') \restriction u$.
\bn
\ub{Case $1$}:   \ub{$k = 4 \ell +1$}.  (This rule of creation is needed for density
arguments.)

A pair $(u,c)$ is in $\Bbb P_{k}$ \ub{iff} (it belong to $\Bbb
P^*_\lambda$ and) there is a $(u_1,c_1) \in 
\dbcu_{m<k} \Bbb P_m$ and $\alpha < \lambda$ satisfying $\alpha \notin
u_1$  such that:
\mr
\item "{$(a)$}"   $u=u_1\cup\{\alpha\}$,
\sn
\item "{$(b)$}"  $c \restriction [u_1]^2=c_1$ and 
\sn
\item "{$(c)$}"    For every $\{\beta,\gamma\}$ and 
$\{\beta',\gamma'\}$ in $[u]^2$ which are not equal, if
$c(\{\beta,\gamma\})$ and $c(\{\beta',\gamma'\})$ are equal, \ub{then}
$\{\beta,\gamma\},\{\beta',\gamma'\}\in [u_1]^2.$
(Hence, $c$ does not add any new equalities except for those already
given by $c_1$.)
\endroster
\bn
\ub{Case $2$}:  \ub{$k = 4 \ell + 2$}.   
(This rule of creation is needed for free amalgamation, 
used in the $\Delta$-system arguments for the proof of the $c.c.c.$.)

A triple $(u,c)$ is in $\Bbb P_{k}$ \ub{iff} (it belongs to $\Bbb
P^*_\lambda$ and) there are $(u_1,c_1),
(u_2,c_2) \in \dbcu_{m<k} \Bbb P_m$ for which we have
\mr
\item "{$(a)$}"   $u= u_1\cup u_2$
\sn
\item "{$(b)$}"   $c \restriction [u_1]^2 = c_1$ and $c \restriction
[u_2]^2=c_2$
\sn
\item "{$(c)$}"   $c$ does not add any unnecessary equalities, i.e.,
if $\{\beta,\gamma\}$ and $\{\beta',\gamma'\}$ are distinct and in $[u]^2$ and
$c(\{\beta,\gamma\})=c(\{\beta',\gamma'\})$, \ub{then}
$\{\{\beta,\gamma\},\{\beta',\gamma'\}\} \subseteq [u_1]^2\cup
[u_2]^2$. \nl 
Note that $[u_1]^2 \cap [u_2]^2 = [u_1 \cap u_2]^2$
\sn
\item "{$(d)$}"   $c \ell_0(u_1 \cap u_2) \cap (u_1 \cup u_2) 
\subseteq u_1$ (usually $c \ell_0(u_1 \cap u_2) \cap (u_1 \cup u_2)
\subseteq u_1 \cap u_2$) is O.K. too for 
present \S2, \S3 but not, it seems, in \scite{4.6}).
\endroster
\bn
\ub{Main rule}:
\sn
\ub{Case $3$}:  \ub{$k = 4 \ell +3$}.   (This rule \footnote{you may
understand it better seeing how it is used in the proof of
\scite{3.2}} is like the
previous one, but the amalgamation is taken over a graph 
$\bold s = (\text{dom}_{\ell(*),m},e) \in \Gamma)$.

A pair $(u,c) \in \Bbb P_k$ iff there are
$\bold s = (\text{dom}_{\ell(*),m(*)},e) \in \Gamma$ and a sequence of conditions
$$
\bar p = \langle p_y:\,y \in Y \rangle \text{ where } Y = 
\{y \in [\text{dom}_{\bold s}]^2:|y/e| > 1\}
$$
\mn
from $\dbcu_{m<k} \Bbb P_m$ 
AND we have a sequence of finite sets 
$\bar v = \langle v_t:t \in Y^+ \rangle$ where

$$
Y^+ = \{t:t \in Y \text{ or } t = \emptyset \text{ or } t =
\{\eta\}, \text{ where } \eta \in \text{ dom}_{\bold s}\}
$$
\sn
such that
\mr
\item "{$(a)$}"   $u = \bigcup\{ u^{p_y}:\,y \in Y\}$
\sn
\item "{$(b)$}"  $(u,c) \in \Bbb P^*_\lambda$ and 
$c \restriction [u^{p_y}]^2= c^{p_y}$ for all $y \in Y$
\sn
\item "{$(c)$}"  if $\alpha_1 \ne \alpha_2,\beta_1 \ne
\beta_2$ are from $u$ and $\{\alpha_1,\alpha_2\} \ne \{\beta_1,\beta_2\}$ and
$c\{\alpha_1,\alpha_2\} = c\{\beta_1,\beta_2\}$ \ub{then} $(\exists
y)[\{\alpha_1,\alpha_2\} \subseteq u^{p_y}]$ and $(\exists
y)[\{\beta_1,\beta_2\} \subseteq u^{p_y}]$
\sn
\item "{$(d)$}"   $v_t \cap v_s \subseteq v_{t \cap s}$ for
$t,s, \in Y^+$ 
\sn 
\item "{$(e)$}"  $c \ell_0(v_t) \cap u^{p_y} \subseteq 
v_t$ for all $y \in Y$ and 
$t \in \{\emptyset\} \cup \{\{\eta\}:\eta \in \text{ dom}_{\bold s}\}$
\sn
\item "{$(f)$}"   $u^{p_y} \subseteq v_y$ for all $y \in Y$
\sn
\item "{$(g)$}"   if $y_1,y_2 \in Y$ and $t \in \{\emptyset\} \cup
\{\{\eta\}:\eta \in \text{ dom}_{\bold s}\}$ and $t = y_1 \cap y_2$,
\ub{then} $p_{y_1} \restriction v_t = p_{y_2} \restriction v_t$;
equivalently: $\{p_\eta:\eta \in Y\}$ has a comon upper bound in $\Bbb
P^*$.
\endroster
\bigskip

\proclaim{\stag{2.4} Claim}  1) $\Bbb P^\lambda_\Gamma$ satisfies the
c.c.c. and even the Knaster condition. \nl
2) For each $\alpha < \lambda$ the set ${\Cal J}_\alpha = \{p \in \Bbb
P^\lambda_\Gamma:\alpha \in u^p\}$ is dense open. \nl
3) $\Vdash_{{\Bbb P}^\lambda_\Gamma} ``\underset\tilde {}\to c =
\cup\{c^p:p \in \underset\tilde {}\to G\}$ is a function from
$[\lambda]^2$ to $\omega"$. 
\endproclaim
\bigskip

\demo{Proof}  1) By Case 2. 

In detail, assume that $p_\varepsilon \in \Bbb P^\lambda_\Gamma$ for
$\varepsilon < \omega_1$ and let $p_\varepsilon =
(u_\varepsilon,c_\varepsilon)$.  As each $u_\varepsilon$ is a finite
subset of $\lambda$, by the $\Delta$-system lemma \wilog \, for some
finite $u^* \subseteq \lambda$ we have: if $\varepsilon < \zeta <
\omega_1$ then $u_\varepsilon \cap u_\zeta = u^*$  By further
shrinking, \wilog \, $\alpha \in u^* \Rightarrow \langle
|u_\varepsilon \cap \alpha|:\varepsilon < \omega_1 \rangle$ is
constant and $\varepsilon < \zeta < \omega_1 \Rightarrow
|u_\varepsilon| = |u_\zeta|$.  Also \wilog \, the set $\{(\ell,m,k)$:
for some $\alpha \in u_\varepsilon$ and $\beta \in u_\varepsilon$ we
have $\ell = |\alpha \cap u_\varepsilon|,m = |\beta \cap
u_\varepsilon|$ and $k = c_\varepsilon\{\alpha,\beta\}\}$ does not
depend on $\zeta$.  We can conclude that $\varepsilon < \zeta <
\omega_1 \Rightarrow \text{ OP}_{u_\zeta,u_\varepsilon}$ maps
$p_\varepsilon$ to $p_\zeta$.  Clearly for $\varepsilon < \omega_1$,
the set $c \ell(u_\varepsilon)$ is countable hence for every $\zeta <
\omega_1$ large enough we have $u_\zeta \cap c \ell_0(u_\varepsilon)
\in u^*$ so restricting $\langle p_\varepsilon:\varepsilon < \omega_1
\rangle$ to a club we get that $\varepsilon < \zeta < \omega_1
\Rightarrow c \ell_0(u_\varepsilon) \cap u_\zeta = u^*$ (this is much
more than needed).  Now for any $\varepsilon < \zeta < \omega_1$ we
can define $q_{\varepsilon,\zeta} =
(u_{\varepsilon,\zeta},c_{\varepsilon,\zeta})$ with
$u_{\varepsilon,\zeta} = u_\varepsilon \cup u_\zeta$ and
$c_{\varepsilon,\zeta}:[u_{\varepsilon,\zeta}]^2 \rightarrow \omega$
is defined as follows: for $\alpha < \beta$ in $u_{\varepsilon,\zeta}$
let $c_{\varepsilon,\zeta}\{\alpha,\beta\}$ be
$c_\varepsilon\{\alpha,\beta\}$ if defined, $c_\zeta\{\alpha,\beta\}$
if defined, and otherwise sup$\{(\text{Rang}(c_\varepsilon)) + 1 +
(|u_{\varepsilon,\zeta} \cap \alpha| + |u_{\varepsilon,\zeta} \cap
\beta|)^2 + |u_{\varepsilon,\zeta} \cap \alpha|$.  Now $q \in \Bbb
P^\lambda_\Gamma$ by case 2, and $p_\varepsilon \le
q_\varepsilon,p_\zeta \le q_{\varepsilon,\zeta}$ by the definition of 
order. \nl
2) By Case 1. 

In detail, let $p \in \Bbb P^\lambda_\Gamma$ and $\alpha < \lambda$
and we shall find $q$ such that $p \le q \in {\Cal I}_\alpha$.  If
$\alpha \in u^p$ let $q=p$, otherwise define $q = (u^q,c^q)$ as
follows $u^q = u^p \cup \{\alpha\}$ and for $\beta < \gamma \in u^q$
we let $c^q\{\beta,\gamma\}$ be: $c^p\{\beta,\gamma\}$ when it is well defined
and sup(Rang$(c^p)) + 1 + (|\beta \cap u^q| + |\gamma \cap u^q|)^2
+ |\beta \cap u^q|$ when otherwise.  Now $q \in \Bbb P^\lambda_\Gamma$ by
case 1 of Definition \scite{2.2}, 
$p \le q$ by the order's definition and $q \in {\Cal
I}_\sigma$ trivially.  
\nl
3) Follows from part (2).  \hfill$\square_{\scite{2.4}}$\margincite{2.4}
\enddemo
\newpage

\head {\S3 Why does the forcing work} \endhead  \resetall \sectno=3
\bigskip

We shall use the following claim for $\mu = \aleph_0$
\proclaim{\stag{3.1} Claim}  1) If $f:[\lambda]^2 \rightarrow \mu$ and $M$
is an algebra with universe $\lambda,|\tau_M| \le \mu$ and 
$w_t \subseteq [\lambda],|w_t| < \aleph_0$ for $t \in
[\lambda]^2$ and $\lambda \ge \beth_2(\mu^+)^+$, \ub{then} for some
$\langle v_t:t \in [W]^{\le 2} \rangle$ we have:
\mr
\item "{$(a)$}"  $W \subseteq \lambda$ is infinite in fact $|W| = \mu^+$
\sn
\item "{$(b)$}"  $f \restriction [W]^2$ is constant
\sn
\item "{$(c)$}"  $t \cup w_t \subseteq v_t \in [\lambda]^{< \aleph_0}$ 
for $t \in [W]^2$
\sn
\item "{$(d)$}"  $v_{t_1} \cap v_{t_2} \subseteq v_{t_1 \cap t_2}$ when
$t_1,t_2 \in [W]^{\le 2}$ but for no $\alpha < \beta < \gamma$ do we
have $\{t_1,t_2\} = \{\{\alpha,\beta\},\{\beta,\gamma\}\}$
\sn
\item "{$(e)$}"  if $t_1,t_2 \in [W]^i$, where $i \in \{1,2\}$ then
$|v_{t_1}| = |v_{t_2}|$ and $OP_{v_{t_2},v_{t_1}}$ maps $t_1$ onto
$t_2$ and $w_{t_1}$ onto $w_{t_2}$ and $v_{t_1}$ onto
$v_{t_2},w_{\{\text{Min}(t_1)\}}$ onto $w_{\{\text{Min}(t_2)\}}$ and
$w_{\{\text{Max}(t_1)\}}$ onto $v_{\{\text{Max}(t_2)\}},v_{\{\text{Min}(t_1)\}}$ onto
$v_{\{\text{Min}(t_2)\}}$, and $v_{\{\text{Max}(t_1)\}}$ onto
$v_{\{\text{Max}(t_2)\}}$
\sn
\item "{$(f)$}"  $v_{\{\alpha,\beta\}} \cap c \ell_M(v_{\{\gamma\}})
\subseteq v_\gamma$ for $\alpha,\beta,\gamma \in W$.
\ermn
2) If $[u \in [\lambda]^{< \aleph_0} \Rightarrow c \ell_M(u) \in
[M]^{< \mu}$, \ub{then} $\lambda = (\beth_2(\mu))^+$ is enough.
\endproclaim
\bigskip

\remark{Remark}  See more in \cite{Sh:289}; this is done for completeness.
\endremark
\bigskip

\demo{Proof}  1) Let 
$w_t \cup t = \{\zeta_{t,\ell}:\ell < n_t\}$ with no repetitions 
and we define the function $c,c_0,c_1$ with domain $[\lambda]^3$ 
as follows: if $\alpha < \beta < \gamma < \lambda$ then

$$
c_0\{\alpha,\beta,\gamma\} = \{(\ell_1,\ell_2):\ell_1 <
n_{\{\alpha,\beta\}},\ell_2 < n_{\{\alpha,\gamma\}} \text{ and }
\zeta_{\{\alpha,\beta\},\ell_1} = \zeta_{\{\alpha,\gamma\},\ell_2}\}
$$

$$
c_1\{\alpha,\beta,\gamma\} = \{ (\ell_1,\ell_2):\ell_1 <
n_{\{\alpha,\gamma\}},\ell_2 < n_{\{\beta,\gamma\}} \text{ and }
\zeta_{\{\alpha,\gamma\},\ell_1} = \zeta_{\{\beta,\gamma\},\ell_2}\}
$$

$$
c\{\alpha,\beta,\gamma\} = 
(c_0\{\alpha,\beta,\gamma\},c_1\{\alpha,\beta,\gamma\},f\{\alpha,\beta\}).
$$
\mn
By Erd\"os-Rado theorem for some $W_1 \subseteq \lambda$ of
cardinality and even order type $\mu^{++}$ for part (1), $\mu^+$ for
part (2) such that $c \restriction [W_0]^3$ is constant.  Let
$\{\alpha_\varepsilon:\varepsilon < \mu^{++}\}$ list $W_0$ in
increasing order.  If $2 < i < \mu^{++}$, let

$$
\align
v_{\{\alpha_i\}} =: \{\zeta_{\{\alpha_i,\alpha_{i+1}\},\ell_1}:&\text{
for some } \ell_2 \text{ we have} \\
  &(\ell_1,\ell_2) \in c_0\{\alpha_i,\alpha_{i+1},\alpha_{i+2}\}\}
\cup \\
  &\{\zeta_{\{\alpha_0,\alpha_i\},\ell_1}:\text{ for some } \ell_2
\text{ we have} \\
  &(\ell_1,\ell_2) \in c_1\{\alpha_0,\alpha_1,\alpha_i\}\}
\endalign
$$
\mn
(clearly $\alpha_i \in v_{\{\alpha_i\}})$.

For $i < j$ in $(2,\mu^{++})$ let $v_{\{\alpha_i,\alpha_j\}} =
v_{\{\alpha_i\}} \cup v_{\{\alpha_j\}} \cup
w_{\{\alpha_i,\alpha_j\}}$.  Now for some unbounded $W_2 \subseteq W_1
\backslash \{\alpha_0,\alpha_1\}$ and $Y \in [\lambda]^{\le \mu}$ we
have:

if $\alpha \ne \beta \in W_2$ then $c \ell_M(v_{\{\alpha\}}) \cap c
\ell_M(v_{\{\beta\}}) \subseteq Y$.

Now by induction on $\varepsilon < \mu^+$ we can choose
$\gamma_\varepsilon \in W_2$ strictly increasing with
$\varepsilon,\gamma_\varepsilon$ large enough.  It is easy to check
that $W = \{\gamma_\varepsilon:\varepsilon < \mu^+\}$ is as
required. \nl
2) The same proof.  \hfill$\square_{\scite{3.1}}$\margincite{3.1}
\enddemo
\bigskip

\proclaim{\stag{3.1a} Claim}  Let $n(*),\ell(*),\lambda$ be as in
Definition \scite{2.2}, and see Definition \scite{2.0}.  
Assume that $\Gamma_1,\Gamma_2,m,p$ satisfies:
\mr
\item "{$(a)$}"  $\Gamma_1,\Gamma_2 \subseteq ID^*_{\ell(*)}$
\sn
\item "{$(b)$}"  if (${\text{\rm dom\/}}_{\ell(*),m},e) \in ID^*_{\ell(*)}$ and
$({\text{\rm dom\/}}_{\ell(*),m},e)$ 
is not $m^*$-nice \ub{then} $({\text{\rm dom\/}}_{\ell(*),m},e) \in
\Gamma_1 \Leftrightarrow ({\text{\rm dom\/}}_{\ell(*),m},e) \in \Gamma_2$
\sn
\item "{$(c)$}"  $p^* \in \Bbb P^*_\lambda$ and $|u^{p^*}| < m^*$.
\ermn
\ub{Then} $p^* \in \Bbb P^\lambda_{\Gamma_1} \Leftrightarrow p^* \in
\Bbb P^\lambda_{\Gamma_2}$.
\endproclaim
\bigskip

\demo{Proof}  We prove by induction on $k < \omega$ that
\mr
\item "{$(*)_k$}"  if $r' \in \Bbb P^{\lambda,\Gamma_1}_k$ (see
Definition \scite{2.2} before case 0) and $r \le
r'$ and $|u^r| < m^*$, \ub{then} $r \in \Bbb P^\lambda_{\Gamma_2}$.  
\ermn
This is enough by the symmetry in our assumptions. \nl
For a fixed $k$ we prove this by induction on $|u^r|$.
The proof splits according to the case in Definition \scite{2.2} which
hold for $r'$.
\enddemo
\bn
\ub{Case 0}:  

Trivial.
\bn
\ub{Case 1}:

Easy.
\bn
\ub{Case 2}:  
  
Should be clear but let us check, so $r' = (u',c')$ is gotten from
$(u'_1,c'_1),(u'_2,c'_2)$ as in clauses (a)-(d) of case 2, and let $r =
(u,c) \le r'$. \nl
Let $u_\ell = u'_\ell \cap u,c_\ell = c'_\ell \restriction
[u_\ell]^2$.  By the induction hypothesis $(u_\ell,c_\ell) \in \Bbb
P^\lambda_{\Gamma_2}$ and it is enough to check that
$(u,c),(u_1,c_1),(u_2,c_2)$ are in case (2) of Definition \scite{2.2}
which is easy, e.g. in
clause (d) we use monotonicity of $c \ell_0$.
\bn
\ub{Case 3}: 

So let $r'$ be gotten from $\bold s = (\text{dom}_{\ell(*),m},e),
\langle p_y:y \in Y \rangle,\langle v_t:t \in Y^+ \rangle$ as
there.  Of course, we have $(\text{dom}_{\ell(*),m},e) \in \Gamma_1$ 
and $p_y \in \dbcu_{\ell < k} \Bbb P^{\lambda,\Gamma_1}_\ell$ so 
by the induction hypothesis clearly $p_y \restriction u^r 
\in \Bbb P^\lambda_{\Gamma_2}$.
\bn
\ub{Subcase 3A}:  nice(dom$_{\ell(*),m(*)},e) < m^*$ (see Definition
\scite{1.3}(1)).  \nl
Hence (dom$_{\ell(*),m},e) \in
\Gamma_2$ and the desired conclusion easily holds. \nl
[Why?  We can find $p^*_y = p_y \restriction u^r = p_y \restriction
(u^{p_y} \cap u^r) = r' \restriction (u^{p_y} \cap u^r)$ 
hence $|u^{p^*_y}| < m^*$.
\nl
By the induction hypothesis $p^*_y$ belongs to $\Bbb
P^\lambda_{\Gamma_2}$ for each $y \in Y$.  Now $r,\langle p^*_y:y \in
Y \rangle$ and $\langle v_t:t \in Y^+ \rangle$ satisfies clauses (a) -
(g) of case 3 of Definition \scite{2.2}.  Hence 
by case 3 of Definition \scite{2.2} $r'' =: 
r' \restriction (\cup\{u^{p^*_y}:y \in Y\})$ belong to $\Bbb
P^\lambda_{\Gamma_2}$ but $r = r''$ so $r \in \Bbb P^\lambda_{\Gamma_2}$.]
\bn
\ub{Subcase 3B}:  Not subcase 3A.

So nice(dom$_{\ell(*),m},e) \ge m^* > |u^r|$.  For $a \subseteq
\text{ dom}_{\ell(*),m}$ let $u_a = \{\alpha \in u^r:\alpha \in
u^{p_y}$ for some $y \in Y$ satisfies $y \subseteq a$ \ub{or} $\alpha \in
v_{\{\eta\}},\eta \in a$ \ub{or} $\alpha \in v_\emptyset\}$.  Now
\mr
\item "{$(*)_0$}"  if $u^r \subseteq v_{\{\eta\}}$ for some $\eta \in
\text{ dom}_{\ell(*),m}$ then $r \in \Bbb P^\lambda_{\Gamma_2}$. \nl
[Why?  By applying case 2 (and 0) of Definition \scite{2.2}.]
\sn
\item "{$(*)_1$}"  if for some $\eta \in \text{ dom}_{\ell(*),m}$ we
have $[(\{\alpha_1,\beta_1\} \ne \{\alpha_2,\beta_2\} \in
[u]^2) \and c^r\{\alpha_1,\beta_1\} = c^r\{\alpha_2,\beta_2\} \Rightarrow
\{\alpha_1,\beta_1,\alpha_2,\beta_2\} \subseteq v_{\{\eta\}}]$
\ub{then} $r \in \Bbb P^\lambda_{\Gamma_2}$. \nl
[Why?  By $(*)_0$ and uses of case 1 of definition \scite{2.2}.]
\sn
\item "{$(*)_2$}"  if $y \in Y$ and $u^r \subseteq v_y$ then $r \in
\Bbb P^\lambda_{\Gamma_2}$. \nl
[Why?  Similarly.]
\ermn
Now
\mr
\item "{$(*)_3$}"  It is enough to 
find $a,b \subseteq \text{ dom}_{\ell(*),m}$ such
that:
{\roster
\itemitem{ $(*)^3_{a,b}$ }   $u_a \ne u_a \cap u_b,u_b \ne u_a \cap u_b,u^r
\subseteq u_a \cup u_b,u^r \nsubseteq u_a,u^r \nsubseteq u_b$ 
and $[\eta_1 \in u_a \backslash u_b \and \eta_2 \in
u_b \backslash u_a \Rightarrow (\{\eta_1,\eta_2\}/e)$ is a singleton]
and $|a \cap b| \le 1$.
\endroster}
\ermn
[Why is this enough? 
As then $r$ is gotten by Case 2 of Definition \scite{2.2} from
$(u_a,c^p \restriction [u_a]^2),(u_b,c^b \restriction [u_b]^2)$.  The
main point is why clause (d) of this case holds; now we shall prove
more  $c \ell_0(u_a \cap u_b) \cap (u_a \cup u_b) \subseteq u_a \cap u_b$; now
by clause $(e)$ of case 3 of Definition \scite{2.2} letting
$t = a \cap b$ (it $\in \{\emptyset\} \cup \{\{\eta\}\}:\eta \in \text{
dom}_{\ell(*),m}\}$ by the last statement in $(*)^3_{a,b}$)
we have $u_a \cap u_b = u_t$ (see $(d),(f)$ Definition
\scite{2.2}, Case 3), hence $c \ell_0(u_a \cap u_b) = c \ell_0(u_t),u_t
\subseteq  v_t$ hence $c \ell_0
(c \ell_0(u_a \cap u_b) \subseteq c \ell_0(v_{a \cap b})$ 
which is disjoint to $u_a \backslash u_b$ and to $u_b
\backslash u_a$ by clause $(e)$ in case 3 of Definition \scite{2.2}
as $u_a \cap v_t = u_t$ and $u_b \cap v_t = u_t$.]
\nl
So now why can we find such $a,b$?

We try to choose $a_i \subseteq \text{ dom}_{\ell(*),m}$ for $i =
2,3,\ldots$ or for $i = 1,2,\ldots,$ such that $|a_i| = i,a_i \subseteq
a_{i+1}$ and $i \le |u_{a_i}|$.  First assume that we cannot find neither $a_2$
nor $a_1$, then $y \in Y \Rightarrow |u^{p_y} \cap u^r| \le 1$ and
$\eta \in \text{ dom}_{\ell(*),m} \Rightarrow |v_{\{\eta\}} \cap u^r| = 0$.
If $(*)_2$ applies we are done, so there are $\langle
(y_\ell,\gamma_\ell):\ell < k \rangle$ satisfying $y_\ell =
\{\eta_{1,\ell},\eta_{2,\ell}\} \in Y$ such that $u^r \cap u^{p_y}
\backslash v_{\{\eta_{1,\ell}\}} \backslash v_{\{\eta_{2,\ell}\}} =
\{\gamma_\ell\}$ and $k \ge 2$ so $u^r \backslash v_\emptyset =
\{\gamma_\ell:\ell < k\}$.  Let $u_1 = (u^r \cap v_\emptyset) \cup
\{\gamma_0\},u_2 = u^r \backslash \{\gamma_0\}$, clearly $r$ is gotten
from $r \restriction u_1,r \restriction u_2$ as in case 2 of
Definition \scite{2.2}.

Second, assume $a_1$ or $a_2$ is defined.
So we are stuck in $a_{i(*)}$ for some $i(*)$, i.e. $a_{i(*)}$ is
chosen but we cannot choose $a_{i(*)+1}$.  If $u_{a_{i(*)}} \ne u^r$, let
$a = a_{i(*)},b = \text{ dom}_{\ell(*),m} \backslash a_{i(*)}$, so we
get $(*)^3_{a,b}$ and we are done.  So $u^r = u_{a_{i(*)}}$, hence $i(*) =
|a_{i(*)}| = |u^r| < m^*$ and we can assume that $(*)_2$ does not
apply.  
By the niceness of (dom$_{\ell(*),m},e)$
the graph $H[e] \restriction a_{i(*)}$ has no cycle so is a tree in
the graph theoretic sense and so for some $c,b
\subseteq a_{i(*)}$ we have $c \cap b = \{\eta\},c \cup b = a_{i(*)},b
\ne \{\eta\},c \ne \{\eta\}$
and $[\eta' \in b \backslash \{\eta\} \and \eta'' \in c \backslash
\{\eta\} \Rightarrow \{\eta',\eta''\}$ not an $H[e]$-edge]; so we get
$(*)^3_{a,b}$ and we are done.
(So if we change slightly the claim demanding only $2|u^r| < m^*$, the proof
is simpler).    \hfill$\square_{\scite{3.1a}}$\margincite{3.1a}
\bn
\centerline{$* \qquad * \qquad *$}  
\bigskip

\proclaim{\stag{3.2} The preservation Claim}  Let
$n(*),\ell(*),\lambda$ be as in Definition \scite{2.2} and assume
$\lambda > \beth_2(\mu^+)$. \nl
1) If $\Bbb P = \Bbb P^\lambda_\Gamma$ and 
$({\text{\rm dom\/}}_{\ell(*),m},e) \in \Gamma \subseteq
ID^*_{\ell(*)}$ \ub{then} in $\bold V^{\Bbb P}$ we have
$({\text{\rm dom\/}}_{\ell(*),m},e) \in ID_2(\lambda,\aleph_0)$. \nl
2) Assume that $\Bbb P = \dsize \prod_{n < \gamma} \Bbb
P^\lambda_{\Gamma_n}$ where $\Gamma_n \subseteq ID^*_{\ell(*)}$ and
$\gamma \le \omega$ and $p^* \in \Bbb P$ forces that $\underset\tilde
{}\to c$ is a function from $[\lambda]^2$ to $\omega$.
\ub{Then} for some finite $d \subseteq \gamma$ for any  
$\bold s \in \dbca_{n \in d} \Gamma_n$ we have $p^* \nVdash_{\Bbb P} ``\bold s
\notin ID_2(\underset\tilde {}\to c)"$.
\endproclaim
\bigskip

\demo{Proof}  1) Follows from (2), letting $\gamma =1,\Gamma_0 =
\Gamma$. \nl
2) Assume $p^* \in \Bbb P$ and $p^*
\Vdash_{\Bbb P} ``\underset\tilde {}\to c$ is a function from 
$[\lambda]^2$ to $\omega$". 
Let $k(*) = 2^{\ell(*)} -1$ and let 
$k(\nu) = |\{\rho \in {}^{\ell(*) >} 2:
\rho <_{\text{lex}} \nu\}|$ for $\nu \in {}^{\ell(*) >}2$.  For $p \in
\Bbb P$ let $u[p] = \cup\{u^{p(n)}:n \in \text{ Dom}(p)\}$, so $u[p]
\in [\lambda]^{< \aleph_0}$ and for any $q \in \Bbb P$ we let
$n[q] = \sup (\cup\{\text{Rang}(c^{q(n)}):n
\in \text{ Dom}(q)\})$.
For any $\alpha < \beta < \lambda$ letting $t = \{\alpha,\beta\}$ we
define, by induction on $k \le k(*)$ the triple $(n_{t,k},w_{t,k},d_{t,k})$
such that:
\mr
\item "{$(*)$}"  $n_{t,k} < \omega,w_{t,k} \in [\lambda]^{< \aleph_0}$ and
$d_{t,k} \subseteq \gamma$ is finite.
\ermn
\ub{Case $1$}: $k=0:n_{t,k} = n[p^*] +2$ and $w_{t,k} = \{\alpha,\beta\}
\cup u^{p^*}$ and $d_{t,k} = \text{ Dom}(p^*)$.
\mn
\ub{Case $2$}:  $k+1$:

Let ${\Cal P}_{t,k} = \{q \in \Bbb P:p^* \le q,u[q] \subseteq w_{t,k}$
and $n[q] \le n_{t,k} \text{ and Dom}(q) \subseteq d_{t,k}\}$; clearly
it is a finite set, and for every $q \in {\Cal P}_{t,k}$ we
choose $p_{t,q}$ such that $q \le p_{t,q} \in \Bbb P$ and $p_{t,q}$
forces a value, say $\zeta_{t,q}$ to $\underset\tilde {}\to c(t)$.  Now we
let 

$$
w_{t,k+1} = w_{t,k} \cup \bigcup\{u[p_{t,k}]:q \in {\Cal P}_{t,k}\}.
$$

$$
d_{t,k+1} = d_{t,k} \cup \{\text{Dom}(q_{t,p}):p \in {\Cal P}_{t,k}\}
$$

$$
n_{t,k+1} = \text{ Max}\{|w_{t,k+1}|^2,n_{t,k} +1,n[p_q]+1:q \in {\Cal P}_{t,k}\}
$$
\mn
We next define an equivalence relation $E$ on $[\lambda]^2:t_1 E t_2$
\ub{iff} letting $t_1 = \{\alpha_1,\beta_1\},t_2 =
\{\alpha_2,\beta_2\},\alpha_1 < \beta_1,\alpha_2 < \beta_2$ and
letting $h =
OP_{w_{\{\alpha_2,\beta_2\},k(*)},w_{\{\alpha_1,\beta_1\},k(*)}}$, we
have
\mr
\widestnumber\item{$(iii)$}
\item "{$(i)$}"  $w_{t_1,k(*)},w_{t_2,k(*)}$ has the same number of
elements
\sn
\item "{$(ii)$}"  $h$ maps $\alpha_1$ to $\alpha_2$ and $\beta_1$ to
$\beta_2$ and $w_{t_1,k}$ onto $w_{t_2,k}$ for $k \le k(*)$ (so $h$ is onto)
\sn
\item "{$(iii)$}"  $d_{t_1,k} = d_{t_2,k}$ for $k \le k(*)$ \nl
(hence $h$ maps ${\Cal P}_{t_1,k}$ onto ${\Cal P}_{t_2,k})$
\sn
\item "{$(iv)$}"  if $q_1 \in {\Cal P}_{t_1,k},k < k(*)$ then $h$ maps
$q_1$ to some $q_2 \in {\Cal P}_{t_2,k}$ and it maps
$p_{t,q_1}$ to $p_{t,q_2}$ and we have $\zeta_{t_1,q_1} = \zeta_{t_2,q_2}$.
\ermn
Clearly $E$ has $\le \aleph_0$ equivalence classes.  So let
$c:[\lambda]^2 \rightarrow \aleph_0$ be such that $c(t_1) = c(t_2)
\Leftrightarrow t_2 E t_2$ and let $w_t = w_{t,k(*)}$.

By Claim \scite{3.1}, recalling that we have assumed $\lambda >
\beth_2(\aleph_1)$ we can find $W \subseteq \lambda$ of
cardinality $\mu^+$ and
$\bar v = \langle v_t:t \in [W]^{\le 2} \rangle$ as there; i.e., we
apply it to an expansion of $M^*_\lambda$ such that $c \ell_0(-) = c
\ell_M(-)$. 

Let $d^*_k = d_{t,k} \subseteq \omega$ for $t \in [W]^2$, now we choose $d =
d^*_{k(*)} \subseteq \gamma$, and we shall show that it is as required
in the claim.
Let $\bold s = (\text{dom}_{\ell(*),m(*)},e) \in \dbca_{\ell \in d}
\Gamma_\ell$ and let $Y_\nu =
Y_{e,\nu} = \{\{\eta_0,\eta_1\}:\nu \char 94 <i> \trianglelefteq
\eta_i \in \text{ dom}_{\bold s}$ for $i=0,1$ and
$\{\eta_0,\eta_1\}/e$ is not a singleton$\}$ for 
$\nu \in {}^{\ell(*) >} 2$ and let $Y = \cup\{Y_{\ell,\nu}:\nu \in
{}^{\ell(*) >} 2\}$.

We now choose $\alpha_\eta \in W$ for $\eta \in \text{
dom}_{\bold s}$ such that $\eta_1 <_{\text{lex}} \eta_2
\Rightarrow \alpha_{\eta_1} < \alpha_{\eta_2}$.  Let $S =
\{\alpha_\eta:\eta \in \text{ dom}_{\bold s}\}$.  For $y \in Y$
let $t(y) = \{\alpha_\eta:\eta \in y\}$.  Let $\langle
\nu^*_\ell:\ell < k(*) = 2^{\ell(*)}-1 \rangle$ list ${}^{\ell(*)>} 2$
in increasing order by $\le_{\ell ex}$.

We now define $q_\ell$ and $q_{\eta,\ell}$ 
for $\eta \in \text{ dom}_{\ell(*),m}$ and
$p_{y,\ell}$ for $y \in Y$ by induction on $\ell \le k(*)$ such that
\mr
\item "{$(a)$}" $p_{y,\ell} \in {\Cal P}_{t(y),\ell}$ hence
$u^{p_{y,\ell}} \subseteq w_{t,\ell}$ for every $y \in Y$
\sn
\item "{$(b)$}" $\Bbb P \models ``p_{y,m} \le p_{y,\ell}"$ for $m \le
\ell$
\sn
\item "{$(c)$}"  if $y \in Y$ and $\eta \in Y$ then 
Dom$(q_{\eta,\ell}) = \text{ Dom}(p_{y,\ell})$ and for each $\beta \in
\text{ Dom}(q_{\eta,\ell})$ we have $q_{\eta,\ell}(\beta) =
p_{y,\ell}(\beta) \restriction v_t$ hence $m \le \ell \Rightarrow 
q_{\eta,m} \le q_{\eta,\ell}$ and Dom$(q_\ell) = \text{
Dom}(q_{\eta,\ell}) \cap v_\emptyset$ and for each $\beta \in \text{
Dom}(q_\ell)$ we have $(q_{\eta,\ell}(\beta) \restriction V_\emptyset
- q_\ell(\beta)$ so $m < \ell \Rightarrow q_m \le q_\ell$
\sn
\item "{$(d)$}"  $\nu^*_\ell \char 94 \langle i \rangle
\trianglelefteq \eta_i \in \text{ dom}_{\ell(*),m(*)}$ for $i = 0,1$
then $p_{y,\ell +1}$ forces a value to $\underset\tilde {}\to
c\{\alpha_\eta:\eta \in y\}$.
\ermn
For $\ell=0$ there is no problem.  For $\ell +1$ choose
$\eta^\ell_0,\eta^\ell_1$ such that $\nu^*_\ell \char 94 \langle i
\rangle \trianglelefteq \eta^\ell_i 
\in \text{ dom}_{\ell(*),m}$ for $i=0,1$ and
$\{\eta^\ell_0,\eta^\ell_1\}/e^{\bold s}$ is not a singleton and let
$y_\ell = \{\eta_0,\eta_1\}$.  As $p^\ell_{y_\ell} \in 
{\Cal P}_{t(y_\ell),\ell}$ by the choice of ${\Cal P}_{t(y),\ell +1}$ there is
$p^\ell_{y_\ell} \in {\Cal P}_{t(y_\ell),\ell +1}$ above
$p_{y_\ell,\ell}$, which forces a value to $\underset\tilde {}\to
c(t(y_\ell))$.  For $p \in \Bbb P$ and $u \subseteq \lambda$ let $q=p
\upharpoonleft u$ means Dom$(p) = \text{ Dom}(q)$ and $\beta \in
\text{ Dom}(p) \Rightarrow q(p) = (p(\beta)) \restriction u$.

Now we define $\langle p_{y,\ell +1}:y \in Y_{\nu_i} \rangle$: 

if $y \in Y_{\nu_i}$ then $p_{y,\ell +1} = \text{
OP}_{v_y,v_{y_\ell}}(p^\ell_{y_\ell})$.

So necessarily
\mr
\item "{$(*)_1$}"   if $y'' \ne y'' \in Y_{\nu_\ell}$ then $p_{y',\ell
+1} \upharpoonleft v_\emptyset  = p_{y'',\ell +1} \upharpoonleft
v_\emptyset$ is above (by $\le_{\Bbb P}$) $q_\ell$
\sn
\item "{$(*)_2$}"  if $y' \ne y'' \in Y_{\nu_\ell}$ and $y' \cap y''
\ne \emptyset$ then for some $\eta \in \text{ dom}_{\bold s}$ we have
$y' \cap y'' = \{\eta\}$ and $p_{y',\ell +1} \upharpoonleft v_{\{\eta\}} =
p_{y'',\ell +1} \upharpoonleft v_{\{\eta\}}$ is above $q_{\eta,\ell}$ \nl
[Why?  As if let $y' = \{\eta'_0,\eta'_1\},y'' =
\{\eta''_0,\eta''_1\},\nu_\ell \char 94 \langle i \rangle
\trianglelefteq \eta'_i,\eta''_i$ for $i=0,1$ then either $\eta'_0 =
\eta''_0,y' \cap y'' = \{\eta'_0\}$ or $\eta'_1 = \eta''_1,y' \cap y''
= \{\eta'_1\}$ but $\eta'_0 \ne \eta''_1 \and \eta'_1 \ne
\eta''_0$. Now use the properties from \scite{3.1} and clause (iv) above.]
\ermn
Let $q_{\emptyset,\ell +1} = p^{\ell +1}_{y_\ell} \upharpoonleft
v_\emptyset$.  The $q_{\eta'_i,\ell}$ is defined as
$q_{\eta'_i,\ell+1} = p_{\{\eta'_0,\eta'_1\},\ell+1}
\upharpoonleft v_{\{\eta'_i\}}$ for $i=0,1$ if 
$\nu_\ell \char 94 \langle i \rangle
\trianglelefteq \eta'_i \and\{\eta'_0,\eta'_1\} \in Y_{\nu_\ell}$.

Let $q_{\eta,\ell +1}$ be the result of free amalgamation
(i.e. case 2 of Definition \scite{2.2}) in each coordinate $\beta$ of 
$q_{\eta,\ell}$ and $q_{\emptyset,\ell +1}$ \ub{if} 
$\eta \in \text{ dom}_{\bold s} \wedge \neg(\nu_\ell \trianglelefteq
\eta)$ and $\eta \in \text{ dom}_{\bold s}$.

Let $p_{y,\ell +1}$ be the result of free amalgamation
(i.e. case 2 of Definition \scite{2.2}) in each coordinate 
(twice) of $p_{y,\ell},q_{\{\eta'_0\},\ell +1},q_{\{\eta_1\},\ell+1}$ 
\ub{if} $y = \{\eta_0,\eta_1\} \in Y \backslash Y_{\nu_\ell}$.

Of course, putting two conditions together using case 2 of Definition
\scite{2.2}, not repeating colours except when absolutely necessary.

Lastly, let $p^+$ be such that Dom$(p^+) = d^*_{k(*)}$ and for each
$\beta \in d^*_{k(*)}$

$$
u^{p^+}(\beta) = \cup\{u^{p_{y,k(*)}}:y \in Y\};
$$
\mn
if $u^{p_{y,k(*)}(\beta)}$ is not defined, it means $\emptyset$

$$
c^{p^+}(\beta) \text{ extend each } c^{p_{y,k(*)}(\beta)} 
\text{ otherwise is } 1-to-1 \text{ with new values}.
$$
\mn
So $p^+ \ge p^*$ forces that $\{\alpha_\eta:\eta \in \text{
dom}_{\ell(*),m(*)}\}$ exemplify $\bold s = (\text{dom}_{\ell(*),m(*)},e) \in
ID(\underset\tilde {}\to f)$, a contradiction.  
${{}}$  \hfill$\square_{\scite{3.2}}$\margincite{3.2}
\enddemo
\bigskip

\proclaim{\stag{3.4} The example Claim}  Let $n(*),\ell(*),\lambda$ be
as in Definition \scite{2.2}.  Assume
\mr
\item "{$(a)$}"  $({\text{\rm dom\/}}_{\ell(*),m(*)},e^*) \in ID^*_{\ell(*)}$
\sn
\item "{$(b)$}"  $\Gamma \subseteq ID^*_{\ell(*)}$
\sn
\item "{$(c)$}"  if $\bold s \in \Gamma$ then
$\bold s$ is $(2^{\ell(*)}m(*))$-nice
\sn
\item "{$(d)$}"  $\Bbb P = \Bbb P^\lambda_\Gamma$
\sn
\item "{$(e)$}"  $\underset\tilde {}\to c$ is the $\Bbb P$-name
$\cup\{c^p:p \in {\underset\tilde {}\to G_{\Bbb P}}\}$
\sn
\item "{$(f)$}"  $\ell(*) \ge n(*)$.
\ermn
\ub{Then} $\Vdash_{\Bbb P} `` \underset\tilde {}\to c$ is a function
from $[\lambda]^2$ to $\mu$ exemplifying (dom$_{\ell(*),m(*)},e^*)$
does not belong to ID$_2(\lambda,\mu)"$.
\endproclaim
\bigskip

\remark{Remark}  The proof is similar to \cite{GcSh:491}.
\endremark
\bigskip
 
\demo{Proof}  So assume toward contradiction that $p \in \Bbb P$ and
$\alpha_\eta < \lambda$ for $\eta \in \text{ dom}_{\ell(*),m(*)}$
are such that $p$ forces that $\eta \mapsto \alpha_\eta$ is a counterexample,
i.e. $\langle \alpha_\eta:\eta \in \text{ dom}_{\ell(*),m(*)} \rangle$
is with no repetitions and $p$ forces that $t_1 e^* t_2 
\Rightarrow \underset\tilde {}\to c(\{\alpha_\eta:\eta
\in t_1\}) = \underset\tilde {}\to c(\{\alpha_\eta:\eta \in t_2\})$.
By \scite{2.4}(2) without loss of generality $\{\alpha_\eta:\eta \in \text{
dom}_{\ell(*),m(*)}\} \subseteq u^p$.

Let $Y = Y_{e^*} = \{y:y \in \text{ Dom}(e)$ and $y/e$ is not a
singleton$\}$ and for $\nu \in {}^{\ell(*)>} 2$ let $Y_\nu =
Y_{\nu,e^*} = \{\{\eta_0,\eta_1\} \in Y_{e^*}:\nu \char 94 <i>
\trianglelefteq \eta_i$ for $i=0,1\}$ as in the previous proof.
We now choose by induction on $\ell \le n(*)$ the objects
$\eta_\ell,\nu_\ell,Z_\ell$ and first order formulas
$\varphi_\ell(x,y_0,\dotsc,y_{\ell -1}),<^\ell_{y_0,\dotsc,y_{\ell
-1}}(x,y)$ in the vocabulary of $M^*_\lambda$ such that:
\mr
\item "{$\boxtimes(a)$}"  $\nu_\ell \in {}^\ell 2,
\eta_\ell \in \text{ dom}_{\ell(*),m(*)}$ and 
$M^*_\lambda \models (\exists^{\le \aleph_{(n(*)-\ell}}x)
\varphi_\ell(x,\alpha_{\eta_0},\dotsc,\alpha_{\eta_{\ell-1}})$
\sn
\item "{$(b)$}" $<^\ell_{\alpha_{\eta_0}},\dotsc,\alpha_{\eta_{\ell
-1}}$ is a well ordering of $\{x:M^*_\lambda \models
\varphi_\ell[x,\alpha_{\eta_0},\dotsc,\alpha_{\eta_{\ell -1}}]\}$ of
order type a cardinal $\le \aleph_{n(*)-\ell}$
\sn
\item "{$(c)$}"  $\nu_0 = <>,\varphi_0 = [x=x]$
\sn
\item "{$(d)$}"   $\nu_{\ell +1} = (\eta_\ell \restriction \ell) \char
94 \langle 1 - \eta_\ell(\ell) \rangle$ and $\nu_\ell \triangleleft \eta_\ell$
\sn
\item "{$(e)$}"  $Z_\ell = \{\eta:\nu_\ell \triangleleft \eta \in
\text{ dom}_{\ell(*),m(*)}$ and $\{\eta_s,\eta\} \in e_{\nu
\restriction s}$ for $s=0,1,\dotsc,\ell-1\}$
\sn
\item "{$(f)$}"  $\eta \in Z_\ell \Rightarrow \alpha_\eta \in
\{\beta:M^*_\lambda \models
\varphi_\ell[\beta,\alpha_{\eta_0},\dotsc,\alpha_{\eta_{\ell -1}}]\}$
\sn
\item "{$(g)$}" $\eta_\ell$ is such that:
{\roster
\itemitem{ $(\alpha)$ }  $\nu_\ell \triangleleft \eta_\ell \in Z_\ell$
\sn
\itemitem{ $(\beta)$ }  if $\nu_\ell \trianglelefteq \eta \in Z_\ell$
then $\alpha_\eta
\le^\ell_{\alpha_{\eta_0},\dotsc,\alpha_{\eta_{\ell-1}}},\alpha_{\eta_\ell}$.
\endroster}
\ermn
(See similar proof with more details in \scite{4.4}).
\sn
Let $\nu^* = \nu_{n(*)},Z = Z_{n(*)},Z^+ = \{\eta_\ell:\ell < n(*)\}
\cup Z$; note that by Definition \scite{1.3}(1), clause (b) and
Definition \scite{2.0} we have $|Z| \ge 2$, i.e., this is part of
$({\text{\rm dom\/}}_{\ell(*),m(*)},e^*)$ being 0-nice.  
For $\nu \in \{\nu_\ell:\ell < n(*)\}$ let $s_\nu$ be such
that:
$\rho_1 \cap \rho_2 = \nu \and \rho_1,\rho_2 \in \{\eta_\ell:\ell <
n(*)\} \cup Z \Rightarrow s_\nu = c\{\alpha_{\rho_1},\alpha_{\rho_2}\}$
(clearly exists).    By case 0 in Definition \scite{2.2}, \wilog

$$
u^p = \{\alpha_\eta:\eta \in Z^+\},
$$
\mn
that is, we may forget the other $\alpha \in u^p$; by 
claim \scite{3.1a} we have $p \in \Bbb P^\lambda_\emptyset$ so for some $k$
we have $r \in \Bbb P^{\lambda,\emptyset}_k$.

So we have
\mr
\item "{$\boxplus$}"   $\langle \eta_\ell:\ell < n(*)
\rangle,Z,Z^+,\langle \nu_\ell:\ell \le n(*) \rangle,\langle
s_{\eta_\ell \restriction \ell}:\ell < n(*) \rangle$ and $p$ 
are as above, that is
{\roster
\itemitem{ $(i)$ }  $\eta_\ell \in \text{ dom}_{\ell(*),m(*)},\nu_0 =
<>,\nu_{\ell +1} = (\eta_\ell \restriction \ell) \char 94 \langle
(1-\eta_\ell(\ell)) \rangle,Z = \{\rho \in \text{
dom}_{\ell(*),m(*)}:\nu_{n(*)} \triangleleft \rho \text{ and }
\{\eta_\ell,\rho\}/e$ is not a singleton for each $\ell < n(*)\}$
hence $|Z| \ge 2$ and $Z^+ = Z \cup \{\eta_\ell:\ell < n(*)\}$
\sn
\itemitem{ $(ii)$ }  $p \in \Bbb P^{\lambda,\emptyset}_k$
\sn
\itemitem{ $(iii)$ }   $\alpha_\eta \in u^p$ for $\eta \in Z^+$
\sn
\itemitem{ $(iv)$ }  $\langle \alpha_\eta:\eta \in Z^+ \rangle$ is with
no repetitions
\sn
\itemitem{ $(v)$ }  $c^p \restriction \{\alpha_\eta:\eta \in Z^+\}$
satisfies: \nl
if $\ell <  n(*)$ and $\nu \in Z \cup \{\eta_t:\ell < t < n(*)\}$
so $\eta_\ell \cap \nu = \eta_\ell \restriction \ell$ 
\ub{then} ($\alpha_\nu \ne \alpha_{\eta_\ell}$ and) 
$c(\{\alpha_\nu,\alpha_{\eta_\ell}\} = s_{\eta_\ell \restriction
\ell}$
\sn
\itemitem{ $(vi)$ }  $\{\alpha_\eta:\eta \in Z\} \subseteq c
\ell_0\{\alpha_{\eta_\ell}:\ell < n(*)\}$
\sn
\itemitem{ $(vii)$ }  $Z$ has at least two members.
\endroster}
\ermn
Among all such examples choose one with $k < \omega$ minimal.  The
proof now splits according to the cases in Definition \scite{2.2}.
\enddemo
\bn
\ub{Case $0$}:  $k=0$.

Trivial.
\bn
\ub{Case $1$}:

Let $p_1,\alpha$ be as there, so recall that $\{\alpha,\beta\} e^{p_1}
\{\alpha',\beta'\} \Rightarrow \{\alpha,\beta\} = \{\alpha',\beta'\}$.
Hence obviously, by clauses (v) and (vii) above, $\eta \in Z^+ \Rightarrow
\alpha_\eta \ne \alpha$, so $\{\alpha_\eta:\eta \in Z^+\} \subseteq
u^{p_1}$, contradicting the minimality of $k$.
\bn
\ub{Case $2$}:

Let $p_i = (u_i,c_i) \in \dbcu_{\ell < k} \Bbb
P^{\lambda,\emptyset}_\ell$ for $i=1,2$ be as there.  We now prove by
induction on $\ell < n(*)$ that $\alpha_{\eta_\ell} \in u_0 \cap u_1$.
If $\ell < n(*)$ and it is true for every $\ell' < \ell$, but
(for some $i \in \{1,2\}$), 
$\alpha_{\eta_\ell} \in u_i \backslash u_{3-i}$, it follows by clause
(v) of $\boxplus$ that the
sequence $\langle c(\{\eta_\ell,\nu\}):\nu \in Z^*_\ell \rangle$ is
constant where we let $Z^*_\ell = \{\eta_{\ell +1},\eta_{\ell +2},
\dotsc,\eta_{n(*)-1}\} \cup Z$, hence $\{\alpha_\nu:\nu \in Z^*_\ell\}$ is
disjoint to $u_{3-i} \backslash u_i$, so $\{\alpha_\nu:\nu \in Z^+\}
\subseteq u_i$, so we get contradiction to the minimality of $k$.

As $\{\alpha_{\eta_\ell}:\ell < n(*)\} \subseteq u_2 \cap u_1$
necessarily (by clause (vi) of $\boxplus$) we have
$\{\alpha_\nu:\nu \in Z^*_{n(*)}\} \subseteq c
\ell_0\{\alpha_{\eta_\ell}:\ell < n(*)\} \subseteq c \ell_0(u_2 \cap
u_1)$.  But $\{\alpha_\nu:\nu \in Z^*_{n(*)}\} \subseteq u_2 \cup u_1$
by $\boxplus(iii)$, and we know that $c \ell_0(u_2 \cap u_1) \cap
(u_2 \cap u_1) \subseteq u_1$ by clause (d) of Definition
\scite{2.2}, Case 2 hence
$\{\alpha_\nu:\nu \in Z^*_0\} \subseteq u_1$ contradiction to ``$k$
minimal".
\bn
\ub{Case $3$}:  

This case never occurs as 
$p \in \Bbb P^{\lambda,\emptyset}_k$.  \hfill$\square_{\scite{3.4}}$\margincite{3.4}
\bigskip

\proclaim{\stag{3.7} Theorem}  1) Let  $n(*)=4$ (or just $n(*) \ge 4$), 
$\lambda = \aleph_{n(*)},\ell(*) = n(*)+1$ 
and $2^{\aleph_\ell} = \aleph_{\ell +1}$
for $\ell < n(*)$.

For some c.c.c. forcing $\Bbb P$ of cardinality $\lambda$ in $\bold
V^{\Bbb P}$ the pair $(\lambda,\aleph_0)$ is not $\aleph_0$-compact. \nl
2) For given $\chi = \chi^{\aleph_0} \ge \lambda$ we can add $\bold
V^{\Bbb P} \models ``2^{\aleph_0} = \chi"$. 
\endproclaim
\bigskip

\demo{Proof}  1) Let $\Gamma_n = \{\bold s \in ID^*_{\ell(*)}:\bold s$ is
$n$-nice$\}$, see Definition \scite{2.0}, 
clearly $\Gamma_{n+1} \subseteq \Gamma_n$ and $\Gamma_n \ne \emptyset$
(see \scite{1.5}) for $n < \omega$ and
$\emptyset = \dbca_{n < \omega} \Gamma_n$ and let
$\Bbb P_n = \Bbb P^\lambda_{\Gamma_n}$
and let ${\underset\tilde {}\to c_n} = \cup\{c^p:p \in
{\underset\tilde {}\to G_{{\Bbb P}_m}}\}$, it is a $\Bbb P_n$-name and
$\Bbb P$ is $\dsize \prod_{n < \omega} \Bbb P_n$ with finite support.
Now the forcing notion 
$\Bbb P$ satisfies the c.c.c. as $\Bbb P_n$ satisfies the Knaster
condition (by \scite{2.4}(1)).
By \scite{3.4} we know that $\Vdash ``ID_2({\underset\tilde {}\to
c_n}) \cap ID^*_{\ell(*)} \subseteq \Gamma_n"$ for $\Bbb P_n$ hence
for $\Bbb P$, in fact it is not hard
to check that equality holds.  If $\aleph_0$-compactness holds then in
$\bold V^{\Bbb P}$ for some $c:[\lambda]^2 \rightarrow \omega$ we have
$ID_2(c) \cap ID^*_{\ell(*)} \subseteq \dbca_n \Gamma_n = \emptyset$
by claim \scite{1.4}.
\nl
But $\bold V^{\Bbb P}$, if $c:[\lambda]^2 \rightarrow \omega$ then by
\scite{3.2}(2) it realizes some 
$\bold s \in \cup\{\Gamma_n:n < \omega\} \subseteq ID^*_{\ell(*)}$ 
(even $k$-nice one for every $k < \omega$).

Together we get that the pair 
$(\lambda,\aleph_0)$ is not $\aleph_0$-compact. \nl
2) We let $\Bbb Q$ be adding $\chi$ Cohen reals, i.e. $\{h:h$ a finite
function from $\chi$ to $\{0,1\}\}$ ordered by inclusion.  Let $\Bbb
P$ be as above and force with $\Bbb P^+ = \Bbb P \times \Bbb Q$, now
it is easy to check that $\Bbb P^+$ is as required.
\hfill$\square_{\scite{3.7}}$\margincite{3.7}
\enddemo
\newpage

\head {\S4 Improvements and additions} \endhead  \resetall \sectno=4
\bigskip

Though our original intention was to deal with the possible
incompactness of the pair $(\aleph_2,\aleph_0)$, we have so far dealt
with $(\lambda,\aleph_0)$ where $2^{\aleph_0}
\ge \lambda = \aleph_{n(*)} \and n(*) \ge 4$.  
For dealing with $(\aleph_3,\aleph_0),(\aleph_2,\aleph_0)$,
that is $n(*) = 3,2$ we need to choose $M^*_\lambda$ more carefully.

What is the problem in \S3 concerning $n(*) = 2$? \nl
On the one hand in the proof of \scite{3.4} we need that there are
many dependencies among ordinals $< \lambda$ by $M^*_\lambda$; so if
$\lambda$ is smaller this is easier, but really just make us use
larger $\ell(*)$ help.

On the other hand, in the proof of \scite{3.2} we use \scite{3.1}, a
partition theorem, so here if $\lambda$ is bigger it is easier; but
instead we can use demands specifically on $M^*_\lambda$.  Along those
lines we may succeed for $n(*)=3$ using \scite{3.1}(1) rather than
\scite{3.1}(1) but we still have problems for the pair
$(\aleph_2,\aleph_0)$; here we change the main definition \scite{2.2},
in case 3 changes $\langle v_y:y \in Y^+ \rangle$, i.e. for $\eta \in
\text{ dom}_{\bold s}$ we have $v^+_{\{\eta\}},v^-_{\{\eta\}}$ instead
$v_{\{\eta\}}$.  For this we have to carefully reconsider \scite{3.2}, but
the parallel of \scite{3.1} is easier.  Note that in \S2, \S3 we could
have used a nontransitive version of $c \ell_M(-)$.
\bigskip

\definition{\stag{4.2} Definition}  We say that $M^*$ is $(\lambda,<
\mu,n(*),\ell(*))$-suitable \ub{if}:
\mr
\item "{$(a)$}"  $M^*$ is a model of cardinality $\lambda$
\sn
\item "{$(b)$}"  $\lambda$ is $> \mu,\le \mu^{+ n(*)}$ and $n(*) < 
\ell(*) < \omega$
\sn
\item "{$(c)$}"  $\tau_{M^*}$, the vocabulary of $M^*$, is of
cardinality $\le \mu$
\sn
\item "{$(d)$}"  for every subset $A$ of $M^*$ of cardinality $< \mu$,
\nl
the set $c \ell_{M^*}(A)$ has cardinality $< \mu$.
\sn
\item "{$(e)$}"  for some $m^* < \omega$ we have: \nl
if $\bold s = (\text{dom}_{\ell(*),m},e) \in ID^*_{\ell(*)}$ and $a_\eta
\in M^*$ for $\eta \in \text{ dom}_{\ell(*),m}$ and $m > m^*$, \ub{then}
we can find $\langle \eta_\ell:\ell < n(*) \rangle$ and $\langle
\nu_\ell:\ell \le n(*) \rangle$ such that
{\roster
\itemitem{ $(\alpha)$ }  $\eta_\ell \in \text{ dom}_{\ell(*),m}$
\sn
\itemitem{ $(\beta)$ }  $\nu_0 = <>,\nu_{\ell +1} = (\eta_\ell
\restriction \ell) \char 94 \langle 1 - \eta_\ell(\ell) \rangle$
\sn
\itemitem{ $(\gamma)$ }  $\nu_\ell \triangleleft \eta_\ell$
\sn
\itemitem{ $(\delta)$ }   $Z = \{\rho \in \text{
dom}_{\ell(*),m}:\nu_{n(*)} \triangleleft \rho$ and in the graph
$H[e],\rho$ is connected to $\eta_\ell$ for $\ell =
0,\dotsc,n(*)-1\}$
\sn
\itemitem{ $(\varepsilon)$ }  $\{\alpha_\rho:\rho \in Z\} \subseteq c
\ell_{M^*}\{\alpha_{\eta_\ell}:\ell < n(*)\}$.
\endroster}
\endroster
\enddefinition
\bigskip

\definition{\stag{4.3} Definition}  1) We say that $M$ is explicitly$^1$
$(\lambda, < \mu,n(*))$-suitable if:
\mr
\item "{$(a)$}"  $M^*$ is a model of cardinality $\lambda$
\sn
\item "{$(b)$}"  $\lambda = \mu^{+ n(*)}$
\sn
\item "{$(c)$}"  $\tau_{M^*}$, the vocabulary of $M^*$, is of
cardinality $\le \mu$
\sn
\item "{$(d)$}"  for $A \subseteq M^*$ of cardinality $< \mu$, the set
$c \ell_{M^*}(A)$ has cardinality $< \mu$ and $A \ne \emptyset \wedge
\mu > \aleph_0 \Rightarrow \omega \subseteq c \ell_{M^*}(A)$
\sn
\item "{$(e)$}"  for some $\langle R_\ell:\ell \le n(*) \rangle$ we have
{\roster
\itemitem{ $(\alpha)$ }  $R_\ell$ is an $(\ell + 2)$-place predicate
in $\tau_{M^*}$; we may write $R_\ell(x,y,z_0,\dotsc,z_{\ell -1})$ as $x
<_{z_0,\dotsc,z_{\ell -1}} y$ or $x <_{\langle z_0,\dotsc,z_{\ell -1} 
\rangle} y$
\sn
\itemitem{ $(\beta)$ }  for any $c_0,\dotsc,c_{\ell -1} \in M^*$, the
two place relation $<_{c_0,\dotsc,c_{\ell -1}}$
(i.e. $\{(a,b):\langle a,b,c_0,\dotsc,c_{\ell -1} \rangle \in R^{M^*}\})$ is a
well ordering of $A_{c_0,\dotsc,c_{\ell -1}} =: A_{\langle c_0,\dotsc,c_{\ell
-1} \rangle} =: \{b:(\exists x)(x <_{c_0,\dotsc,c_{\ell -1}} b \vee b
<_{c_0,\dotsc,c_{\ell -1}}x)\}$ of order-type a cardinal
\sn
\itemitem{ $(\gamma)$ }  $R^{M^*}_0$ is a well ordering of $M^*$ of
order type $\lambda$
\sn
\itemitem{ $(\delta)$ }  if $\bar c = \langle c_\ell:\ell < k \rangle$
and $<_{\bar c}$ is a well ordering of $A_{\bar c}$ of order type
$\mu^{+m}$ then for every $c_k \in M^*_{\bar c}$ we have $A_{\bar c
\char 94 \langle c_k \rangle} = \{a \in A_{\bar c}:a <_{\bar c} c_k\}$
so is empty if $c_k \notin A_{\bar c}$
\sn
\itemitem{ $(\varepsilon)$ }  if $\bar c = \langle c_\ell:\ell < k
\rangle \in {}^k(M^*)$ and $|A_{\bar c}| < \mu$ then $A_{\bar c} \subseteq c
\ell_{M^*}(\bar c)$.
\endroster}
\ermn
2) We say that $M$ is explicitly$^2$ \, $(\lambda,<
\mu,n(*))$-suitable if
\mr
\item "{$(a)-(d)$}"  as in part (1)
\sn
\item "{$(e)$}"  for some $\langle R_\ell:\ell \le n(*) \rangle$ we
have (like (e) but we each time add $z$'s and see clause $(\delta)$)
{\roster
\itemitem{ $(\alpha)$ }  $R_\ell$ is a $(2\ell+2)$-place predicate in
$\tau_\mu$; we may write $R_\ell(x,y,z_0,\dotsc,z_{z \ell-1})$ or $x
<_{z_0,\dotsc,z_{2 \ell-1}}y$ or $x <_{\langle z_0,\dotsc,z_{2 \ell-1}
\rangle}y$
\sn
\itemitem{ $(\beta)$ }  for any $c_0,\dotsc,c_{2 \ell-1} \in M^*$ the
two-place relation $<_{c_0,\dotsc,c_{2 \ell-1}}$ (i.e.,
$\{(a,b):\langle a,b,c_0,\dotsc,c_{2 \ell-1} \rangle \in
R^{M^*}_\ell\}$) is a well
ordering of $A_{c_0,\dotsc,c_{2 \ell-1}} = A_{(c_0,\dotsc,c_{2
\ell-1})} = \{b:$ for some $a,\langle a,b,c_0,\dotsc,c_{2 \ell-1}
\rangle \in R^{M^*}_\ell$ or $\langle b,a,c_0,\dotsc,c_{2 \ell-1}
\rangle \in R^{M^*}_\ell$
\sn
\itemitem{ $(\gamma)$ }  $R^{M^*}_0$ is a well ordering of $M^*$ of
order type $\lambda$; for simplicity $R^{M^*}_0 = c \restriction \lambda$
\sn 
\itemitem{ $(\delta)$ }  if $\bar c = \langle c_\ell:\ell < 2k
\rangle$ and $<_{\bar c}$ is a well ordering of $A_{\bar c}$ of order type
$\mu^{+m}$ \ub{then}
for any $c_{2k},c_{2k+1} \in M^*$ we have $A_{\bar c \char 94 \langle
c_{2k},c_{2k+1}\rangle}$ is empty if $\{c_{2k},c_{2k+1}\} \nsubseteq
A_{\bar c}$ and otherwise is $\{a \in A_{\bar c}:a <_{\bar c} c_{2k}$
and $a < c_{2k+1}\}$.
\endroster}
\endroster
\enddefinition
\bigskip

\demo{\stag{4.4} Observation}  1) If $M$ is an explicitly$^1$
$(\lambda,< \mu,n(*))$-suitable model, \ub{then} $M$ is a $(\lambda,<
\mu,n(*)+1,\ell(*))$-suitable model if $\ell(*) > n(*)+1$. \nl
2) If $M$ is an explicitly$^2 \, (\lambda,< \mu,n(*))$-suitable
model, \ub{then} $M$ is a $(\lambda,< \mu,2n(*)+2,2n(*)+3)$-suitable model.
\enddemo
\bigskip

\demo{Proof}  1) Straightforward, similar to inside the proof of
\scite{3.4} and as we shall use part (2) only and the proof of (1) is
similar but simpler, we do not elaborate. \nl
2) Clearly clauses (a) - (d) of Definition \scite{4.2} holds, so we
deal with clause (e). 
So assume $\ell(*) \ge 2n(*)$ and $\bold s =
(\text{dom}_{\ell(*),m},e) \in ID^*_{\ell(*)}$ and $\alpha_\eta \in M$
for $\eta \in \text{ dom}_{\ell(*),m}$ are pairwise distinct.  
We choose by induction on
$\ell \le n(*)$ the objects $\eta_{2 \ell},\nu_{2 \ell +1},Z_{2
\ell},\eta_{2 \ell +1},\nu_{2 \ell +2},Z_{2 \ell +1}$ such that node
$\nu_\ell = <>$ and $\nu_{2 \ell +2}$ is chosen in stage $\ell$
\mr
\item "{$\boxtimes(a)$}"  $\nu_\ell \in {}^\ell 2,
\eta_\ell \in \text{ dom}_{\ell(*),m(*)}$ and 
$M \models (\exists^{\le \aleph_{(n(*)-\ell}}x)
\varphi_\ell(x,\alpha_{\eta_0},\dotsc,\alpha_{\eta_{2 \ell-1}})$
\sn
\item "{$(b)$}" $<^\ell_{\alpha_{\eta_0},\dotsc,\alpha_{\eta_{2 \ell
-1}}}$ is a well ordering of $A_{\langle
\alpha_{\eta_0},\dotsc,\alpha_{2 \ell-1} \rangle} =: \{x:M \models
\varphi_\ell[x,\alpha_{\eta_0},\dotsc,\alpha_{\eta_{\ell -1}}]\}$ of
order type a cardinal $\le \aleph_{n(*)-\ell}$
\sn
\item "{$(c)$}"  $\nu_0 = <>,\varphi_0 = [x=x]$
\sn
\item "{$(d)$}"   $\nu_{\ell +1} = (\eta_\ell \restriction \ell) \char
94 \langle 1 - \eta_\ell(\ell) \rangle$
\sn
\item "{$(e)$}"  $Z_\ell = \{\eta:\nu_{2 \ell} \triangleleft \eta \in
\text{ dom}_{\ell(*),m(*)}$ and $\{\eta_s,\eta\} \in e_{\nu
\restriction s}$ for $s=0,1,\dotsc,\ell-1\}$
\sn
\item "{$(f)'$}"  $\eta \in Z_\ell \Rightarrow \alpha_\eta \in
A_{\langle a_{\eta_k}:k < 2 \ell \rangle}$
\sn
\item "{$(g)$}" $\eta_\ell$ is such that:
{\roster
\itemitem{ $(\alpha)$ }  $\nu_\ell \triangleleft \eta_\ell \in Z_\ell$
\sn
\itemitem{ $(\beta)$ }  if $\nu_\ell \trianglelefteq \eta \in Z_\ell$
then [$\ell$ even $\Rightarrow \alpha_\eta 
\le_{\alpha_{\eta_0},\dotsc,\alpha_{\eta_{\ell-1}}},\alpha_{\eta_\ell}$]
and [$\ell$ odd $\Rightarrow \alpha_\eta \le \alpha_{\eta_\ell}$].
\endroster}
\ermn
How do we do the induction step?  Arriving to $\ell$ we have already
defined $\langle \nu_k:k \le 2 \ell \rangle,\langle \eta_k:k < 2 \ell
\rangle$ and $\langle Z_k:k < 2 \ell \rangle$, recalling $\nu_0 =<>$.
So by the definition of $Z_k$ also $Z_{2 \ell}$ is well defined and
$\{\alpha_\eta:\eta \in Z_{2 \ell}\}$ is included in 
$A_{\langle a_{\eta_k}:k < 2 \ell \rangle}$ and let
$\eta_{2 \ell} \in Z_{2 \ell}$ be such that $\eta \in Z_{2 \ell}
\Rightarrow a_\eta \le_{\langle a_{\eta_k}:k < 2 \ell \rangle} a_{\eta_{2
\ell}}$ and $\nu_{2 \ell +1} = \nu_{2 \ell} \char 94 \langle 1 - \ell,\eta_{2
\ell}(2 \ell) \rangle = (\eta_{2 \ell} \restriction (2 \ell)) \char 94
\langle 1 - \eta_{2 \ell}(2 \ell) \rangle$ so $Z_{2 \ell +1}$ is well
defined.  Let $\eta_{2 \ell +1} \in Z_{2 \ell +1}$ be such that $\eta
\in Z_{2 \ell} \Rightarrow a_\eta \le \alpha_{\eta_{2 \ell +1}}$ and
$\nu_{2 \ell +2} = \nu_{2 \ell +1} \char 94 \langle 1 - \eta_{2 \ell +1}(2
\ell +1) \rangle$ and we have carried the induction.
\hfill$\square_{\scite{4.4}}$\margincite{4.4}
\enddemo
\bn
Are there such models?  We shall use \scite{4.5}(2), the others are
for completeness (i.e. part (3) is needed for $\lambda = \aleph_3$ and
part (4) says concerning $\lambda = \aleph_2$ it suffices to use $ID^*_3$):
\demo{\stag{4.5} Observation}  1) For $\mu$ regular uncountable, there is
an explicitly$^1$ $(\mu^{+2},< \mu,2)$-suitable model. \nl
2) If $\mu = \aleph_0$, then there is an explictly$^2$ $(\mu^{+2},<
\mu,2)$-suitable model. \nl
3) If $\mu$ is regular uncountable, $t=1$ or $\mu = \aleph_0 \and t=2$ 
and $n \in [3,\omega)$, \ub{then} there is an explicitly$^t$ 
$(\mu^{+n},< \mu,n)$-suitable model. \nl
4) If $2^{\aleph_0} = \aleph_1,\mu = \aleph_0$ \ub{then} for some
$\aleph_2$-c.c., $\aleph_1$-complete forcing notion $\Bbb Q$ of
cardinality $\aleph_2$ in $\bold V^{\Bbb Q}$ there 
is an explicitly $(\aleph_2,< \aleph_0,2)$-suitable model.
\enddemo
\bigskip

\remark{\stag{4.5a} Remark}  1) It should be clear that if $\bold V = \bold
L$ (or just $\neg 0^\#$), then this works also for singular $\mu$
\ub{but} more reasonable is to use nontransitive closure. 
\endremark
\bigskip

\demo{Proof}  1), 2)  Let $t=1$ for part (1) and $t=2$ for part (2).  Let
$n(*)=2$ and $\lambda = \mu^{+2}$.  We choose $M_\alpha$ by
induction on $\alpha \le \lambda$ such that:
\mr
\widestnumber\item{$(\varepsilon)(ii)$}
\item "{$(\alpha)$}"  $M_\alpha$ is a $\tau^-$-model where $\tau^- =
\{R_0,R_1,R_2\}$ with $R_\ell$ is $(t,\ell +2)$-predicate and $x
<_{\bar z} y$ means $R_\ell(x,y,\bar z)$
\sn
\item "{$(\beta)$}"  $M_\alpha$ is increasing with $\alpha$ and has
universe $1 + \alpha$
\sn
\item "{$(\gamma)$}"  $R^{M_\alpha}_0$ is $< \restriction \alpha$ (and
$A^{M_\alpha}_{<>} = \alpha$)
\sn
\item "{$(\delta)$}"  for $\bar c \in {}^{tk}(M_\alpha),k = 0,1,2$ we
have $<_{\bar c}$ is a well ordering of $A^{M_\alpha}_{\bar c} =:
\{a:M_\alpha \models (\exists x)(a <_{\bar c} x \vee x <_{\bar c} a)\}$ of
order type a cardinal $< \mu^{+(n(*) + 1-k)}$
\sn
\item "{$\quad(\varepsilon)(i)$}"  if 
$t=1,\bar c \in {}^k(M_\alpha),k=0,1$ and $d
\in A^{M_\alpha}_{\bar c}$ then $A^{M_\alpha}_{\bar c \char 94 <d>} =
\{a \in A^{M_\alpha}_{\bar c}:M_\alpha \models a <_{\bar c} d\}$
\sn
\item "{${{}}(ii)$}"  if 
$t=2,\bar c \in {}^{2k}(M_\alpha),k=0,1$ and $d_0,d_1
\in A^{M_\alpha}_{\bar c}$ then $A^{M_\alpha}_{\bar c \char 94 \langle
d_0,d_1 \rangle} =
\{a \in A^{M_\alpha}_{\bar c}:M_\alpha \models ``a <_{\bar c} d_0 \and
a < d_1"\}$
\sn
\item "{$(\zeta)$}"  if $A$ is a subset of $M_\alpha$ of cardinality
$< \mu$ then $c \ell^*_{M_\alpha}(A)$ is of cardinality $< \mu$ and $c
\ell^*_{M_\alpha}(c \ell^*_{M_\alpha}(A)) = c \ell^*_{M_\alpha}(A)$ where
{\roster
\itemitem{ $\boxtimes$ }  for $A \subseteq M_{\alpha'},c
\ell^*_{M_\alpha}(A)$ is the minimal set $B$ 
such that: $A \subseteq B$ and $(\forall \bar c \in
{}^{2t} B)(|A^{M_\alpha}_{\bar c}| < \mu \rightarrow A_{\bar c} \subseteq B)$;
clearly $B$ exists and $c \ell^*_{M_\alpha}(\emptyset) = \emptyset$
\endroster}
\sn
\item "{$(\eta)$}"   for every $\beta < \alpha,k = 1,2$ and $\bar
c \in {}^k(M_\beta)$ we have $A^{M_\alpha}_{\bar c} =
A^{M_\beta}_{\bar c}$
\sn
\item "{$(\theta)$}" if $A \subseteq \beta < \alpha$ then $c
\ell^*_{M_\beta}(A) = c \ell^*_{M_\alpha}(A)$ 
\sn
\item "{$(\iota)$}"  if $t =2$ and $\mu = \aleph_0$ and $A \subseteq
\alpha$ is finite, $\beta$ is the last element in $A$, \ub{then} for some
finite $B \subseteq \beta$ we have $c \ell^*_{M_\alpha}(A) = \{\beta\}
\cup c \ell^*_{M_\beta}(B)$.
\ermn
We leave the cases $\alpha < \mu$ and $\alpha$ a limit ordinal to the
reader (for $(\zeta)$ we use $(\theta)$) and assume $\alpha = \beta +1$
and $M_\gamma$ for $\gamma \le \beta$ are defined.  We can choose
$\langle B_{\beta,i}:i < \mu^+ \rangle$, a (not necessarily strictly)
increasing sequence of subsets of 
$\beta$, each of cardinality $\le \mu,B_{\beta,0} =
\emptyset$ and $\cup \{B_{\beta,i}:i < \mu^+\} = \beta$ and $c
\ell^*_{M_\beta}(B_{\beta,i}) = B_{\alpha,i}$.

For each $i < \mu^+$ let $\langle B_{\beta,i,\varepsilon}:\varepsilon
< \mu \rangle$ be (not necessarily strictly) 
increasing sequence of subsets of $B_{\beta,i}$
with union $B_{\beta,i}$ such that $c
\ell^*_{M_\beta}(B_{\beta,i,\varepsilon}) =
B_{\beta,i,\varepsilon},B_{\beta,0} = \emptyset$.  Let $<^*_\beta$ be
a well ordering of $\{\gamma:\gamma < \beta\}$ such that each
$B_{\beta,i}$ is an initial segment so it has order type $\mu^+$. 
For $\gamma \in B_{\beta,i+1}
\backslash B_{\beta,i}$ let $<^*_{\beta,\gamma}$ be a well ordering of
$A^*_{(\beta,\gamma)} = \{\xi:\xi <^*_\beta \gamma\}$ of order type
$\le \mu$ such that $(\forall
\varepsilon < \mu)(B_{\beta,i+1,\varepsilon} \cap A^*_{(\beta,\gamma)}$
is an initial segment of $A^*_{(\beta,\gamma)}$ by $<^*_{\beta,\gamma}$.

Now we define $M_\alpha$:

universe is $\alpha$

$R^{M_\alpha}_0 = < \restriction \alpha$
\mn
\ub{Case 1}:  $t=1$.

$R^{M_\alpha}_1 = R^{M_\beta}_1 \cup \{(a,b,\beta):a <^*_\beta b\}$

$R^{M_\alpha}_2 = R^{M_\beta}_2 \cup \{(a,b,\beta,\gamma):\gamma <
\beta$,  and $a <^*_{\beta,\gamma} b$ hence $a <^*_\beta \gamma \and b
<^*_\beta \gamma$ and $a,b \in B_{\beta,i+1}$ for the unique $i$ such
that $\gamma \in B_{\beta,i+1} \backslash B_{\beta,i}\}$.
\mn
\ub{Case 2}:  $t=2$.
\block
$R^{M_\alpha}_1 = R^{M_\beta}_1 \cup \{(a,b,\beta,\gamma):a <^*_\beta
b$ and $a < \gamma,b < \gamma$ and, of course, \nl

\hskip75pt $a,b,\beta \in \alpha\}$.
\endblock
\sn
\block
$R^{M_\alpha}_2 = R^{M_\beta}_2 \cup
\{(a,b,\beta,\gamma_0,\beta_1,\gamma_1):a,b,\gamma_0,\beta_1 \in
\alpha$ and $a < \beta$, \nl

\hskip120pt $b < \beta,a < \gamma_0,b <
\gamma_0,a,b,\beta_1,\gamma_1 \in A_{<\beta,\gamma_0>}$ \nl

\hskip120pt and $a <^*_{\beta,\gamma_0} b$ and $a < \gamma_1,b < \gamma_1\}$.
\endblock
\sn
To check for clause $(\zeta)$ is easy if $\mu = \text{ cf}(\mu) >
\aleph_0$ and follows by clause $(\iota)$ if $\mu = \aleph_0$. \nl
Having carried the induction we define $M$: it is $M_\lambda$ expanded
by $\langle F^M_i:i < \mu \rangle$ such that: if $\bar c \in {}^{2t}
\lambda = {}^{2t}(M_\lambda)$ and $A_{\bar c}$ 
is a non empty well defined and of cardinality  $< \mu$ (which follows)
then $\{F^M_i(\bar c):i < \mu\}$ list $A_{<c_0,c_1,c_2>} \cup \{0\}$
otherwise $\{F^M_i(\bar c):i < \mu\}$ is $\{0\}$. \nl
3) Similar and used only for $(\aleph_3,\aleph_0)$ so we do not elaborate. \nl
4) Let $\Bbb Q$ be defined as follows: \nl
$p \in \Bbb Q$ \ub{iff}
\mr
\item "{$(\alpha)$}"  $p$ is a $\tau^-$-model, as in $(\alpha)$ of the
proof of part (1)
\sn
\item "{$(\beta)$}"  the universe univ$(p)$ of $p$ is a countable
subset of $\lambda$, we let $A^p_{<>} = \text{ univ}(p)$
\sn
\item "{$(\gamma)$}"  $R^p_0 = < \restriction \text{ univ}(p)$ and
$<_{<>} = R^p_0$
\sn
\item "{$(\delta)$}"  if $\bar c \in {}^{tk}(\text{univ}(p)),k=1,2$
then $<_{\bar c} = <^p_{\bar c}$ is a well ordering of $A^p_{\bar c} = \{a \in
p:p \models (\exists x)(a <_{\bar c} x \vee x <_{\bar c} a)$ and for $d \in A^p_{\bar c}$ we
let $A^p_{\bar c \char 94 <d>} = \{a \in A^p_{\bar c}:a <^p_{\bar c} d\}$
\sn
\item "{$(\varepsilon)$}"  $(A^p_{\bar c},<_{\bar c})$ has order type $\omega$ if $k =2$
\sn
\item "{$(\zeta)$}"   if $A \subseteq \text{ univ}(p)$ is finite,
then $c \ell^*_p(A)$ is finite (is defined as in(2)).
\ermn
\ub{the order}:
\sn
$\Bbb Q \models p \le q$ \ub{iff}
\mr
\widestnumber\item{$(iii)$}
\item "{$(i)$}"  $p$ is a submodel of $q$
\sn
\item "{$(ii)$}"  if $\bar c \in {}^2(\text{univ}(p))$ then $A^p_{\bar
c} = A^q_{\bar c}$
\sn
\item "{$(iii)$}"  if $\bar c \in {}^1(\text{univ}(p))$ then
$A^p_{\bar c}$ is an initial segment of $A^q_{\bar c}$ by $<_{\bar
c}$.
\ermn
The rest should be clear.  \hfill$\square_{\scite{4.5}}$\margincite{4.5}
\enddemo
\bigskip

\proclaim{\stag{4.6} Claim}  Assume (main case is $n(*)=2$)
\mr
\item "{$(*)$}"  $2 \le n(*) < \omega,\lambda = \aleph_{n(*)},\ell(*) =
2n(*)+3$ and $\lambda \le \chi = \chi^{\aleph_0}$. 
\ermn
\ub{Then} for some $\Bbb P^*$ we have
\mr
\item "{$(a)$}"   $\Bbb P^*$ is a forcing notion of cardinality
$\chi$
\sn
\item "{$(b)$}"   $\Bbb P^*$ satisfies the c.c.c.
\sn
\item "{$(c)$}"    in $\bold V^{{\Bbb P}^*}$ the pair
$(\aleph_{n(*)},\aleph_0)$ is not compact 
\sn
\item "{$(d)$}"  in $\bold V^{{\Bbb P}^*}$ we have $2^{\aleph_0} = \chi$.
\endroster
\endproclaim
\bigskip

\remark{Remark}  We intend to prepare a full version.
\endremark
\bigskip

\demo{Proof} We repeat \S2, \S3 with the following changes.

If $n(*) \ge 3$, we need change (A) below and using \scite{3.1}(2)
instead of \scite{3.1}(1).  For $n(*)=2$ we need all the changes below
\mr
\item "{$(A)$}"  we replace $M^*_\lambda$ by any model as in
\scite{4.5}(2) if $n(*)=2$, \scite{4.5}(3) if $n(*) \ge 3$
\sn
\item "{$(B)$}"  in Definition \scite{2.2}, case 3: we add $\langle
v^+_{\{\eta\}},v^-_{\{\eta\}}:\eta \in \text{ dom}_{\ell(*),m}
\rangle,v^+_\emptyset,v^-_\emptyset$
{\roster
\itemitem{ $(d)'(i)$ }  $v_y \supseteq u^{p_y}$ for $y \in Y$
\sn
\itemitem{ ${{}}(ii)$ }    if $\eta_1 <_{\text{lex}} \eta_2
<_{\text{lex}} \eta_3$ are from dom$_{\ell(*),m}$ and
$\{\eta_1,\eta_2\},\{\eta_1,\eta_3\} \in Y$, \ub{then} \nl

$v_{\{\eta_1,\eta_2\}} \cap v_{\{\eta_1,\eta_3\}} = v^+_\eta$
\sn
\itemitem{ ${{}}(iii)$ }  if $\eta_1 <_{\text{lex}} \eta_2
<_{\text{lex}} \eta_3$ are from dom$_{\ell(*),m}$ and
$\{\eta_1,\eta_3\},\{\eta_2,\eta_3\} \in Y$, \ub{then} \nl

\hskip10pt $v_{\{\eta_1,\eta_3\}} \cap v_{\{\eta_2,\eta_3\}} = v^-_\eta$ \nl

\hskip10pt $p_{\{\eta_1,\eta_3\}} \upharpoonleft v^-_\eta =
 p_{\{\eta_2,\eta_3\}} \upharpoonleft v^-_\eta$
\sn
\itemitem{ ${{}}(iv)$ }  $v_{\{\eta\}} = v^+_\eta \cup v^-_\eta$  
\sn
\itemitem{ ${{}}(v)$ }  if $\eta_1 \ne \eta_2$ then $v^+_{\eta_1} \cap
v^+_{\eta_2} = v^+_\emptyset$ and $v^-_{\eta_2} \cap v^-_{\eta_2} =
v^-_\emptyset$
\sn 
\itemitem{ ${{}}(vi)$ }  if $\eta_1 <_{\text{lex}} \eta_2
<_{\text{lex}} \eta_2$ are from dom$_{\ell(*),m}$, \ub{then} $p
\restriction v_\eta \in \dbcu_{r<k} \Bbb P_r$
\sn
\itemitem{ $(e)\,\,\,\,\,\,\,$ }    if $\eta_1 <_{\text{lex}} \eta_2$ 
are from dom$_{\ell(*),m}$ and $t = \{\eta_1,\eta_2\} \in Y$ then $c
\ell(v^+_{\eta_1}) \cap v_{\{\eta_1,\eta_2\}} = v^+_{\eta_1},c
\ell(v^-_{\eta_2}) \cap v_{\{\eta_1,\eta_2\}} = v^-_{\eta_2}$ if
$\{\eta_1,\eta_2\} \in Y,\eta \in \text{ dom}_{\ell(*),m} \backslash
\{\eat_1,\eta_1\}$ then $c \ell(v^\pm_\eta) \cap v_{\{\eta_1,\eta_2\}}
\subseteq v^\pm_{\eta_1}$
\sn
\itemitem{ $(f)\,\,\,\,\,\,\,$ } the functions $\langle c^{p_\eta}:\eta \in y
\rangle$ are pairwise compatible
\endroster}
\item "{$(C)$}"  in \scite{3.1}
{\roster
\itemitem{ $(a)$ }  $\lambda \ge (2^\mu)^+,\mu = \mu^{< \mu},(\forall
A \in [M]^{< \mu})(|c \ell_M(A)| < \mu)$
\sn
\itemitem{ $(b)$ }  the conclusion: change as in Definition
\scite{2.2}, case 3
\sn
\itemitem{ $(c)$ }  proof: 
\endroster}
\ermn
\ub{Case 1}:  $\mu = \aleph_0$:  let $g:[\lambda]^2 \rightarrow
\omega$ be $g(t) = |c \ell_M(\alpha,t \cup w_t)| < \omega$.
\nl
Let $W_1 \in [\lambda]^{\mu^+}$ be such that $g \restriction [W]^2$ is
constant say $k(*)$ and $f \restriction [W]^2$ is constantly $\gamma$.
Let $c \ell_M(t) = \{\zeta_{t,\ell}:\ell <
g(t)\}$.  By Ramsey theorem, there is an infinite $W \subseteq W_1$
such that:
\mr
\item "{$\circledast$}"  the truth value on
$\zeta_{\{\alpha_1,\beta_2\},\ell_1} =
\zeta_{\{\alpha_2,\beta_2\},\ell_2}$ depend just on 
$\ell_1,\ell_2$, T.V.$(\alpha_i,\beta_j)$, T.V.$(\beta_j < \alpha_i)$
for $i,j \in \{1,2\}$. 
\ermn
The conclusion should be clear.
\mr
\item "{$(D)$}"  $p$ in the proof of \scite{3.1a}: only case 3B need
case, assuming $m(*) > 2|u^{p^*}|$, the relevant subgraph has no cycle
by clause (e) of case 3 of \scite{2.2} we are done
\sn
\item "{$(E)$}"  in the proof of \scite{3.2}, we will have
$q^+_{\eta,\ell},q^-_{\eta,\ell}$ with domain 
$\subseteq v^-_\eta,v^+_\eta$ respectively and $q^+_\ell,q^-_\ell$
such that if $\eta_1 <_{\text{lex}} \eta_2$ and $\{\eta_1,\eta_2\} \in Y$
then $p_{\{\eta_1,\eta_2\},\ell} \restriction v^+_{\eta_1} =
q^+_{\eta_1,\ell}$ and
$p_{\{\eta_1,\eta_2\},\ell} \restriction v^-_{\eta_2} =
q^-_{\eta_2,\ell}$ and $q^+_{\eta,\ell},q^-_{\eta,\ell}$ are
compatible, $q^+_{\eta,\ell} \restriction v^+_\emptyset =
q^+_\ell,q^-_{\eta,\ell} \restriction v_\emptyset = q^-_\ell$
\sn
\item "{$(F)$}"  \scite{3.4}: part of the work has already been done
in \scite{4.1} - \scite{4.4}.
\endroster
\enddemo
\newpage

\head {\S5 Open Problems and concluding remarks} \endhead  \resetall \sectno=5
\bigskip

We finish the paper by listing some problems (some are old, see \cite{CK}).
\bn
\margintag{5.1}\ub{\stag{5.1} Question}:     Suppose that $\lambda$ is 
strongly inaccessible, $\mu > \aleph_0$ is regular not Mahlo 
and $\square_\mu$. \ub{Then}
$\lambda \rightarrow \mu$ in the $\lambda$-like model sense, i.e.
if a first order $\psi$ has a $\lambda$-like model then it has a
$\mu$-like model. \nl
If $\lambda$ is $\omega$-Mahlo, the answer is yes, see
\cite{ScSh:20} by appropriate partition theorems.  
The assumption that $\mu$ is Mahlo is necessary by
Schmerl, see \cite{Sch85}.
\bn
\margintag{5.2}\ub{\stag{5.2} Question}:     (Maybe under $\bold V = \bold L$.) Suppose that
$\lambda^{\beth_\omega(\kappa)}=\lambda$ and
$\lambda_1^{<\lambda_1}=\lambda_1>\kappa_1$.   \ub{Then}
$(\lambda^+,\lambda,\kappa) \rightarrow(\lambda_1^+,\lambda_1,\kappa_1)$.
\bn
\margintag{5.3}\ub{\stag{5.3} Question}:   $(GCH)$ If $\lambda$ and $\mu$ are strong limit singulars
and $\lambda$ is a limit of supercompacts, \ub{then}
$(\lambda^+,\lambda) \rightarrow (\mu^+,\mu)$.
\bn
\margintag{5.4}\ub{\stag{5.4} Question}:   Find a universe with
$(\beth_2(\aleph_0),\aleph_0) \rightarrow (2^{2^\lambda},\lambda)$
for every $\lambda$. \nl
(The author has a written sketch of a result
which is close to this one. He starts with
$\aleph_0=\kappa_0<\kappa_1<\ldots<\kappa_m$  which are 
supercompacts and let $\Bbb P_n$ be the forcing which adds
$\kappa_{n+1}$ Cohen subsets to $\kappa_n$ in
$V^{{\Bbb P}_0\ast {\Bbb P}_1\ldots {\Bbb P}_{n-1}}$ for $n<m$.  The
idea is using the partition on trees from \cite[\S4]{Sh:288}).
\bn
\margintag{5.5}\ub{\stag{5.5} Question}:     Are all pairs in the set
$$
\{(\lambda,\mu):\,2^\mu = \mu^+ \and \mu = \mu^{<\mu}
\and \lambda \le 2^{\mu^+}\}
$$
such that there is $\lambda^+$-tree with $\ge\mu$ branches,
equivalent for the two cardinal problem?
\bn
More related to this particular work are \nl
\margintag{5.6}\ub{\stag{5.6} Question}:
\sn
1) Can we find $n < \omega$ and an infinite set $\Gamma^*$ of identities (or
2-identities) such that for any $\Gamma \subseteq \Gamma^*$ for some
forcing notion $\Bbb P$ in $\bold V^{\Bbb P}$ we have $\Gamma =
\Gamma^* \cap ID(\aleph_n,\aleph_0)$. \nl
2) In (1) we can consider $(\lambda,\mu)$ with $\mu =
\mu^{\aleph_0},\lambda = \mu^{+n}$, so we ask: can we find a forcing
notion $\Bbb P$ not adding reals such that for every $\Gamma \subseteq
\Gamma^*$ for some $\mu = \mu^{< \mu}$ we have $\Gamma = \Gamma^* \cap
ID(\mu^{+n},\mu)$. 
\bn
\margintag{5.7}\ub{\stag{5.7} Question}:  1) Can we get results 
parallel to \scite{3.7} for $(\aleph_2,\aleph_1)
+ 2^{\aleph_0} \ge \aleph_2$ (so we should start with a large
cardinal, at least a Mahlo). \nl
2) The parallel to \scite{5.6}(1),(2). 
\bn
\margintag{5.8}\ub{\stag{5.8} Question}:  1) Can we get results parallel to \scite{3.7} for
$(\aleph_{\omega +1},\aleph_\omega) +$ G.C.H. (or $(\mu^+,\mu),\mu$ strong limit
singular + G.C.H. \nl
2) The parallel to \scite{5.6}(1),(2). 
\bn
\margintag{5.9}\ub{\stag{5.9} Question}:   How does assuming MA $+ 2^{\aleph_0} >
\aleph_n$ influence $ID(\aleph_n,\aleph_0)$? (see below).
\bn
\ub{We end with some comments}:
\definition{\stag{t.1} Definition}  1) For $k \le \aleph_0$, we 
say $(\lambda,\mu)$ has $k$-simple identities \ub{when} $(a,e) \subseteq
ID(\lambda,\mu) \Rightarrow (a,e') \in ID(\lambda,\mu)$ whenever:
\mr
\item "{$(*)_k$}"  $a \subseteq \omega,(a,e)$ is an identity of
$(\lambda,\mu)$ and $e'$ is defined by
$$
b e'c \text{ iff } |b| = |c| \and (\forall b',c')[b' \subseteq b \and
|b'| \le k \and c' = OP_{c,b}(b') \rightarrow b' ec'].
$$
recalling $OP_{A,B}(\alpha) = \beta$ iff $\alpha \in A \and \beta \in
B \and \text{ otp}(\alpha \cap A) = \text{ otp}(\beta \cap B)$.
\endroster
\enddefinition
\bigskip

\proclaim{\stag{t.1a} Claim}  1) If $(\lambda_1,\mu_1)$ has $k$-simple
identities and there is $f:[\lambda_2]^{\le k} \rightarrow \mu_2$ such
that $ID_{\le k}(f) \subseteq ID_{\le k}(\lambda_1,\mu_1)$, \ub{then}
$(\lambda_1,\mu_1) \rightarrow (\lambda'_1,\mu'_1)$. \nl
2) If cf$(\lambda_1) > \mu$, then we can 
use $f$ with domain $[\lambda'_1]^{\le k} \backslash
[\lambda'_1]^{\le 1}$.
\endproclaim
\bigskip

\demo{Proof}  Should be easy.
\enddemo
\bigskip

\proclaim{\stag{t.2} Claim}  1) [MA $+ 2^{\aleph_0} > \aleph_n$].  The
\footnote{of course the needed version of MA is quite weak; going more
deeply in \cite{Sh:522}}
pair $(\aleph_n,\aleph_0)$ has 2-simple identities. \nl
2) If $\mu = \mu^{< \mu}$ and $\gamma \le \omega$ \ub{then} for some
$\mu^+$-c.c., $(< \mu)$-complete forcing notions, $\Bbb P$ in $\bold
V^{\Bbb P}$ we have $2^\mu \ge \mu^{+ \gamma}$ and $n \le \gamma \and
n < \omega \Rightarrow (\mu^{+n},\mu)$ has 2-simple identities. \nl
3) If $m < n < \omega,\mu = \mu^{< \mu}$, \ub{then}
$[\mu^{+n},\mu^{+m})$ has $(m+2)$-simple identities in $\bold V^{\Bbb
P}$ for appropriate $\mu^+$-c.c. $(< \mu)$-complete forcing notion. 
\endproclaim
\bigskip

\demo{Proof}  1) For any $c:[\aleph_n]^{< \aleph_0} \rightarrow \omega$ we
define a forcing notion $\Bbb P = \Bbb P_c$ as follows: \nl
$p \in \Bbb P$ \ub{iff}:
\mr
\item "{$(a)$}"  $p = (u,f) = (u^p,f^p)$
\sn
\item "{$(b)$}"  $u$ is a finite subset of $\aleph_n$
\sn
\item "{$(c)$}"   $f$ is a function from $[u]^2$ to $\omega$
\sn
\item "{$(d)$}"   if $k < \omega,k \ge 2$ and $\alpha_0 < \ldots < \alpha_{k-1}$ are from
$u,\beta_0 < \ldots < \beta_{k-1}$ are from $u$ and $[\ell(1) < \ell(2) < k
\Rightarrow f(\{\alpha_{\ell(1)},\alpha_{\ell(2)}\}) =
f(\{\beta_{\ell(1),\beta_{\ell(2)}}\})]$, \ub{then}
$c(\{\alpha_0,\dotsc,\alpha_{k-1}\}) =
c(\{\beta_0,\dotsc,\beta_{k-1}\})$.
\ermn
The rest should be clear. \nl
2), 3)  Similar (use e.g. \cite{Sh:546}).  \hfill$\square_{\scite{t.2}}$\margincite{t.2}
\enddemo
\bn
We can give an alternative proof of \cite{Sh:49}, note that by
absoluteness the assumption MA is not a real one; it can be eliminated
and $(\mu^{++},\mu) \rightarrow' (2^{\aleph_0},\aleph_0)$ can be deduced.
\proclaim{\stag{t.3} Claim}  Assume $MA + 2^{\aleph_0} >
\aleph_\omega$.

\ub{Then} $(\aleph_\omega,\aleph_0) \rightarrow (2^{\aleph_0},\aleph_0)$.
\endproclaim
\bigskip

\demo{Proof}  Let $\{\eta_\alpha:\alpha < 2^{\aleph_0}\}$ list
${}^\omega 2$, and define $f:[2^{\aleph_0}]^2 \rightarrow {}^{\omega
>}2$ by:
\mr
\item "{$(*)$}"  $f \{\alpha_0,\alpha_1\} = \eta_{\alpha_0} \cap
\eta_{\alpha_1} \in {}^{\omega >}2$ for $\alpha_0 \ne \alpha_1$.
\ermn
So by \scite{t.1a}, \scite{t.2} it is enough to prove that $ID_2(f)
\subseteq ID_2(\aleph_\omega,\aleph_0)$. \nl
Clearly
\mr
\item "{$(*)_1$}"  if $\lambda \le 2^{\aleph_0},(u,e) \in
ID_2(\lambda,\aleph_0)$ then $(u,e) \in ID_2(f \restriction
\lambda)$ hence $(u,e) \in ID_2(f)$ hence for some $n,(u,e)$ can be
embedded (in the natural sense) into
$({}^n 2,e^*_n)$ where $(\{\eta_1,\eta_2\} e^*_n\{\nu_1,\nu_2\})
\equiv (\eta_1 \cap \eta_2 = \nu_1 \cap \nu_2)$. 
\ermn
So it is enough to prove
\mr
\item "{$(*)_2$}"   $({}^n 2,e^*_n) \in ID_2(\mu^{+n},\mu)$.
\ermn
We prove this by induction on $n$.
\enddemo
\bn
\ub{$n=0$}:  Trivial.
\bn
\ub{$n+1$}:  Let $c:[\mu^{+n+1}]^2 \rightarrow \mu$, choose $M \prec
({\Cal H}(\mu^{+n+2}),\in)$ of cardinality $\mu^{+n}$ such that $\mu^{+n}+1
\subseteq M,c \in M$, so let $\delta = M \cap \mu^{+n}$. \nl
Define $c_n:\mu^{+n} \rightarrow \mu$ by $c_n \{\alpha,\beta\} =
(c\{\alpha,\beta\},c\{\delta,\alpha\},c\{\delta,\beta\})$ for $\alpha
< \beta < \mu^{+n}$.  By the induction hypothesis there is a sequence
$\langle \beta_\eta:\eta \in {}^n 2 \rangle$ of distinct ordinal $<
\mu^{+n}$ such that $\{\eta_1,\eta_2\} \,e^*_n\, \{\nu_1,\nu_2\}
\Rightarrow c_n\{\beta_{\eta_1},\beta_{\eta_2}\} =
c_n\{\beta_{\nu_1},\beta_{\nu_2}\}$. \nl
Let

$$
\align
A = \bigl\{ \gamma < \mu^{+n+1}:&\gamma \notin \{\beta_\eta:\eta \in {}^n 2\}
\text{ and for every} \\
  &\eta \in {}^n 2 \text{ we have } c\{\beta_\eta,\gamma\} =
c\{\beta_\eta,\delta\} \bigr\}.
\endalign
$$
\mn
Clearly $A \in M$ and $\delta \in A$ so $A \nsubseteq M$, 
hence necessarily $|A| =
\mu^{+n+1}$.  So by the induction hypothesis we can find a sequence
$\langle \gamma_\eta:\eta \in {}^n 2 \rangle$ of distinct members of
$A$ such that

$$
\{\eta_1,\eta_2\} e^*_n\{\nu_1,\nu_2\} \Rightarrow
c\{\gamma_{\eta_1,\eta_2},\gamma_{\eta_2}\} =
c\{\gamma_{\nu_1},\gamma_{\nu_2}\}.
$$
\mn
Now we define

$$
\alpha_\eta = \cases \beta_{\langle \eta(1 + \ell):\ell < n \rangle}
&\text{ \ub{if} } \eta(0)=0 \\
\gamma_{\langle \eta(1 + \ell):\ell < n \rangle}
&\text{ \ub{if} } \eta(0)=1 \endcases
$$
\mn
It is easy to check that $\langle \alpha_\eta:\eta \in {}^{n+1}2
\rangle$ is as required.  \hfill$\square_{\scite{t.3}}$\margincite{t.3}
\bn
We further can ask: \nl
\margintag{t.5}\ub{\stag{t.5} Question}:  Assume $\Gamma_i \subseteq ID^*$ for $i <
i^*,\Bbb P$ is $\Pi\{\Bbb P^\lambda_{\Gamma_i}:i < i^*\}$ with finite
support, $c:[\aleph_{n(*)}]^2 \rightarrow \omega$ in $\bold V^{\Bbb
P}$ \ub{then} $ID(c)$ is not too far from some $\dbcu_{i \in w}
\Gamma_i,w \subseteq i^*$ finite.
\bn
\margintag{t.7}\ub{\stag{t.7} Discussion}:  We can look more at ordered identities (recall)
\mr
\item "{$(*)_1$}"  for $\bold c_i:[\lambda]^{< \aleph_0} < \mu$ let
OID$(c) = \{(a,e):a$ a set of ordinals and there is an ordered
preserving $f:a \rightarrow \lambda\}$ such that $b_1 e b_2
\Rightarrow \bold c\{f''(b_1)) = \bold c(f''(b_2))$ and
OID$(\lambda,\mu) = \{(n,e):(n,e) \in \text{ OID}(\bold c)$ for every
$\bold c:[\lambda]^{< \aleph_0} \rightarrow \mu$, and similarly
OID$_2$, OID$_k$.
\ermn
Of course,
\mr
\item "{$(*)_2$}"  ID$(\lambda,\mu)$ can be computed from
OID$(\lambda,\mu)$.
\endroster
\bigskip
\newpage
\noindent

     \shlhetal 

\newpage
    
REFERENCES.  
\bibliographystyle{lit-plain}
\bibliography{lista,listb,listx,listf,liste}

\enddocument

\enddocument